\documentclass[a4paper,11pt]{article}
\usepackage{stmaryrd}
\usepackage{amsfonts}
\usepackage{bbm}
\usepackage{float}
\usepackage{amscd}
\usepackage{mathrsfs}
\usepackage{latexsym,amssymb,amsmath,amscd,amscd,amsthm,amsxtra,xypic}
\usepackage[dvips]{graphicx}
\usepackage[utf8]{inputenc}
\usepackage[T1]{fontenc}
\usepackage{lmodern}
\usepackage{amssymb}
\usepackage[all]{xy}
\usepackage{nicefrac,mathtools,enumitem}
\usepackage{microtype}
\usepackage{lipsum}
\usepackage{mathtools}
\usepackage{xfrac}

\usepackage{geometry}
\usepackage{hyperref}

\usepackage{times}
\usepackage{amsthm}
\usepackage{amsmath}
\usepackage{amsfonts}
\usepackage{amssymb}
\usepackage{mathdots}
\usepackage{color}
\usepackage{mathrsfs}
\usepackage{stmaryrd}
\SetSymbolFont{stmry}{bold}{U}{stmry}{m}{n}

\usepackage{bm}
\usepackage{tikz-cd}
\usepackage{tikz}
\usetikzlibrary{matrix,arrows,decorations.pathmorphing}
\usepackage[all]{xy}
\usepackage{graphicx}
\usepackage{subfigure}
\usepackage[toc,page]{appendix}
\usepackage[usestackEOL]{stackengine}

\usepackage{theoremref}
\newtheorem{thm}{Theorem}[section]
\newtheorem{lem}[thm]{Lemma}
\newtheorem{cor}[thm]{Corollary}
\newtheorem{pro}[thm]{Proposition}

\newtheorem{conj}[thm]{Conjecture} 

\theoremstyle{definition}

\newtheorem{rmk}[thm]{Remark}
\newtheorem{defi}[thm]{Definition}

\newtheorem{assumption}{Assumption}

\newcommand{\be }{\begin{equation}}
\newcommand{\ee }{\end{equation}}

\newcommand{\mathe}{\mathrm{e}}

\newcommand{\mathi}{\mathrm{i}}

\newcommand{\h}{\mathfrak h}

\def\qed{\hfill ~\vrule height6pt width6pt depth0pt}

\newcommand{\br}[1]{   [ \cdot,    \cdot  ]   }

\newcommand{\g}{\mathfrak g}

\newcommand{\gl}{\mathfrak {gl}}

\newcommand {\IR}{\mathbb{R}}

\newcommand {\Sect}{{\rm Sect}}

\newcommand {\IC}{\mathbb{C}}

\title{Representations of quantum groups arising from the Stokes phenomenon}
\author{Xiaomeng Xu}
\date{}

\newcommand{\Addresses}{{
  \bigskip
  \footnotesize
\noindent \textsc{
School of Mathematical Sciences \& Beijing International Center
for Mathematical Research, Peking University, Beijing 100871, China}\par\nopagebreak
  \textit{E-mail address}: \texttt{xxu@bicmr.pku.edu.cn}
}}
\begin{document}

\maketitle
    
\begin{abstract}
In this paper we prove that the quantum Stokes matrices of the quantum differential equation at a second order pole give rise to representations of the quantum group $U_q(\frak{gl}_n)$. We explain our results from the viewpoint of deformation quantization of the classical Stokes matrices at a second order pole. As a consequence, we can get a dictionary between the theory of Stokes phenomenon and the theory of quantum groups. We briefly discuss several such correspondences, and outline the generalization of our results to all classical types of Lie algebras and to the quantum differential equation at an arbitrary order pole.
\end{abstract}

\section{Introduction}
The Stokes phenomenon states that a solution of a meromorphic linear ordinary differential equation may have different asymptotic expressions as $z$ approaches an irregular singularity $z_0$ from different sectorial regions around $z_0$. See the books \cite{Balser, LR, Wasow} for a detailed introduction to the Stokes phenomenon. The discontinuous jump of the asymtotics around $z_0$ can be measured by the Stokes matrices. In the past decades, the Stokes phenomenon/matrices of meromorphic linear systems of differential equation at a second order pole have played an important role in many subjects of mathematics and physics. However the Stokes matrices are highly transcendental and thus hard to study. 
In this paper, we find the algebraic structures hidden behind the Stokes phenomenon. The importance of the result is that it enables us to apply the algebraic methods/ideas to solve some open analysis problems for the linear systems and the associated isomonodromy deformation equations. In particular, in later works, we solve the long time asymptotics and the associated Riemann-Hilbert problem \cite{Xu, TangXu}, the nonlinear connection problems \cite{Xu2}, and give the characterization of the WKB analysis \cite{ANXZ, Xu}, the classification of the soliton solutions and so on. In the following, let us explain the main result of this paper.

Let us take the Lie algebra ${\frak {gl}_n}$ over the field of complex numbers, and its universal enveloping algebra $U({\frak {gl}}_n)$ generated by $\{e_{ij}\}_{1\le i,j\le n}$ subject to the relation $[e_{ij},e_{kl}]=\delta_{jk}e_{il}-\delta_{li} e_{kj}$. Let us take the $n\times n$ matrix $T=(T_{ij})$ with entries valued in $U({\frak {gl}}_n)$
\begin{equation*}
T_{ij}=e_{ij}, \ \ \ \ \ \text{for} \ 1\le i,j\le n.
\end{equation*}
Let $\h_{\rm reg}$ denote the set of $n\times n$ diagonal matrices with distinct eigenvalues. Given any finite-dimensional irreducible representation $L(\lambda)$ of ${\gl}_n$ with highest weight $\lambda$, let us consider the quantum confluent hypergeometric system
\begin{equation}\label{introqeq}
\frac{dF}{dz}=h\Big( u+\frac{T}{z}\Big)\cdot F,
\end{equation}
for $F(z)\in {\rm End}(L(\lambda))\otimes {\rm End}(\mathbb{C}^n)$ an $n\times n$ matrix function with entries in ${\rm End}(L(\lambda))$. Here $h$ is a complex parameter, $u\in\h_{\rm reg}$ is seen as an $n\times n$ matrix with scalar entries in $U({\frak {gl}}_n)$, and the action of the coefficient matrix on $F(z)$ is given by matrix multiplication and the representation of ${\frak {gl}}_n$. 

\subsection{Main results}
\subsubsection*{Nonresonant case $h\notin \mathbb{Q}$}
Let us first assume $h\notin \mathbb{Q}$. The equation \eqref{introqeq} is then nonresonant and thus has a unique formal solution $\hat{F}$ around
$z = \infty$. See Proposition \ref{uniformal} for a proof. The standard theory of resummation (see e.g, \cite{Balser, LR, Wasow}) states that there exist certain  sectorial regions around $z=\infty$, such that on each of these sectors there is a unique (therefore canonical) holomorphic solution with the prescribed asymptotics $\hat{F}$. These solutions are in general different (that reflects the Stokes phenomenon), and the transition between them can be measured by a pair of Stokes matrices $S_{h\pm}(u)\in{\rm End}(L(\lambda))\otimes{\rm End}(\mathbb{C}^n)$. See Section \ref{beginsec} for more details. 

\begin{thm}(c.f. Theorem \ref{thm1})\label{introthm1}
For any fixed $h\notin\mathbb{Q}$ and $u\in\h_{\rm reg}$, the map (with $q=e^{\pi\mathi h}$)
\begin{equation}\label{mapbasis}
\begin{split}
\mathcal{S}_q(u): U_q(\frak{gl}_n)&\rightarrow {\rm End}(L(\lambda))~;\\
e_i&\mapsto \frac{S_{h+}(u)_{i,i}^{-1}\cdot S_{h+}(u)_{i,i+1}}{q^{-1}-q}, \\
f_i&\mapsto \frac{S_{h-}(u)_{i+1,i}\cdot S_{h-}(u)_{i,i}^{-1}}{q-q^{-1}}, \\
q^{h_i}&\mapsto S_{h+}(u)_{i,i}
\end{split}
\end{equation}
defines a representation of the Drinfeld-Jimbo quantum group $U_q(\frak{gl}_n)$ on the vector space $L(\lambda)$. Here recall that $U_q(\frak{gl}_n)$ is a unital associative algebra with generators $q^{\pm h_i}, e_j, f_j,$ $1\le j\le n-1, 1\le i\le n$ and relations:
\begin{itemize}
    \item for each $1\le i\le n$, $1\le j\le n-1$,
\[q^{h_i}q^{-h_i}=q^{-h_i}q^{h_i}=1, \
q^{h_i}e_jq^{-h_i}=q^{\delta_{ij}}q^{-\delta_{i,j+1}}e_j, \ q^{h_i}f_jq^{-h_i}=q^{-\delta_{ij}}q^{\delta_{i,j+1}}f_j;
\]
\item for each $1\le i,j\le n-1$,
\[
[e_i,f_j] = \delta_{ij} \frac{q^{h_i-h_{i+1}}-q^{-h_i+h_{i+1}}}{q-q^{-1}};
\]
\item for $|i-j|=1$, 
\[
e_i^2e_j - (q+q^{-1})e_ie_je_i + e_je_i^2=0, 
\]
\[
f_i^2f_j - (q+q^{-1})f_if_jf_i + f_jf_i^2=0,
\]
and for $|i-j|\ne 1$, $[e_i,e_j]=0=[f_i,f_j]$.
\end{itemize}
\end{thm}

The theorem associates to any representation $L(\lambda)$ of $U(\frak{gl}_n)$ a representation $\mathcal{S}_q(u)$ of $U_q(\frak{gl}_n)$ on the same vector space $L(\lambda)$. 

\subsubsection*{Resonant case $h\in\mathbb{Q}$}
It is convenient to think of \eqref{introqeq} as a blocked system with each block a ${\rm dim}(L(\lambda))\times {\rm dim}(L(\lambda))$ matrix.
Its irregular term $u$ has $n$ diagonal blocks, each of them has repeated eigenvalues (each $u_i$ repeat ${\rm dim}(L(\lambda))$ times). Now let us assume $h=h_0\in\mathbb{Q}$. If for some $k=1,...,n$ the corresponding blocks $h_0e_{kk}$ in the sub-leading term $T$ have two eigenvalues differed by a positive integer, the corresponding differential equation \eqref{introqeq} becomes resonant. In this case, the uniqueness of the formal fundamental solution $\hat{F}(z)$ is not valid. Instead, the equation \eqref{introqeq} has a family of formal solution around
$z = \infty$ depending on a finite set $c$ of complex parameters. Accordingly, there are a family of Stokes matrices $S_{h_0\pm}(u;c)$ depending on the same set of parameters. 
However, there is a pair of distinguished Stokes matrices, denoted by $S_{h_0\pm}^o(u)$, among the family $S_{h_0\pm}(u;c)$. They are actually the continuous extension of $S_{h\pm}(u)$ from $h\in\mathbb{C}\setminus \mathbb{Q}$ to $h_0\in\mathbb{Q}$. See Section \ref{seclast} for more details, and Section \ref{2by2} for an explicit example.

\begin{thm}\label{introthm1re}
For any fixed $h_0\in\mathbb{Q}$ and $u\in\h_{\rm reg}$, the map $\mathcal{S}_{q_0}(u)$ given in \eqref{mapbasis} defines a representation of $U_{q_0}(\frak{gl}_n)$ at $q_0=e^{\pi\mathi h_0}$ a root of unity.
\end{thm}

By pursuing the heuristics of the above algebraic structure of the Stokes matrices, we are led to rather interesting results. On the one hand, Theorem \ref{introthm1} indicates a characterization of the Stokes phenomenon in the WKB approximation via the theory of canonical basis introduced by Lusztig \cite{Lu} and Kashiwara \cite{Ka1, Ka2}, see \cite{Xu} or Section \ref{app} for more details. On the other hand, Theorem \ref{introthm1re} indicates an algebraic characterization of the soliton/polynomial solutions of confluent hypergeometric type equations. These seem to us striking examples of the applications of the algebraic understanding of the Stokes matrices to analysis problems.

\subsection{Motivation: deformation quantization of the classical Stokes matrices at second order poles}\label{app}
In this subsection, we explain our original motivation and explain the main theorem from the viewpoint of deformation quantization. 

Let us consider the classical confluent hypergeometric system for a function $f(z)\in {\rm GL}_n$
\begin{equation}\label{introeq}
\frac{df}{dz} = \left(u+\frac{A}{z}\right)\cdot f,\end{equation}
where $u\in \h_{\rm reg}$, and $A=(a_{ij})\in\frak{gl}_n$. 

In \cite{Boalch1}, Boalch interpreted the dual Poisson Lie group $GL_n^*$ of $GL_n$ as the space of Stokes matrices $S_\pm(u,A)\in GL_n$ of \eqref{introeq}, and proved that the irregular Riemann-Hilbert map (also known as Riemann-Hilbert-Birkhoff map) from the de Rham to the Betti moduli space
\begin{equation}\label{RHB}
    \nu(u):\frak{gl}_n^*\cong \frak{gl}_n\rightarrow GL_n^*;~A\mapsto (S_-(u,A), S_+(u,A))
\end{equation}
associating the Stokes matrices to the confluent hypergeometric system \eqref{introeq} parameterized by $A\in\frak{gl}_n$, is a Poisson map. 

In terms of the r-matrix formulation of the Poisson brackets on the dual Poisson Lie group, Boalch's theorem can be rewritten via the classical RLL formulation (see {\cite[Formula (235)]{AM1}}):
\begin{thm}\cite{Boalch1}\label{sclrll}
For any fixed $u\in\h_{\rm reg}$, the	matrix valued function $M(u,A):=S_-(u,A) S_+(u,A)\in GL_n$ of $A\in \frak{gl}_n$ satisfies
	\begin{eqnarray}\label{StokesPoisson} 
		\{M^1 \otimes M^2 \} = r_+ M^1 M^2 + M^1 M^2 r_- - M^1 r_+M^2 -M^2 r_-M^1
		.  \label{cRLL}
	\end{eqnarray}
	Here
	\begin{align*}
		r_+ & =  \frac{1}{2}  \sum_{i = 1}^n E_{ii} \otimes E_{ii} + \sum_{1 \le
			i < j \le n} E_{ij} \otimes E_{ji} \in \mathrm{End} (\mathbb{C}^n) \otimes
		\mathrm{End} (\mathbb{C}^n), \\
		r_- & =  - \frac{1}{2}  \sum_{i = 1}^n E_{ii} \otimes E_{ii} - \sum_{1
			\le j < i \le n} E_{ij} \otimes E_{ji} \in \mathrm{End} (\mathbb{C}^n)
		\otimes \mathrm{End} (\mathbb{C}^n), 
	\end{align*}
	and the tensor notation is used: both sides are $\mathrm{End} (\mathbb{C}^n)
	\otimes \mathrm{End} (\mathbb{C}^n)$ valued functions on
	$A\in\mathfrak{g}\mathfrak{l}_n$, $M^1 := M(u) \otimes 1$, $\nu^2 := 1
	\otimes \nu(u)$, and the $ij, kl$ coefficient of the matrix $\{M^1,M^2 \}$ is
	defined as the Poisson bracket $\{M(u)_{ij}, M(u)_{kl} \}$ of the functions on the
	canonical linear Poisson space $\mathfrak{g}\mathfrak{l}_n^{\ast} \cong
	\mathfrak{g}\mathfrak{l}_n$.
\end{thm}
The theorem gives rise to a distinguished (complex version of the) Ginzburg-Weinstein linearization of the dual Poisson Lie group \cite{GW}. In the literature, there are many proofs of the Ginzburg-Weinstein linearization theorem, from the perspectives of cohomology calculation, Moser's trick in symplectic geometry, Stokes phenomenon, Gelfand-Tsetlin integrable systems, the quantum algebras and so on, see e.g., \cite{GW, Anton,AM, Boalch1,Boalch2,EEM}.   
It is natural to ask if these seemingly rather different methods are related to each other in some deep manner. We have been working on this problem with the belief that the connection between the different methods will bring new insights into the above mentioned subjects. In a series of works, we clarified the relations between these constructions. In particular, in \cite{Xu0}, we pointed out a relation between the Boalch's construction and the Enriquez-Etingof-Marshall construction. It motivates us to study the deformation of the classical Stokes matrices. In the following, let us explain the idea.

On the one hand, it is known that (see e.g., Proposition \ref{sclrmatrix}) the classical equation \eqref{introeq} has a unique formal solution
\begin{eqnarray}\label{sconeformal}
\widehat{f}(z)=\widehat{h}(z) e^{{uz}}z^{\delta A}, \ \ \ \text{with} \ \widehat{h}=1+h_1z^{-1}+h_2z^{-1}+\cdot\cdot\cdot, \end{eqnarray}
where $\delta A$ is the diagonal part of the matrix $A=(a_{ij})_{i,j=1,...,n}$, and the $(k,l)$-entry $h_{m,kl}\in {\mathbb{C}}$ of the $m$-th coefficient $h_{m}\in {\rm End}(\mathbb{C}^n)$ is given by a degree $2m$ polynomial in the variables $\{a_{ij}\}$. The formal power series $\hat{h}(z)$ itself already encodes all the information. Thus the observation is that the Stokes phenomenon of \eqref{introeq} for all $A\in \frak{gl}_n$ is encoded in the sequence of functions
\begin{equation} \label{scldata}\Big\{h_{m}(u,A)\in Sym(\frak{gl}_n)\otimes {\rm End}(\mathbb{C}^n) \Big\}.
\end{equation}

On the other hand, following Proposition \ref{Yangsol}, the equation \eqref{introqeq} has a formal solution
\begin{eqnarray}\label{oneformal}
\widehat{F}(z)=\widehat{H}(z) e^{{huz}}z^{h\delta T}, \ \ \ \text{with} \ \widehat{H}=1+H_1z^{-1}+H_2z^{-1}+\cdot\cdot\cdot, \end{eqnarray}
where $\delta T$ is the diagonal part of $T$, and the $(k,l)$-entry $H_{m,kl}\in {\rm End}(L(\lambda))$ of the $m$-th coefficient $H_{m}$ is given by (the evaluation on the representation $L(\lambda)$ of) a filtered degree $2m$ element in $U(\frak{gl}_n)$. Thus the Stokes phenomenon of \eqref{introqeq} is encoded in the sequence 
\[\Big\{H_{m}(hu,hT)\in U(\frak{gl}_n)\otimes {\rm End}(\mathbb{C}^n) \Big\}.\]

Following Proposition \ref{Yangsol} and \ref{sclrmatrix}, the quantum coefficients $H_m$ satisfy the same recursive relation as the classical coefficients $h_m$ (provided replacing $A$ and $\delta A$ by $T$ and $\delta T$ for all $l=1,...,k$).
Thus, the quantum equation \eqref{introqeq} provides us a canonical deformation quantization of the sequence \eqref{scldata} from the Poisson algebra $Sym(\frak{gl}_n)$ to the quantized algebra $U(\frak{gl}_n)$ (see Section \ref{seclast}): for all positive integer $m$ and $k,l=1,...,n$, we have a canonical lift
\begin{equation}\label{lift}
    \Big\{h_{m,kl}(u,A)\in Sym(\frak{gl}_n)\Big\}\rightarrow \Big\{H_{m,kl}(hu,hT)\in U(\frak{gl}_n)\Big\}.
\end{equation}
The Stokes matrices $S_\pm(u,A)$ and $S_{h\pm}(u)$ of \eqref{introeq} and \eqref{introqeq} are obtained via the Borel-Laplace transforms of the power series $\sum_m h_mz^{-m}$ and $\sum_m H_mz^{-m}$ respectively. And the Borel-Laplace transforms are compatible with the deformation quantization in \eqref{lift}. Therefore, after taking
completion, the Stokes matrices $S_{h\pm}(u)$ of \eqref{hqeq}, seen as $n\times n$ matrices with entries valued in $\hat{U}(\frak{gl}_n)\llbracket h\rrbracket$ (or in any ${\rm End}(L(\lambda))$) are understood as the deformation quantization of the classical Stokes matrices $S_\pm(u,A)\in \widehat{Sym}(\frak{gl}_n)\otimes {\rm End}(\mathbb{C}^n)$ of \eqref{introeq}, seen as $n\times n$ matrix valued analytic functions of the parameter $A\in {\frak {gl}}_n^*$. 
\begin{rmk}
More directly, using the associated function method, see e.g.,\cite[Chapter 9]{Balser} and \cite{LuS,Yokoyama1988}, the Stokes matrices can be obtained from the regularized limit of the matrix $h_m$ and $H_m$ as $m\rightarrow \infty$ respectively (as a consequence, the Poisson brackets between the entries of Stokes matrices $S_\pm(u,A)$, as functions on $A\in\frak{gl}_n^*$, can be obtained by taking the limit of the Poisson brackets between the entries of the coefficient matrices $h_m(u,A)$ as $m\rightarrow \infty$, see \cite{TXu}).
\end{rmk} 
That is considering the Stokes matrices of \eqref{introqeq} leads to a canonical lift for all $i,j=1,...,n$
\begin{equation}\label{liftS}
    \Big\{S_\pm(u,A)_{ij} \in \widehat{Sym}(\frak{gl}_n) \Big\}\rightarrow \Big\{S_{h\pm}(u)_{ij} \in \widehat{U}(\frak{gl}_n)\llbracket h\rrbracket \Big\}.
\end{equation}
By this reason, we also call $S_{h\pm}(u)$ the quantum Stokes matrices.

Now the irregular Riemann-Hilbert map $\nu(u)$ in \eqref{RHB} pulls back the entry coordinates on the space of Stokes matrices to the highly transcendental functions $\{S_\pm(u,A)_{ij}\}$ on $A\in\frak{gl}_n^*$. This collection of functions is the data of the irregular Riemann-Hilbert map. By \eqref{liftS}, the (quantum) Stokes matrices of \eqref{introqeq} is seen as a quantization/deformation of the irregular Riemann-Hilbert map in \eqref{RHB}. 

For fixed $u$, this collection of functions $\{S_\pm(u,A)_{ij}\}$ is a new coordinate system around a small neighbourhood of $0\in \frak{gl}_n^*$ (here we take the upper and lower entries of the upper and lower triangular Stokes matrices respectively). Theorem \ref{sclrll} states that in terms of the new coordinates $\{S_\pm(u,A)_{ij}\}$, the canonical linear Poisson bracket on $\frak{gl}_n^*$ takes the nonlinear (but still algebraic) classical r-matrix form. That is the Poisson bracket is linear under the entry coordinates of $\frak{gl}_n^*$, and nonlinear but still has a nice form under the transcendental coordinates $\{S_\pm(u,A)_{ij}\}$. 

After the deformation quantization, the coordinates $\{S_\pm(u,A)_{ij}\}$ become elements $\{S_{h\pm}(u)_{ij}\}$ in $\hat{U}(\frak{gl}_n)\llbracket h\rrbracket$ or by taking a representation in ${\rm End}(L(\lambda))$. Now the question is if the associative algebra product in $U(\frak{gl}_n)$ or the ordinary composition in ${\rm End}(L(\lambda))$, between these elements $S_{h\pm}(u)_{ij}$, takes a nice form. Given Theorem \ref{sclrll}, we expect that the entries of (quantum) Stokes matrices $S_{h\pm}(u)$ satisfy the quantum RLL relation. In the paper we prove such a statement as in Theorem \ref{introthm1} (see Section \ref{Uqg} for the RLL relation and its equivalence with Theorem \ref{introthm1}).

The important idea conveyed in the discussion above is that, the deformation/quantization of the irregular Riemann-Hilbert map, thought of as a collection of transcendental functions on $\frak{gl}_n^*$, is just a collection of generators in the associative algebra or ${\rm End}(L(\lambda))$. In a sense, it already quantizes the Poisson structures involved. We stress that to see if this collection of elements satisfy natural algebraic relations (like the RLL relation) is a different question. 

Forgetting the transcendental Riemann-Hilbert-Birkhoff map, one can study the quantization of the algebraic Poisson structures on the moduli space in a purely algebraic way. The quantization of the moduli space of flat connections on a Riemann surface
has been extensively studied, mainly motivated by a relation to the Chern-Simons theory and knot invariants \cite{RT, Tu, Witten1, Witten2}. The quantization of Poisson structures on certain moduli spaces closely related to the dual Poisson Lie group (i.e., the space of Stokes matrices) was studied in \cite{CMR, CS}, with the help of the Fock-Goncharov coordinates \cite{FG0}. In particular, the monodromy matrices are explicitly expressed in terms of the Fock-Goncharov coordinates. Then the quantum monodromy matrices, defined by the same expression provided the Fock-Goncharov coordinates are replaced by their natural quantization and the proper ordering is accounted for (see \cite{CS, CGT, FG}), are proved to satisfy the quantum R-matrix type relations \cite{CMR,CS}. It is seen as a quantization of the Betti moduli space. To the best of our knowledge, the quantization of the Riemann-Hilbert(-Birkhoff) maps from the Betti to the de Rham moduli spaces, as a (highly transcendental) homomorphism between associative algebras, is missing in the literature. In particular, it is an interesting question to see if the quantum monodromy matrices in \cite{CMR, CS} have a transcendental construction, i.e., can be realized as the monodromy data of certain differential equations, and then to see if they are related to the Stokes matrices of the quantum confluent hypergeometric systems.

\subsection{The main application of Theorem \ref{introthm1}: crystal basis arising from the WKB approximation}
In this subsection, we explain how Theorem \ref{introthm1} and \ref{introthm1re} can provide a dictionary between the Stokes phenomenon and the theory of quantum groups, by discussing several striking examples of this kind. Some of them will appear in our follow up work \cite{TX, Xu}, and the others are left for future study. In particular, we will explain the part (4) in the following table. The correspondence (3), (5), (6) are given in \cite{Xu}.

\begin{table}[h]
\caption{A dictionary between the Stokes phenomenon at second order poles and the theory of quantum groups}
\begin{tabular}{ |c|c|c| } 
 \hline
 & Stokes phenomenon of equation \eqref{introqeq} & Quantum group $U_q(\frak{gl}_n)$ with $q=e^{\pi \mathi h}$\\ 
\hline
1 & Nonresonant case $h\notin \mathbb{Q}$ & Representation of $U_q(\frak{gl}_n)$ at a generic $q$\\ 
 \hline
2 & Resonant case $h\in \mathbb{Q}$ &  Representation at roots of unity \\  
 \hline
3 & Long time $u_1\ll \cdots\ll u_n$ behaviour & Braching rules/ Gelfand-Tsetlin theory\\
\hline
4 & WKB approximation as $h\rightarrow+\infty$ & $\frak{gl}_n$-crystals \\
 \hline
5 & Wall-Crossing formula in the WKB approximation & cactus group actions on $\frak{gl}_n$-crystals\\
 \hline
6 & Whitham dynamics & HKRW cover by eigenlines\\
 \hline
7 & Involution of equations & Quantum symmetric pairs
\\
 \hline
8 & Formal power series solutions & Yangians/ Trigonometric R-matrix\\
\hline 
\end{tabular}
\end{table}

On the one hand, the WKB method, named after Wentzel, Kramers, and Brillouin, is for approximating solutions of a differential equation whose highest derivative is multiplied by a small parameter (other
names, including Liouville, Green, and Jeffreys are sometimes attached to this method). 

On the other hand, a ${\gl}_n$-crystal is a finite set equipped with some operators called crystal operators satisfying certain conditions, where the finite set models a weight basis for a representation of ${\gl}_n$, and crystal operators indicate the leading order behaviour of the simple root vectors on the basis under the crystal limit $q\rightarrow \infty$ in the quantum group $U_q({\gl}_n)$. See e.g., \cite{Jo, HKbook} for more details.

Note that in the realization of the quantum $U_q(\frak{gl}_n)$ via the quantum Stokes matrices $S_{h\pm}(u)$ given in Theorem \ref{introthm1}, the parameters $q$ and $h$ are related by $q=e^{h/2}$. Therefore, the WKB limit as $h\rightarrow +\infty$ (along the positive real axis) should correspond to the crystal limit $q\rightarrow \infty$,
\[\fbox{WKB approximation of the equation \eqref{introqeq}} \longleftrightarrow \fbox{crystal limit of the quantum group $U_q({\gl}_n)$}\] 
The correspondence is made precise in our follow up work \cite{Xu}. Here we only outline the correspondence, and refer to \cite{Xu} for full details. 
Let us take the shift of argument subalgebra $\mathcal{A}(u)$ of $U({\gl}_n)$, which is a maximal commutative subalgebras parameterized by $u\in \h_{\rm reg}(\mathbb{C})$, see e.g., \cite{HKRW, Ry}. It was proved \cite{FFR, HKRW} that $\mathcal{A}(u)$ acts cyclically on the representation $L(\lambda)$ for any $u\in\widetilde{h_{\rm reg}}$. In particular, the action of $\mathcal{A}(u)$ on $L(\lambda)$ has simple spectral for the real $u\in\h_{\rm reg}(\mathbb{R})$ (the set of $n\times n$ diagonal matrices with distinct real eigenvalues). We denote by $B(u;\lambda)$ an orthonormal eigenbasis of the action of $\mathcal{A}(u)$ on $L(\lambda)$. 
Let us now make the following WKB asymptotic expansion assumption. 
\begin{assumption}\label{assumption} For $h\in\mathbb{R}_{>0}$ and $u\in\h_{\rm reg}(\mathbb{R})$, the action of the off-diagonal entry $S_{h+}(u)_{k,k+1}$ of the quantum Stokes matrix on the eigenbasis vectors $\{v_a(u)\}_{a\in I}\in {\rm End}(L(\lambda))$ of $\mathcal{A}(u)$ have the WKB type asymptotic behaviour
\begin{equation}\label{WKBasy}
S_{h+}(u)_{k,k+1} \cdot v_a(u)=\sum_{b\in I}e^{h\phi^{(k)}_{ab}(u)+\mathi h g_{ab}^{(k)}(u)} v_b(u)\left(1+\mathcal{O}(h^{-1})\right),
\end{equation}
as $h\rightarrow +\infty$ uniformly with respect to $u\in U_\sigma$,
where for all $a,b\in I$, $\phi^{(k)}_{ab}(u)$ and $g_{ab}^{(k)}(u)$ are some real valued functions of $u$. \end{assumption}
An element $v_a(u)$ of $B(u;\lambda)$ is called generic if there exists only one index $c\in I$ such that $\phi^{(k)}_{ac}(u)$ is the biggest in the collection $\{\phi^{(k)}_{ab}(u)\}_{b\in I}$ of real numbers.
Thus, the WKB approximation of $S_{h+}(u)_{k,k+1}$ naturally defines an operator $\widetilde{e_k}$ on the generic elements of $B(u;\lambda)$ by picking the unique leading term in \eqref{WKBasy}, i.e.,
\begin{eqnarray}
\widetilde{e_k}v_a(u):= v_c(u), \ \ \ \text{if} \ \ \ \phi^{(k)}_{ac}(u)={\rm max}\{\phi^{(k)}_{ab}(u)~|~{b\in I}\}.
\end{eqnarray}
Similarly, by considering the WKB approximation of $S_{h-}(u)_{k+1,k}$, one defines an operator $\widetilde{f_k}$ on (some other) generic elements of $B(u;\lambda)$. In a universal sense, the operators $\widetilde{e_k}$ and $\widetilde{f_k}$ extend from generic elements to all elements of $B(u;\lambda)$, see \cite{Xu} for more details.
In the heuristic spirit, the correspondence between the WKB approximation and the crystal limit predicts that the finite set $B(u;\lambda)$ equipped with the operators $\{\widetilde{e_k}(u), \widetilde{f_k}(u)\}_{k=1,...,n-1}$ is a $\frak{gl}_n$-crystal. It is proved in \cite{Xu} that
\begin{thm}\cite{Xu}\label{conj2}
Under Assumption \ref{assumption}, for any $u\in\h_{\rm reg}(\mathbb{R})$ and each $k=1,...,n-1$, there exists canonical operators $\widetilde{e_k}(u)$ and $\widetilde{f_k}(u)$ acting on the finite set $B(u;\lambda)$, and real valued functions $c_{kj}(\xi(u))$, $\theta_{kj}(h,u,\xi(u))$ with $j=1,2$ such that for any generic element $\xi(u)\in B(u;\lambda)$,
\begin{eqnarray} \label{eku}
\mathop{\rm lim}\limits_{h\rightarrow +\infty}\Big(S_{h+}(u)_{k,k+1}\cdot e^{c_{k1}(\xi) h+\mathi  \theta_{k1}(h,u,\xi)} \xi(u)\Big)=\widetilde{e_k}(\xi(u)),\\ \label{fku}
\mathop{\rm lim}\limits_{h\rightarrow +\infty}\Big(S_{h-}(u)_{k+1,k}\cdot e^{c_{k2}(\xi) h+\mathi  \theta_{k2}(h,u,\xi)} \xi(u)\Big)=\widetilde{f_k}(\xi(u)).
\end{eqnarray}
Furthermore, the WKB datum $(B(u;\lambda), \widetilde{e_k}(u), \widetilde{f_k}(u))$ is a $\frak{gl}_n-$crystal.
\end{thm}
Here we remark that the $\frak{gl}_n$-crystals are unique, see e.g., \cite[6.4.21]{Jo}, and there are a number of ways to construct them: combinatorially using Littelmann’s path
model \cite{Li}, representation theoretically using crystal bases of a quantum group representation \cite{Ka2},
and geometrically using the affine Grassmannian \cite{BG}. And Theorem \ref{conj2} gives a transcendental construction of them. 
More importantly, it characterizes the WKB approximation of the (quantum) Stokes matrices of \eqref{introqeq} via the $\frak{gl}_n$-crystals. In the classical case, one can characterize the WKB approximation of the (classical) Stokes matrices of \eqref{introeq} via closely related algebraic structures \cite{ANXZ}. We remark that the description of the WKB approximation of linear system of differential equations with rank higher than two is still an open problem. One of the finial goal of our project is to characterize the WKB approximation in the Stokes phenomenon (at an arbitrary order pole) via algebraic method.

\subsubsection*{The application of Theorem \ref{introthm1re}: soliton solutions of confluent hypergeometric equations}
Let us consider the case $q_0=e^{\pi\mathi h_0}$ a primitive $p$-th root of unity where $p$ is odd. 
The representation theory of the quantum group $U_{q_0}(\frak{gl}_n)$ is much richer than the generic case. We refer the reader to the book \cite{KS} for an introduction to the representation
theory in the roots of unity case. In particular, any irreducible representation of $U_{q_0}(\frak{gl}_n)$ has dimension less than or equal to $\frac{n(n-1)}{2}$, and an irreducible representation can be neither a highest nor a lowest weight representation. Now let us focus on the representation $\mathcal{S}_{q_0}(u)$ of $U_{q_0}(\frak{gl}_n)$ on the vector space $L(\lambda)$ obtained in Theorem \ref{introthm1re}.

On the one hand, following \cite{HKRW} (or from the viewpoint of the isomonodromy deformation in the WKB approximation, see \cite{Xu}), for $u\in\h_{\rm reg}(\mathbb{R})$, there exists a parametrization of the eigenbasis $B(u;\lambda)$, of the action of the shift of argument algebra $\mathcal{A}(u)$ on $L(\lambda)$, by the Gelfand-Tsetlin patterns $\Lambda$ (up to a cactus group action). Such a pattern $\Lambda$ (for fixed $(\lambda^{(n)}_1,...,\lambda^{(n)}_n):=\lambda$) is a collection of numbers $\{\lambda^{(i)}_j(\Lambda)\}_{1\le i\le j\le n-1}$ satisfying the interlacing conditions
\begin{equation}\label{interineq}
\lambda_j^{(i+1)}(\Lambda)-\lambda^{(i)}_j(\Lambda)\in\mathbb{Z}_{\ge 0}, \hspace{5mm} \lambda^{(i)}_j(\Lambda)-\lambda^{(i+1)}_{j+1}(\Lambda)\in\mathbb{Z}_{\ge 0}.
\end{equation} 
Let us denote by $\xi_{\Lambda}(u)\in B(u;\lambda)$ the basis vector labelled by the pattern $\Lambda$.

On the other hand, let us assume the highest weight $\lambda$ is such that $\lambda^{(n)}_{i-1}-\lambda^{(n)}_{i}$ is bigger than $p$ for $i=2,...,n$. Since the representation $\mathcal{S}_{q_0}(u)$ of $U_{q_0}(\frak{gl}_n)$ at the root of unity is the continuous extension of the irreducible representation $\mathcal{S}_q(u)$ as $q\rightarrow q_0$, one can argue that generically the actions of $S^o_{h_0\pm}(u)_{k,k+1}$ for all $k$ on $L(\lambda)$ generate invariant subspaces of dimension $p^{{n(n-1)}/{2}}$. It motivates
\begin{conj}\label{lastconj} At the primitive $p$-th root of unity, each invariant subspace of dimension $p^{{n(n-1)}/{2}}$ of the representation $\mathcal{S}_{q_0}(u)$ of $U_{q_0}(\frak{gl}_n)$ on $L(\lambda)$ is spanned by the basis vectors $\xi_{\Lambda}(u)$ whose corresponding patterns $\Lambda$ satisfy \[(m_{ij}-1)p\le \lambda^{(i+1)}_j(\Lambda)-\lambda^{(i)}_{j}(\Lambda)\le m_{ij}p-1\] for all $1\le j\le i\le n-1$ and some fixed positive integer $m_{ij}$. In particular, if $\Lambda^o$ is the pattern determined by
\begin{equation}
\lambda^{(i)}_j(\Lambda^o)-\lambda^{(i+1)}_{j+1}(\Lambda^o)=m_{ij}p-1, \hspace{5mm} \text{for } 1\le j\le i\le n-1.
\end{equation}
then the actions of $S^o_{h_0+}(u)_{k,k+1}$ on $\xi_{\Lambda^o}(u)$ are zero for all $k=1,...,n-1$. 
\end{conj}
It follows directly from the explicit expression of the Stokes matrices given in the formula \eqref{Stokes22} in Section \ref{2by2} that
\begin{pro}\label{prop2}
The conjecture is true for $n=2$.
\end{pro}
\begin{rmk}
In order to prove this conjecture, in \cite{TX} we put the study of the Stokes phenomenon of \eqref{introqeq} in the framework of the quantum inverse scattering method, and prove that its Stokes matrices can be obtained as the infinite products of representations of Yangians. Then the idea is to prove Conjecture \ref{lastconj} via the representation of Yangians, which we leave to a next paper. 
\end{rmk}
As far as the quantum–classical correspondence is concerned, 
Conjecture \ref{lastconj} can be used to characterize those $u$ and $A$ such that the corresponding classical equation \eqref{introeq} has a soliton solution taking the form
\[f(z)=(1+h_1 z^{-1}+\cdots + h_m z^{-m})\cdot e^{uz}z^{\delta A},\]
where $h_i$ is an $n\times n$ matrix, $m$ is a positive integer and $\delta A$ is the diagonal part of $A$. It is because the existence of soliton solutions amounts to determining when the Stokes matrices $S_\pm(u,A)$ are diagonal. Here the name soliton solution is after \cite[Section 6]{JMU}. We leave the study for a future work.

\subsection{Outlook}
\subsubsection*{Generalization to quantum groups of classical types}
Let us discuss the generalization of Theorem \ref{introthm1} to other types. Let $\g_n$ denote the rank $n$ simple complex Lie algebra of type $B, C$, or $D$, i.e.,
\[\g_n=\frak{o}_{2n+1}, \ \ \frak{sp}_{2n}, \ \ 
 \text{or} \ \ \frak{o}_{2n}.\]
For $-n\le i,j\le n$, let us introduce the generators of $U(\mathfrak{g}_n)$
\[K_{ij}=e_{ij}-\theta_{ij} e_{ji}, \]
where
\begin{equation}
\theta_{ij}=\left\{\begin{array}{lr} 1,  & \text{ in the orthogonal case} \\ {\rm sgn}(i)\cdot {\rm sgn}(j), & \text{ in the symplectic case}.
             \end{array}
\right.
\end{equation}

Given any finite-dimensional irreducible representation $L(\lambda)_{\g_n}$ of $\g_n$ with a highest weight $\lambda$, let us consider the quantum confluent hypergeometric type equation
\begin{eqnarray}\label{introqeqBCD}
\frac{dF}{dz}=h\Big( u+\frac{K}{z}\Big)\cdot F,
\end{eqnarray}
for $F(z)\in {\rm End}(L(\lambda)_{\g_n})\otimes {\rm End}(\mathbb{C}^{m})$ with $m=2n$ or $2n+1$. Here the $m\times m$ matrix $K=(K_{ij})$ has entries valued in $U(\g_n)$, and the $m\times m$ matrix $u$ is diagonal with entries 
\begin{equation}
u=\left\{\begin{array}{lr} {\rm diag}(- u_n, \ldots, -u_1, 0, u_1, \ldots, u_n),  & \mathfrak{s}\mathfrak{o}_{2 n + 1} \\ {\rm diag}(- u_n, \ldots, -u_1, u_1, \ldots, u_n),  & \mathfrak{s}\mathfrak{o}_{2 n} \text{ or } \mathfrak{s}\mathfrak{p}_{2 n}.
             \end{array}
\right.
\end{equation}

We make the assumption that $u$ is regular, i.e., $u_1,...,u_n$ and $0$ are distinct. Then the proof of Theorem \ref{thmRLL} can be applied to show that for any $h\notin\mathbb{Q}$ and regular $u$, the Stokes matrices $S_{h\pm}(u)_{\g_n}$ of \eqref{introqeqBCD} also satisfy the RLL relation
\begin{align*}
R_{\g_n}^{12}S^{(1)}_{\pm}(u)_{\g_n}S^{(2)}_{\pm}(u)_{\g_n} &=S^{(2)}_{\pm}(u)_{\g_n}S^{(1)}_{\pm}(u)_{\g_n} R_{\g_n}^{12}\in {\rm End}(L(\lambda)_{\g_n})\otimes {\rm End}(\mathbb{C}^n)\otimes {\rm End}(\mathbb{C}^m),\\
R^{12}_{\g_n}S^{(1)}_{+}(u)_{\g_n}S^{(2)}_{-}(u)_{\g_n} &=S^{(2)}_{+}(u)_{\g_n}S^{(1)}_{-}(u)_{\g_n} R^{12}_{\g_n}\in {\rm End}(L(\lambda)_{\g_n})\otimes {\rm End}(\mathbb{C}^n)\otimes {\rm End}(\mathbb{C}^m),
\end{align*}
where we take the convention in Theorem \ref{thmRLL}. Here similar to Lemma \ref{conn+3+2}, the matrix $R^{12}_{\g_n}\in {\rm End}(\mathbb{C}^m)\otimes {\rm End}(\mathbb{C}^m)$ is given by the Stokes matrices of \eqref{introqeqBCD}, provided specializing $L(\lambda)_{\g_n}$ to the natural representation of $\g_n$ on $\mathbb{C}^{m}$. In this specialization, the system \eqref{introqeqBCD} becomes a $m^2\times m^2$ (ordinary) system that decomposes into multiple $1\times 1$, $2\times 2$ systems and only one $m\times m$ system. Therefore, under the isomorphism of the two different realizations of quantum groups (see Section \ref{tworealization}), in order to prove an analog of Theorem \ref{introthm1} for type B, C and D, one only needs to compute explicitly the Stokes matrices of the $m\times m$ system and then to check that $R^{12}_{\g_n}$ coincides with the standard R-matrix $R\in {\rm End}(\IC^m)\otimes {\rm End}(\IC^m)$ for type BCD (see e.g., \cite{Jimbo2}\cite{FRT}).

\subsubsection*{Deformation of the classical Stokes matrices: higher order poles and quantum irregular Riemann-Hilbert maps}
The Poisson geometric nature of the irregular Riemann-Hilbert correspondence, beyond the second order pole case, was developed in \cite{Boalch2}. Let us take a positive integer $k$, and consider the differential equations for a function $f(z)\in {\rm GL}(n,\mathbb{C})$
\begin{eqnarray}\label{heq}
\frac{df}{dz} = \left(\frac{u}{z^{k+1}}+\frac{A_k}{z^k}+\cdots+\frac{A_2}{z^2}+\frac{A_1}{z}\right)\cdot f,\end{eqnarray}
where $u\in \h_{\rm reg}$, and $A_i\in\frak{gl}_n$ for all $i=1,...,k$. The differential equation has a pole of order $k+1$ at $z=0$. 

On the one hand, for fixed $u$, the moduli space of the differential equations \eqref{heq} can be described as follows. Let us take the Lie group ${\rm GL}_n(\mathbb{C}[t]/t^{k})$ and let $\frak{gl}_n(\mathbb{C}[t]/t^{k})={\rm Lie}({\rm GL}_n^{(k)})$ be its Lie algebra, which contains elements of the form $X(t) = X_0+X_1t+\cdots + X_{k-1}t^{k-1}$ with $A_i\in\frak{gl}_n$. Then the elements $A(t)$ of the dual $\frak{gl}_n(\mathbb{C}[t]/t^{k})^*$ will be written as
\begin{equation}\label{parameter}
    A(t)=\frac{A_k}{t^k}+\cdots+\frac{A_2}{t^2}+\frac{A_1}{t}
\end{equation}
via the pairing with $\frak{gl}_n(\mathbb{C}[t]/t^{k})$ given by $\langle A(t), X(t)\rangle := {\rm Res}_0(A, X) =\sum_{j=1}^k
(A_j, X_{j-1})$. Note that the equation \eqref{heq} is parameterized by such $A(t)$. 

On the other hand, in this higher order pole case, the Stokes phenomenon of the equation \eqref{heq} at the irregular singularity $z=0$ is characterized by $k$ pairs of Stokes matrices, denoted by $\{S_\pm^{(j)}(u,A(t))\}_{j=1,...,k}$ (with $S_+^{(j)}$ and $S_-^{(l)}$ invertible upper and lower triangular and having the same diagonal part for all $j,l=1,...,k$). See e.g., \cite{BJL, Boalch2}. Let us introduce the space of Stokes matrices
\[\mathcal{M}^{(k)}:=\{(S_-^{(j)}, S_+^{(j)})_{j=1,...,k}\in (B_-\times B_+)^k~|~\delta(S_\pm^{(l)})=\delta(S_-^{(j)}) \ \text{for all } l,j\}.\]

Now the dual $\frak{gl}_n(\mathbb{C}[t]/t^{k})^*$ is equipped with the canonical linear Poisson structure induced from the Lie algebra structure, and $\mathcal{M}^{(k)}$ is equipped with the Poisson structure inherited from the irregular Atiyah-Bott construction \cite{Boalch2}. Then
\begin{thm}\cite{Boalch2}\label{kthm}
The irregular Riemann-Hilbert (also known as Riemann-Hilbert-Birkhoff) map 
\[\nu(u): \frak{gl}_n(\mathbb{C}[t]/t^{k})^*\rightarrow \mathcal{M}^{(k)}~;~A(z)\mapsto \Big(S_-^{(1)},S_+^{(1)},...,S_-^{(k)},S_+^{(k)}\Big)\]
associating the Stokes matrices to the differential equation \eqref{heq}, is a locally analytic Poisson isomorphism. 
\end{thm}
In the case $k=1$, Theorem \ref{kthm} recovers Theorem \ref{sclrll}. 
Here like in the $k=1$ case, it is useful to think that there is only one Poisson structure, i.e., the canonical liner one on $\frak{gl}_n(\mathbb{C}[t]/t^{k})^*$. The map $\nu(u)$ pulls back the entry coordinates on the space of Stokes matrices to the transcendental functions $\{S^{(j)}_\pm(u,A(t))_{il}\}$ on $A(t)\in\frak{gl}_n(\mathbb{C}[t]/t^{k})^*$. Theorem \ref{kthm} then states that for fixed $u$, under the new coordinate system $\{S^{(j)}_\pm(u,A(t))_{il}\}$ around a neighbourhood of $0\in \frak{gl}_n(\mathbb{C}[t]/t^{k})^*$, the canonical linear Poisson bracket takes an algebraic form (the explicit form was computed in \cite{Boalch3} and \cite{Krichever}). From this viewpoint, a quantization of the map $\nu(u)$ is a lift of the functions $\{S^{(j)}_\pm(u,A(t))_{il}\}$ on $A(t)\in\frak{gl}_n(\mathbb{C}[t]/t^{k})^*$ to elements in the completion of $U(\frak{gl}_n(\mathbb{C}[t]/t^{k})).$ Motivated by Theorem \ref{introthm1}, let us discuss how to do it in a canonical way.

\begin{rmk}
In \cite{Boalch3}, an explicit description of the Poisson structure on $\mathcal{M}^{(k)}$ was derived in the framework of quasi-Hamiltonian geometry \cite{AMM}. In particular, the explicit quasi-Hamlitonian two form on an enlarged (framed) moduli space of Stokes matrices, induced from the Atiyah-Bott construction, was given. After the ${\rm GL}_n$-reduction of the global framing, one gets the explicit algebraic Poisson structure on the space $\mathcal{M}^{(k)}$ of Stokes matrices. In the meanwhile, the expression of the Poisson/symplectic structure on $\mathcal{M}^{(k)}$ was derived in \cite{Krichever} via a Hamiltonian approach. 
\end{rmk}
In a next work, we will consider a quantum analog of the equation \eqref{heq}. Recall that the equation \eqref{heq} was parameterized by the dual ${\frak {gl}}_n(\mathbb{C}[t]/t^{k})^*$ equipped with the Lie-Poisson (canonical linear Poisson) bracket, the above discussion in the case $k=1$ then implies that the quantum analog should encode the generators of the universal enveloping algebra $U\big({\frak {gl}}_n(\mathbb{C}[t]/t^{k})\big)$.
Thus we will consider the equation
\begin{equation}\label{hqeq}
\frac{dF}{dz} = h\left(\frac{u}{z^{k+1}}+\frac{T_{[k]}}{z^k}+\cdots+\frac{T_{[2]}}{z^2}+\frac{T_{[1]}}{z}\right)\cdot F,\end{equation}
where each $T_{[m]}$ is an $n\times n$ matrix with entries valued in the universal enveloping algebra $U\big({\frak {gl}}_n(\mathbb{C}[t]/t^{k})\big)$
\begin{eqnarray*}
(T_{[m]})_{ij}=e_{ij}t^{m-1}, \ \ \ \ \ \text{for} \ 1\le i,j\le n, \ \ 1\le m\le k.
\end{eqnarray*}

Let us outline in which sense \eqref{hqeq} is a quantization of \eqref{heq}. First, the equation \eqref{heq} has a formal power series solution
\[\hat{w}=({\rm Id}+w_1z+w_2z^2+\cdots)\cdot  e^{\int^z\big(\frac{u}{t^{k+1}}+\frac{d_k}{t^{k}}+\cdots+\frac{d_1}{t}\big)dt},\]
where $d_k,...,d_1$ are diagonal, and for all positive integers $m$ the entries $\{w_{m,ij}\}_{i,j=1,...,n}$ of the coefficients $w_m\in {\rm End}(\mathbb{C}^n)$ are understood as polynomial functions of the parameter $A(t)\in {\frak {gl}}_n(\mathbb{C}[t]/t^{k})^*$ given in \eqref{parameter}. 

Second, the equation \eqref{hqeq} has a formal power series solution
\[\hat{W}=({\rm Id}+W_1z+W_2z^2+\cdots)\cdot  e^{\int^z\big(\frac{u}{t^{k+1}}+\frac{D_k}{t^{k}}+\cdots+\frac{D_1}{t}\big)dt},\]
where $D_k,...,D_1$ are diagonal with entries in $U({\frak {gl}}_n(\mathbb{C}[t]/t^{k})$, and the entries $\{W_{m,ij}\}_{i,j=1,...,n}$ of the coefficients $W_m$ are elements in $U({\frak {gl}}_n(\mathbb{C}[t]/t^{k}))$. We remark that just like in the proof of Proposition \eqref{oneformal}, the quantum coefficients $W_m$ satisfy the same recursive relation as the classical coefficients $w_m$, provided replacing $A_l$, $d_l$ by $T_{[l]}$, $D_l$ for all $l=1,...,k$.

Then by considering the quantum equation, we actually get a canonical lift of all the data (the analytic data are encoded in the coefficients of the formal solution) from the Poisson algebra to the quantized algebra, i.e., for all positive integers $m$ we have the lift
\[\begin{CD}
\Big\{w_m=(w_{m,ij})_{i,j=1}^n\in {\rm End}(\mathbb{C}^n)\otimes Sym\big({\frak {gl}}_n(\mathbb{C}[t]/t^{k})\big)\Big\} \\
    @V h \text{ deformation } VV \\
    \Big\{W_m=(W_{m,ij})_{i,j=1}^n\in {\rm End}(\mathbb{C}^n)\otimes U\big({\frak {gl}}_n(\mathbb{C}[t]/t^{k})\big)\Big\}.
\end{CD}\]
Taking the Borel-Laplace transform (in the completion of $Sym\big({\frak {gl}}_n(\mathbb{C}[t]/t^{k})\big)$ and the $\mathbb{C}\llbracket h\rrbracket$ completion, with respect to the $h$-adic topology, of $U\big(\frak{{gl}}_n(\mathbb{C}[t]/t^{k})\big)$ respectively), the classical Stokes matrices $\{S_\pm^{(j)}\}_{j=1,...,k}$, seen as $n\times n$ matrix valued analytic functions of the parameter $A(t)\in {\frak {gl}}_n(\mathbb{C}[t]/t^{k})^*$, is then quantized to the quantum Stokes matrices $\{S_{h\pm}^{(j)}\}_{j=1,...,k}$ of \eqref{hqeq}, seen as $n\times n$ matrices with entries valued in $U\big(\frak{{gl}}_n(\mathbb{C}[t]/t^{k})\big)\llbracket h\rrbracket$. 

In this way, a quantization of the irregular Riemann-Hilbert map in Theorem \ref{kthm} is given by the data of quantum Stokes matrices $\{S_{h\pm^{(j)}}\}_{j=1,...,k}$. 
Furthermore, similar to the $k=1$ case in Theorem \ref{introthm1}, we can deduce that the $k$ pairs of (quantum) Stokes matrices of \eqref{hqeq} satisfy some natural algebraic relation, certain $RL\cdots L=L\cdots LR$ relation, see our next work \cite{Xu3}.

When forgetting the transcendental origin from the quantum Stokes matrices, the abstract algebraic relation itself (like the RLL relation in the $k=1$ case) defines an (abstract) associative algebra $\mathcal{A}^{(k)}$ that quantizes the algebraic Poisson structures on the space $\mathcal{M}^{(k)}$ of (classical) Stokes matrices of \eqref{heq}. In this way, the fact that the quantum Stokes matrices $\{S_{h\pm}^{(j)}\}_{j=1,...,k}$ satisfy the same relation amounts to an associative algebra homomorphism $\nu_h(u)$ from $\mathcal{A}^{(k)}$ to $U\big(\frak{{gl}}_n(\mathbb{C}[t]/t^{k})\big)$. Thus the semiclassical limit of the algebra homomorphism $\nu_h(u)$ automatically gives rise to a morphism between Poisson algebras, which coincides with the map $\nu(u)$ as in Theorem \ref{kthm}. Therefore, it interprets the Poisson geometric nature of the classical irregular Riemann-Hilbert map $\nu(u)$, and we get the following diagram
\[
\begin{CD}
\nu_h(u):~\mathcal{A}^{(k)} @> \text{quantum Stokes matrices}>> U\big(\frak{{gl}}_n(\mathbb{C}[t]/t^{k})\big)\\
@V \text{semiclassical limit} VV @V \text{semiclassical limit} VV\\
\nu(u)^*:~ \mathcal{O}(\mathcal{M}^{(k)}) @> \text{classical Stokes matrices}>> \mathcal{O}({\frak {gl}}_n(\mathbb{C}[t]/t^{k})^*).
\end{CD}\]

We stress that, like in $k=1$ case, the quantization of the Poisson structures on $\mathcal{M}^{(k)}$ can be studied purely algebraically. It is independent with the quantization of the irregular Riemann-Hilbert map, which is highly transcendental. 

Let us mention that the Knizhnik–Zamolodchikov (KZ) type equations with irregular singularities have been introduced from various perspectives. For example, the irregular Knizhnik–Zamolodchikov (KZ) equation was introduced in \cite{Re}, and was given a
representation-theoretic interpretation in \cite{FR}. From the perspective of isomonodromy, 
the Casimir connection of De Concini/Millson–Toledano
Laredo \cite{MTL, TL2, TL} governs the isomonodromic deformations of the irregular KZ equation at a second order pole. And the connection matrix of the irregular KZ equation was used to construct the Toledano Laredo twist killing the KZ associator.
The works \cite{Rem, Rem2, Ya} construct a flat connection via quantisation of the
Hamiltonian systems of \cite{Boalch5}, and that flat connection reduces to KZ in the logarithmic
case—so it is also an irregular version of KZ equation. 
The work \cite{GMR} studies the deformation quantisation of isomonodromy equations, in the viewpoint of truncated current
Lie algebras \cite{Wilson}, also known as ‘generalised Takiff Lie algebras’ \cite{Ta, Geo1, Geo2}.
These resulting differential equations can have arbitrary order pole and therefore has Stokes phenomenon. Just like equation \eqref{hqeq}, it is interesting to study their Stokes phenomenon, and to see how the understanding of their Stokes matrices can be used, in the end, to solve the general nonlinear isomonodromy equations \cite{JMU} via the Riemann-Hilbert method.

\section{Stokes matrices of the equation \texorpdfstring{\eqref{introqeq} in the nonresonant cases}{d}}\label{beginsec}
In Section \ref{s11}-\ref{s13}, we introduce the unique formal solution, the canonical solutions and Stokes matrices of the equation \eqref{introqeq} respectively.

\subsection{The unique formal fundamental solution of equation \eqref{introqeq} in nonresonant case}
\label{s11}
\begin{pro}\label{uniformal}
For any $h\notin \mathbb{Q}$ and $u\in\h_{\rm reg}$, the ordinary differential equation
\begin{eqnarray}\label{Stokeseq}
\frac{dF}{dz}=h\Big(u+\frac{T}{z}\Big)\cdot F
\end{eqnarray}
has a unique formal fundamental solution taking the form \begin{eqnarray}\label{formalsum}
\widehat{F}(z)=\widehat{H}(z) e^{{huz}}z^{h\delta T}, \ \ \ {\it for} \ \widehat{H}=1+H_1z^{-1}+H_2z^{-2}+\cdot\cdot\cdot, \end{eqnarray}
where each coefficient $H_m\in{\rm End}(L(\lambda))\otimes{\rm End}(\mathbb{C}^n)$, and $\delta T$ denotes the diagonal part of $T$, i.e., $\delta T=\sum_{k} e_{kk}\otimes E_{kk}.$
\end{pro}
\begin{proof}
Plugging \eqref{formalsum} into the equation \eqref{Stokeseq} gives rise to the equation for $\widehat{H},$
\begin{eqnarray}
\frac{1}{h}\frac{d\hat{H}}{dz}+\hat{H}\cdot \Big(u+
\frac{\delta T}{z}\Big)=\Big( u +
\frac{T}{z}\Big)\cdot \hat{H}.\end{eqnarray}
Comparing the coefficients of $z^{-m-1}$, we see that $H_m$ satisfies 
\begin{eqnarray}\label{simHm}
[H_{m+1}, u]=(\frac{m}{h}+T)\cdot H_{m}-H_{m}\cdot   \delta T.\end{eqnarray}
Let $\{E_{kl}\}_{1\le k,l\le n}$ be the standard basis of ${\rm End}(\mathbb{C}^n)$. Then \[T=\sum_{k,l} e_{kl}\otimes E_{kl}, \ \text{ and } u= \sum_i 1\otimes u_iE_{ii}.\] Plugging $H_m=\sum_{k,l} H_{m, kl}\otimes E_{kl}$, with $H_{m,kl}\in{\rm End}(L(\lambda))$, into the equation \eqref{simHm} gives rise to 
\begin{align} \nonumber
&\sum_{k,l}(u_l-u_k) H_{m+1, kl}\otimes E_{kl}\\ \label{recuH}
=&\sum_{k,l}\frac{m}{h}H_{m,kl}\otimes E_{kl}+
\sum_{k,l,j} e_{kj} H_{m, jl} \otimes E_{kl}-\sum_{k,l}  H_{m, kl} e_{ll}\otimes E_{kl}.
\end{align}
Here $e_{kl}'s$ are understood as elements in ${\rm End}(L(\lambda))$ via the given representation. That is for $k\ne l$ 
\begin{equation}\label{Hknel}
(u_l-u_k) H_{m+1, kl}=\frac{m}{h}H_{m,kl}+\sum_{j=1}^n e_{kj} H_{m, jl}- H_{m, kl} e_{ll} \ \in {\rm End}(L(\lambda)),
\end{equation}
and for $k= l$ (replacing $m$ by $m+1$ in \eqref{recuH}), 
\begin{eqnarray}\label{Hkel}
0=\sum_{j\ne k} e_{kj} H_{m+1, jk}+\frac{m+1}{h} H_{m+1, kk}+[e_{kk}, H_{m+1, kk}]\ \in {\rm End}(L(\lambda)).
\end{eqnarray}

Suppose $H_m$ is given, let us check that the above recursive relation has a unique solution $H_{m+1}$. First note that, since $u_k\ne u_l$ for $k\ne l$, the identity \eqref{Hknel} uniquely defines the "off-diangonal" part $H_{m+1, kl}$ ($k\ne l$) of $H_{m+1}$ from $H_m$. Furthermore, since $h\notin \mathbb{Q}$, we have $\frac{m+1}{h}{\rm Id}+{\rm ad}_{e_{kk}}$ is invertible on ${\rm End}(L(\lambda))$ for any integer $m+1$. Thus, the condition \eqref{Hkel} uniquely defines the "diagonal" part $H_{m+1, kk}$ of $H_{m+1}$ from the off diagonal part. \end{proof} 

\begin{rmk}
The recursive relation \eqref{simHm} can be solved explicitly via the quantum minors of the matrix $T$. See Proposition \ref{Yangsol}.
\end{rmk}

\subsection{The canonical solutions with a prescribed asymptotics}\label{canonicalF}
The radius of convergence of the formal power series $\hat{H}(z)$ is in general zero. However, it follows from the general principle of differential equations with irregular singularities that (see e.g., \cite{Balser,LR, MR} or the proof of Proposition \ref{keypara}) the Borel resummation (Borel-Laplace transform) of $\hat{H}$ gives a holomorphic function in each Stokes supersector around $z=\infty$. In this way, one gets actual solutions of \eqref{Stokeseq} on each Stokes supersector. These sectors are determined by the irregular term $hu$ of the differential equation as follows.

\begin{defi}\label{Stokesrays}
The anti-Stokes rays of the equation \eqref{Stokeseq} are the directions along which $e^{h(u_i-u_j)z}$ decays most rapidly as $z\rightarrow \infty$ for some $u_i\ne u_j$. 
\end{defi}
Equivalently, the anti-Stokes rays
are the directions between pairs of eigenvalues of $hu$, when plotted in the ${z}$-plane. Therefore, there is an even number of such rays.
In this paper we denote a ray by its argument, thus the anti-Stokes rays of \eqref{Stokeseq} are $-{\rm arg}(hu_i-hu_j)$ for all the possible $i\ne j$. 
For any two real numbers $a, b$, an open sector and a closed sector with opening
angle $b-a>0$ are respectively denoted by
\[S(a,b):=\{z\in\mathbb{C}~|~a<{\rm arg}(z)<b\}, \hspace{5mm} \overline{S(a,b)}:=\{z\in\mathbb{C}~|~a\le {\rm arg}(z)\le b\}. \]

\begin{defi}
The Stokes supersectors
are the open regions of $\mathbb{C}$ 
\[{\rm Sect}_i:=S\Big(d_i-\frac{\pi}{2},d_{i+1}+\frac{\pi}{2}\Big),\]
bounded by two adjacent rays $d_i-\frac{\pi}{2}$ and $d_{i+1}+\frac{\pi}{2}$. Here we label the anti-Stokes rays $d_1,d_2,...,d_{2l}$ in anticlockwise order (the indices are taken modulo $2l$), and use the same letter $d_i$ to denote the argument of the ray/direction $d_i$.
\end{defi}
Let us make a choice of branch of ${\rm log}(z)$ on $\Sect_0$, and by convention the branch of ${\rm log}(z)$ on $\Sect_0$ is extended to
the other sectors. Then the following theorem follows from the general principle of differential equations with irregular singularities (see e.g., \cite{Balser,LR, MR} or the proof of Proposition \ref{keypara}).

\begin{thm}\label{thmcanonicalsol}
For fixed $h\notin \mathbb{Q}$ and $u\in\h_{\rm reg}$, there is a unique (therefore canonical) holomorphic function $H_i:{\rm Sect}_i\to {\rm End}(L(\lambda))\otimes {\rm End}(\mathbb{C}^n)$ such that the function
\[F_{i}(z):=H_i (z)\cdot  e^{h uz}z^{h{\delta T}}\]
satisfies the equation \eqref{introqeq}, and moreover
\begin{eqnarray}\label{simH}
H_i(z)\sim \hat{H}(z), \hspace{3mm} \text{as } z\rightarrow\infty \text{  within  } \Sect_i. 
\end{eqnarray}
\end{thm}
We stress that, the actual solution of the differential equation, with the prescribed asymptotic expansion as $z\rightarrow\infty$
in a sector whose opening is larger than $\pi$, is unique and this fact is important in what follows later.

\subsection{Stokes matrices}\label{s13}
\begin{defi}\label{defiStokes}
For any $h\notin \mathbb{Q}$ and $u\in\h_{\rm reg}$, the Stokes matrices $S_{h\pm}(u)\in {\rm End}(L(\lambda))\otimes {\rm End}(\mathbb{C}^n)$ of the equation \eqref{introqeq} (with respect
to the chosen sector $\Sect_0$) are the unique matrices such that:
\begin{itemize}
    \item  If $F_0$ is continued in a positive sense to $\Sect_l$ then $F_l(z)=F_0 (z)\cdot S_{h+}(u)e^{\pi\mathi h\delta T}$, and 
    \item If $F_l$
is continued in a positive sense to $\Sect_0$ then $F_0(ze^{-2\pi \mathi })=F_l(z)\cdot  e^{-\pi\mathi h\delta T} S_{h-}(u)^{-1}$.
\end{itemize}
\end{defi} 

\begin{rmk}
The definition of the Stokes matrices involves the auxiliary choices of
initial sector $\Sect_0$ and a choice of branch of ${\rm log}(z)$ on $\Sect_0$, and the formal monodromy $e^{2\pi\mathi h \delta T}$ usually is separated from the Stokes matrices (thus the latter has one's along the diagonal). Here we take a different convention that the definition of Stokes matrices includes the data of formal monodromy. The advantage is that in our convention the Stokes matrices themselves satisfy the RLL relations. See Section \ref{Uqg}.
\end{rmk}

It is convenient to think of $S_{h\pm}(u)$ as $n\times n$ matrices with entries in ${\rm End}(L(\lambda))$. Let us introduce the permutation matrix $1\otimes J\in {\rm End}(L(\lambda))\otimes {\rm End}(\mathbb{C}^n)$ associated to the choice of $\Sect_0$. Here the $n\times n$ matrix $J$ is
defined by $J_{ij} = \delta_{\sigma(i)j}$, where $\sigma$ is the permutation of $\{1, . . . , n\}$ corresponding to the
dominance ordering of $\{e^{hu_1z},...,e^{hu_nz}\}$ along the direction $\theta$ bisecting the sector $\Sect(d_1, d_l)$. That is
$\sigma(i) < \sigma(j)$ if and only if $e^{h(u_i-u_j)z}\rightarrow 0$ as $z\rightarrow \infty$ along $\theta$. Then the prescribed asymptotics of $F_0(z)$ and $F_l(z)$ at $z=\infty$ ensure that
\begin{lem}\label{upper}
The matrices $S_{h\pm}(u)$ are triangular, up to conjugating
by the permutation matrix $1\otimes J$. Furthermore, their diagonal elements are 
\[S_{h+}(u)_{kk}^{-1}=S_{h-}(u)_{kk}=\mathe^{h\pi \mathi e_{kk}} \in {\rm End}(L(\lambda)), \ \ \text{  for } k=1,...,n.\]
\end{lem}
\begin{proof} By the asymptotics \eqref{simH} of $H_\pm(z)$, Definition \ref{defiStokes} leads to 
\[e^{ huz} z^{h\delta T}\cdot \Big(e^{-\pi\mathi h \delta T} S_{h+}(u)^{-1}\Big) \cdot e^{- huz} z^{-h\delta T} =H_l(z)^{-1} H_0(z) \rightarrow 1, \] 
as $z\rightarrow \infty \text{ within } \Sect_0\cap \Sect_l.$
Let us write $S_{h+}(u)^{-1}=\sum_{ij} (S_{h+}^{-1})_{ij}\otimes E_{ij}\in {\rm End}(L(\lambda))\otimes {\rm End}(\mathbb{C}^n)$. Since the exponential $e^{huz}$ dominate, we must have, for any $i\ne j$, $e^{h(u_i-u_j)z}(S_{h+}^{-1})_{ij}\otimes E_{ij} \rightarrow 0.$
It implies that  $(e^{-\pi\mathi h  \delta T} S_{h+}^{-1})_{ij}=\delta_{ij}$ unless $e^{z h(u_i-u_j)}\rightarrow 0$ as $z\rightarrow \infty$ within $\Sect_0\cap \Sect_l$. Thus $S_{h+}(u)^{-1}$ is triangular, up to the permutation matrix $1\otimes J$. The argument for $S_{h-}$ is the same once the change of the branches of ${\rm log}(z)$ is accounted for.
\end{proof} 

\vspace{2mm}
For the sake of convenience, in the rest of the paper, we will assume that $u\in\h_{\rm reg}$ and the sector $\Sect_0$ are chosen such that the corresponding permutation
matrix is the identity. Under this assumption, the Stokes matrices $S_{h+}(u)$ and $S_{h-}(u)$ are upper and lower triangular respectively. That is, $S_{h+}(u)_{ij}=0\in {\rm End}(L(\lambda))$ if $i>j$, and  $S_{h-}(u)_{ij}=0\in {\rm End}(L(\lambda))$ if $i<j$.

\begin{rmk}
There is a factorization of the Stokes matrices $S_{h\pm}(u)$ as an ordered product of the Stokes factors $S_i(u)$ for $i=0,1,...,l-1$ (up to a multiplication by $e^{\pi\mathi h\delta T}$). Here each $S_i(u)\in {\rm End}(L(\lambda))\otimes {\rm End}(\mathbb{C}^n)$ is determined by $F_{i+1}(z)=F_i (z)\cdot S_i(u)$ within $\Sect_i\cap \Sect_{i+1}$. It is interesting to see if the factorization of the Stokes matrices is related to the factorization of
the universal R–matrices of Drinfeld-Jimbo quantum groups \cite{KR,LS}.
\end{rmk}

\section{Stokes matrices satisfy the RLL relation}\label{Uqg}
Let us take the standard R-matrix $R\in {\rm End}(\IC^n)\otimes {\rm End}(\IC^n)$, see e.g., \cite{Jimbo2}\cite{FRT},
\begin{equation}\label{sRmatrix}
R=\sum_{i\ne j, i,j=1}^n E_{ii}\otimes E_{jj}+e^{\pi\mathi h}\sum_{i=1}^n E_{ii}\otimes E_{ii}+(e^{\pi\mathi h}-e^{-\pi\mathi h})\sum_{1\le j<i\le n}E_{ij}\otimes E_{ji}.
\end{equation} 
The following theorem is the main result of this section. 
Note that the Stokes matrices $S_{h\pm}(u)$ are $n\times n$ matrices with entries $S_{h\pm}(u)_{ij}$ in ${\rm End}(L(\lambda))$, i.e., $S_{h\pm}(u)=\sum_{i,j}S_{h\pm}(u)_{ij}\otimes E_{ij}\in {\rm End}(L(\lambda))\otimes {\rm End}(\mathbb{C}^n)$ . Let us introduce \[S^{(1)}_{h\pm}(u):= \sum_{i,j}S_{h\pm}(u)_{ij}\otimes E_{ij}\otimes 1, \ S^{(2)}_{h\pm}(u):= \sum_{i,j}S_{h\pm}(u)_{ij}\otimes 1\otimes E_{ij}, \text{ and } \ R^{12}:=1\otimes R\] as elements in ${\rm End}(L(\lambda))\otimes  {\rm End}(\mathbb{C}^n)\otimes {\rm End}(\mathbb{C}^n)$.

\begin{thm}\label{thmRLL}
For any $h\notin\mathbb{Q}$ and $u\in\h_{\rm reg}$, the Stokes matrices $S_{h\pm}(u)$ of \eqref{introqeq} satisfy
\begin{equation}\label{RLL1}
R^{12}S^{(1)}_{h\pm}(u)S^{(2)}_{h\pm}(u) =S^{(2)}_{h\pm}(u)S^{(1)}_{h\pm}(u) R^{12},
\end{equation}
\begin{equation}\label{RLL2}
R^{12}S^{(1)}_{+}(u)S^{(2)}_{-}(u) =S^{(2)}_{+}(u)S^{(1)}_{-}(u) R^{12}.
\end{equation}
\end{thm}
This theorem is closely related to the result in \cite{Xu1}. To prove it, we introduce a partial differential equation for a function $Y(z,t)$ with two complex variables. We then introduce six holomorphic solutions $Y_{\pm k}(z,t)$ of the equation in six different connected regions $D_{\pm k}$ for $k=1,2,3$, specified by a common prescribed asymptotics. In the end, we prove that after the analytic continuation of $Y_{+3}$ from the domain $D_{+3}$ to the domain $D_{-3}$, along two homotopic paths $D_{+3}\rightarrow D_{+2}\rightarrow D_{+1}\rightarrow D_{-3}$ and $D_{+3}\rightarrow D_{-1}\rightarrow D_{-2}\rightarrow D_{-3}$, we have respectively
\begin{eqnarray*}
Y_{-1}(z,t)=Y_{+3}(z,t)\cdot S^{(2)}_{h+}(u)S^{(1)}_{h+}(u)R^{(12)}, \text{  for  } (z,t)\in D_{-3}.
\end{eqnarray*}
and
\begin{eqnarray*}
Y_{-1}(z,t)=Y_{+3}(z,t)\cdot R^{(12)}S^{(2)}_{h+}(u)S^{(1)}_{h+}(u), \text{  for  } (z,t)\in D_{-3}.
\end{eqnarray*}
It implies $S^{(2)}_{h+}S^{(1)}_{h+}R^{(12)}=R^{(12)}S^{(2)}_{h+}S^{(1)}_{h+}$. Similarly, we obtain the other identities in Theorem \ref{thmRLL}.

\subsection{Solutions of the Knizhnik–Zamolodchikov equations with irregular singularities}\label{sec3first}
Again, let $L(\lambda)$ be the finite dimensional highest weight
$U(\frak{gl}_n)$-module. Let us consider a system of equations, for an ${\rm End}(L(\lambda))\otimes {\rm End}(\IC^n)\otimes {\rm End}(\IC^n)$ valued function $Y(z_1,z_2)$ of two complex variables, 
\begin{eqnarray}\label{gcKZ1}
&&\frac{1}{h}\frac{{\partial Y}}{\partial z}= \Big( u^{(1)}+tu^{(2)}+
\frac{T^{(1)}+T^{(2)}+P}{z}\Big)\cdot Y, \\ \label{gcKZ2}
&&\frac{1}{h}\frac{{\partial Y}}{\partial t}= \Big(z u^{(2)}+\frac{T^{(2)}}{t}+\frac{P}{t- 1}\Big)\cdot Y.
\end{eqnarray}
Here $h$ is a complex parameter, $u={\rm diag}(u_1,...,u_n)\in\h_{\rm reg}(\mathbb{R})$, and
\begin{align*}
&u^{(1)}= \sum_i 1\otimes u_i E_{ii}\otimes 1, \ u^{(2)}= \sum_i 1\otimes1 \otimes u_i E_{ii},& \\ 
&T^{(1)}=\sum_{k,l} e_{kl}\otimes E_{kl}\otimes 1, \ T^{(2)}=\sum_{k,l} e_{kl}\otimes 1\otimes E_{kl}, \ P=-\sum_{k,l} 1\otimes E_{kl}\otimes E_{lk},&
\end{align*}
are elements in $U(\frak{gl}_n)\otimes {\rm End}(\IC^n)\otimes {\rm End}(\IC^n)$. The action of the coefficient matrix of the system on $Y(z_1,z_2)$ is given by matrix multiplication and the representation of ${\frak {gl}}_n$.
\begin{pro}
The system of equations \eqref{gcKZ1} and \eqref{gcKZ2} is compatible.
\end{pro}
\begin{proof} It follows from a direct computation and the identities $[T^{(1)}+T^{(2)}, P]=0$, $[T^{(1)}, T^{(2)}+P]=0$. Here we remark that there is a minus sign in the expression of $P$. \end{proof}
\begin{rmk}
The system is equivalent to the two variables Knizhnik–Zamolodchikov (KZ) equations with irregular singularities \cite{FMTV}. Under the assumption that the eigenvalues of $hu$ are purely imaginary, the canonical solutions of the KZ equation with prescribed asymptotics was first studied in a formal setting by Toledano Laredo \cite{TL}. Furthermore, the Stokes phenomenon was used to construct a Drinfeld twist killing the KZ associator. The proof in this section is mainly motivated by \cite{TL}.
\end{rmk}
In the rest of this section, we fix a $u={\rm diag}(u_1,...,u_n)\in\h_{\rm reg}$. For any fixed $t$, the first equation \eqref{gcKZ1} becomes a meromorphic ordinary differential equation with an irregular singularity at $z=\infty$. In the following, we first follow the standard way to produce solutions of equation \eqref{gcKZ1} with a prescribed asymptotics at $z=\infty$, and then apply the resulting solutions to solve the system of equations \eqref{gcKZ1} and \eqref{gcKZ2}.

\subsubsection{The unique formal solution of equation \eqref{gcKZ1} in the nonresonant case}
For any $u\in\h_{\rm reg}$ let us introduce the $u$ dependent region in $\IC$,
\begin{eqnarray*}
D_{u}:=\Big\{t\in \IC~|~t\ne \frac{u_i-u_j}{u_k-u_l}, \text{ for any } i,j, k,l\Big\}.\end{eqnarray*} 
In the following, we will first consider the case $t\in D_u$, and then treat the leftover case by means of analytic continuation. 

\begin{pro}\label{formalsol}
Let us assume that $h\notin \mathbb{Q}$. Then for any fixed $u\in\h_{\rm reg}$ and $t\in D_u$, the ordinary differential equation \eqref{gcKZ1} has a unique formal fundamental solution taking the form \begin{eqnarray}\label{2formalsum}
\widehat{Y}(z,t)=\widehat{Q}(z;t) e^{ hz(u^{(1)}+tu^{(2)})}z^{h(\delta T^{(1)}+\delta T^{(2)}+\delta P)}, \end{eqnarray}
where $\widehat{Q}=1+Q_1z^{-1}+Q_2z^{-2}+\cdot\cdot\cdot$, and each coefficient $Q_m(t)\in {\rm End}(L(\lambda))\otimes {\rm End}(\IC^n)\otimes {\rm End}(\IC^n)$. Here $\delta$ takes the diagonal part, i.e.,
\[\delta T^{(1)}=\sum_{k} e_{kk}\otimes E_{kk}\otimes 1, \ \ \ \delta T^{(2)}=\sum_{k} e_{kk}\otimes 1\otimes E_{kk}, \ \ \ \delta P=-\sum_{k} 1\otimes E_{kk}\otimes E_{kk}.\]
\end{pro}
\begin{proof}
Plugging $\hat{Y}$ in \eqref{2formalsum} into the equation \eqref{gcKZ1} gives rise to the equation for the formal power series $\widehat{Q}$,
\begin{align}\label{Hhat}\begin{split}
&\frac{1}{h}\frac{d\hat{Q}}{dz}+\hat{Q}\cdot \Big(u^{(1)}+t u^{(2)} +
\frac{\delta T^{(1)}+\delta T^{(2)}+\delta P}{z}\Big)\\
=&\Big(u^{(1)}+t u^{(2)}+
\frac{T^{(1)}+T^{(2)}+P}{z}\Big)\cdot \hat{Q}.\end{split}\end{align}
Comparing the coefficients of $z^{-m-1}$ in the two sides of \eqref{Hhat}, we get that the recursion relation
\begin{align}\label{simQm}\begin{split}
&[Q_{m+1}(t), u^{(1)}+t u^{(2)}]\\
=&\Big( m/h+T^{(1)}+T^{(2)}+P\Big)\cdot Q_{m}(t)-Q_{m}(t)\cdot  \Big(\delta T^{(1)}+\delta T^{(2)}+\delta P\Big).\end{split}\end{align}
Let us write $Q_m=\sum_{i,j,k,l} Q_{m, ijkl}\otimes E_{ij}\otimes E_{kl}$ in terms of the basis $\{E_{ij}\otimes E_{kl}\}_{1\le i,j,k,l \le n}$ of ${\rm End}(\mathbb{C})^{\otimes 2}$, where each $Q_{m, ijkl}\in{\rm End}(L(\lambda))$. Then a direct computation shows that the relation \eqref{simQm} gives rise to 

$(1)$ For $i\ne j$ or $k\ne l$: 
\begin{align}\label{offdia}\begin{split}
&(u_j-u_i+tu_l-tu_k) Q_{m+1, ijkl}\\ 
=&\sum_{j'} e_{ij'} Q_{m, j'jkl}+\sum_{l'} e_{kl'} Q_{m, ijl'l}+\delta_{lj}Q_{m, ijkl}-Q_{m, kjil}\\
&-Q_{m, ijkl}e_{jj}- Q_{m, ijkl}e_{ll}+\frac{m}{h}Q_{m, ijkl}.
\end{split}
\end{align}

$(2)$ For $i=j$ and $k=l$: 
\begin{align}\label{diagQ}\begin{split}
0=&\sum_{i'\ne i} e_{ii'} Q_{m, i'ikk}+\sum_{k'\ne k} e_{kk'} Q_{m, iik'k}+\delta_{ik} Q_{m, ikki}-Q_{m, kiik}\\
&+[e_{ii}+e_{kk}, Q_{m, iikk}]+\frac{m}{h}Q_{m, iikk}.\end{split}
\end{align}

Let us check that the above recursive relation have a unique solution. First, by the assumption $t\in D_u$, i.e., $t\ne \frac{u_i-u_j}{u_k-u_l}$ for all $i,j, k,l$, the identity \eqref{offdia} uniquely defines the "off-diangonal" part $Q_{m+1, ijkl}$ ($i\ne j$ or $k\ne l$) of $Q_{m+1}$ from $H_m$. Second, replacing $m$ by $m+1$ in the identity \eqref{diagQ} (the solvability condition of the recursion relation), we get a linear system
\begin{eqnarray} \nonumber
&&[ Q_{m+1, iikk},e_{ii}+e_{kk}]- (m+1)/h \cdot Q_{m+1, iikk}\\ \label{diarecursion}
&=&\sum_{i'\ne i} e_{ii'} Q_{m+1, i'ikk}+\sum_{k'\ne k} e_{kk'} Q_{m+1, iik'k}+\delta_{ik}Q_{m+1, ikki}-Q_{m+1, kiik}
\end{eqnarray}
for the unknown $Q_{m+1, iikk}\in{\rm End}(L(\lambda))$. Note that the off-diagonal part of $Q_{m+1}$ on the right hand side of \eqref{diarecursion} has been determined by the identity \eqref{offdia}.
Since $h\notin \mathbb{Q}$, we have $(m+1)/h{\rm Id}+{\rm ad}_{e_{ii}+e_{kk}}$ is invertible on ${\rm End}(L(\lambda))$. That is, the linear system for $Q_{m+1,iikk}$ is nondegenerate and has a unique solution. This concludes the proof. \end{proof}

\subsubsection{Anti-Stokes rays and canonical solutions of the equation \eqref{gcKZ1} with prescribed asymptotics}\label{canonicalY}

In the rest of this section, let us fix $h\notin\mathbb{Q}$ and $u\in\h_{\rm reg}$. 
Similar to Section \ref{canonicalF}, the following definitions and theorem follow from the general principle of the Borel resummation of the formal solutions of the meromorphic differential equations. 
\begin{defi}\label{tStokesray}
The anti-Stokes rays of the equation \eqref{gcKZ1} are the directions along which $e^{h(u_i-u_j+tu_k-tu_l)z}$ decays most rapidly as $z\rightarrow \infty$ for some $u_i\ne u_j$ or $u_k\ne u_l$. 
\end{defi}
Thus the rays have the arguments $-{\rm arg}\big(h(u_i-u_j)+ht(u_k-u_l)\big)$ for some index $i,j,k,l$.
\begin{defi}\label{tsector}
Let $d(t)$ and $d'(t)$ be two adjacent anti-Stokes rays (in anticlockwise order). Then the corresponding Stokes supersectors are the open regions of the complex plane
\[{\rm Sect}_{d(t),d'(t)}(t):=S\Big(d(t)-\frac{\pi}{2},d'(t)+\frac{\pi}{2}\Big).\]
\end{defi}
Henceforth, to simply the notation, we write $d,d'$ for $d(t), d'(t)$ when they appear in the subscript.
We stress that the anti-Stokes rays and Stokes supersectors depend on $t$. For fixed $t$ the Borel resummation of $\hat{Q}(z)$ defines holomorphic functions in each Stokes supersector, which in turn determines the holomorphic solutions of \eqref{gcKZ1} with the prescribed asymptotics. In summary, we have 
\begin{thm}\label{BorelY}
For fixed $t\in D_u$, there is a unique (therefore canonical) holomorphic function $Q_{d,d'}(z;t):{\rm Sect}_{d,d'}(t)\to {\rm End}(L(\lambda))\otimes {\rm End}(\mathbb{C}^n)$ such that the function
\[\widetilde{Y}_{d,d'}(z,t):=Q_{d,d'}(z;t)\cdot  e^{hz(u^{(1)}+tu^{(2)})}z^{h(\delta T^{(1)}+\delta T^{(2)}+\delta P)}\]
satisfies the equation \eqref{gcKZ1}, and moreover
\begin{eqnarray}\label{simQ}
Q_{d,d'}(z;t)\sim \hat{Q}(z;t), \hspace{3mm} \text{ as } z\rightarrow \infty \text{  within  } \Sect_{d,d'}(t). 
\end{eqnarray}
\end{thm}
There are two kinds of anti-Stokes rays. One is the set of rays that vary with respect $t$, another is the set of fixed rays with arguments $-{\rm arg}\big(h(u_i-u_j))$ for some index $i,j$ (independent of $t$). Note that the fixed rays coincide with the anti-Stokes rays of equation \eqref{introqeq}, and like in Section \ref{canonicalF} we label these rays by $d_0,...,d_{2l-1},d_{2l}=d_0$. 
\begin{figure}[h]
\begin{center}
\includegraphics[scale=0.35]{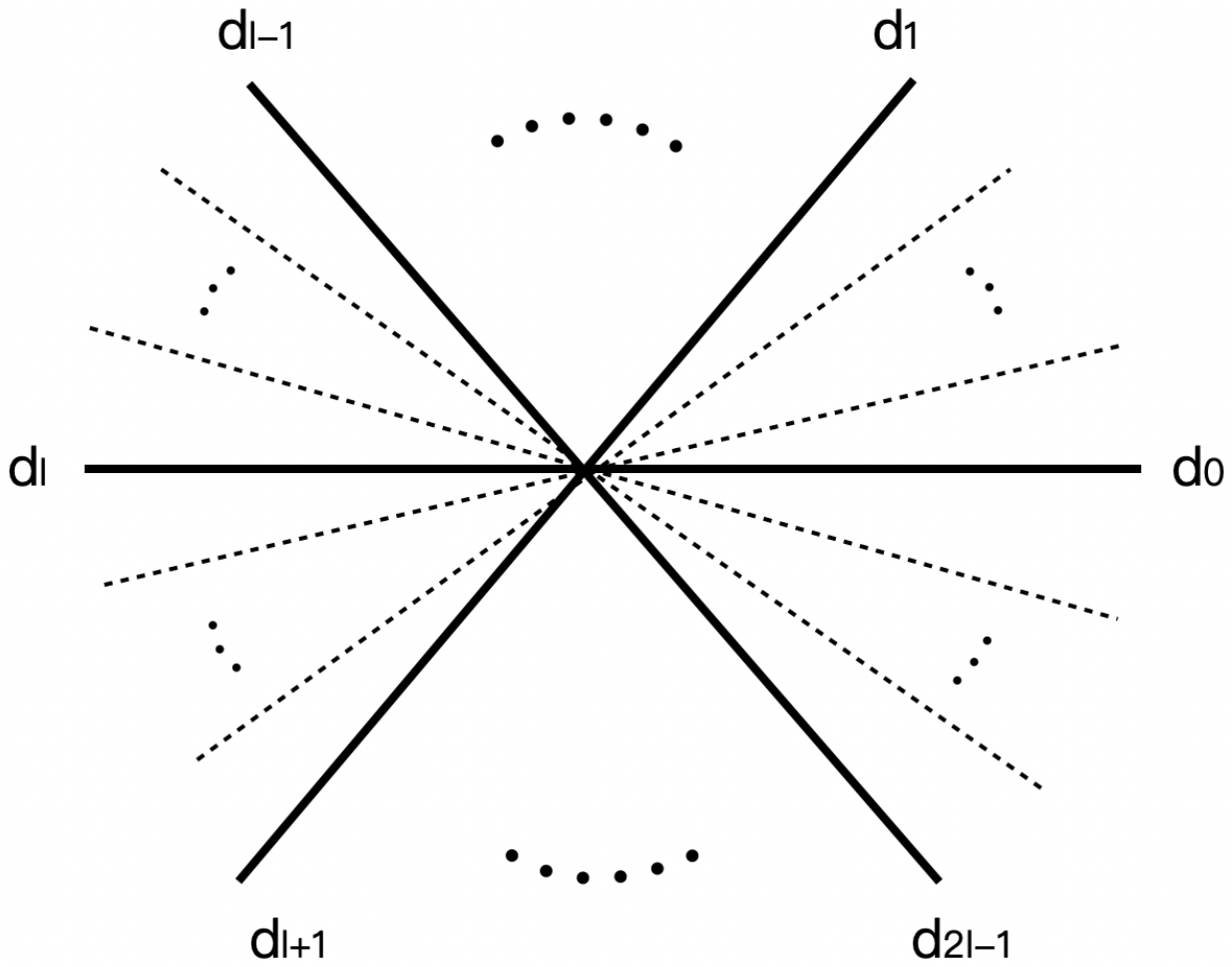}
\end{center}
\caption{The anti-Stokes rays of the equation \eqref{gcKZ1}. The solid lines are the fixed rays (independent of $t$) labelled by $d_0,...,d_{2l-1}, d_{2l}=d_0$. The dashed lines are the rays varying with respect to $t$.}
\par\end{figure}

\subsubsection{Solutions of the system of equations \eqref{gcKZ1} and \eqref{gcKZ2} on different domains}
Let $\mathcal{R}\in D_u$ be any region of the complex $t$-plane such that there exist two adjacent anti-Stokes rays $d(t)$ and $d'(t)$ that continuously dependent on $t\in D_u$ and never collide. For any $t\in \mathcal{R}$, let $\widetilde{Y}_{d,d'}$ be the function given in Theorem \ref{BorelY} defined on the corresponding $t$-dependent sector $\Sect_{d,d'}(t)$. Then we have
\begin{thm}\label{solution}
The function
\begin{eqnarray}\label{Ypm}
Y_{d,d'}(z,t)=\widetilde{Y}_{d,d'}(z,t)\cdot t^{h\delta T^{(2)}}(t-1)^{h\delta P}, 
\end{eqnarray}
defined for $t\in\mathcal{R}$ and $z\in\Sect_{d,d'}(t)$,
satisfies the system of equations \eqref{gcKZ1} and \eqref{gcKZ2}. 
\end{thm}
\begin{proof} 
By the compatibility of the equations \eqref{gcKZ1} and \eqref{gcKZ2}, we have that the functions
\begin{eqnarray}\label{Ezt}
U_{d,d'}(z,t):= \frac{dY_{d,d'}(z,t)}{dt}-h\Big(z u^{(2)}+\frac{T^{(2)}}{t}+\frac{P}{t- 1}\Big)\cdot Y_{d,d'}
\end{eqnarray}
satisfy the equation \eqref{gcKZ1} in the corresponding domains. It implies that the ratio 
\begin{equation}\label{B_k}
C:=Y_{d,d'}^{-1}\cdot U_{d,d'}=Y_{d,d'}^{-1} \frac{dY_{d,d'}}{dt}-Y_{d,d'} ^{-1}\cdot  h\Big(z u^{(2)}+\frac{T^{(2)}}{t}+\frac{P}{t- 1}\Big)\cdot Y_{d,d'}
\end{equation}
is independent of $z$.
To show that $C$ is zero, let us rewrite 
\begin{align}\label{QtoY}\begin{split}
Y_{d,d'}&=\widetilde{Y}_{d,d'}(z,t)\cdot t^{h\delta T^{(2)}}(t-1)^{h\delta P}\\
&=Q_{d,d'}(z,t)  e^{hz(u^{(1)}+tu^{(2)})} z^{h{\delta T^{(1)}}}(tz)^{h{\delta T^{(2)}}}(tz-z)^{h{\delta P}}.
\end{split}
\end{align}
Differentiating the above identity with respect to $t$ gives
\begin{equation*}
\frac{dY_{d,d'}}{dt} \cdot Y_{d,d'} ^{-1}= \frac{dQ_{d,d'}}{dt} Q_{d,d'} ^{-1}+Q_{d,d'}\cdot h\left(zu^{(2)}+\frac{\delta T^{(2)}}{t}+\frac{\delta P}{t-1} \right)\cdot Q_{d,d'}^{-1}.
\end{equation*}
Following from \eqref{simQ}, we have
\begin{equation*}
Q_{d,d'}(z,t)\sim \hat{Q}(z,t)= 1+Q_1(t) z^{-1}+Q_2(t) z^{-2}+\cdot\cdot\cdot,
\end{equation*}
as $z\rightarrow \infty$ within $\Sect_{d,d'}(t)$.
By \eqref{simQm}, the first term $Q_1(t)$ should satisfy \begin{eqnarray}\label{recH1}
[Q_1(t),u^{(1)}+t u^{(2)}]=T^{(1)}+T^{(2)}+P-\delta T^{(1)}-\delta T^{(2)}-\delta P.\end{eqnarray}
Solving \eqref{recH1} gives the off diagonal part of $Q_1(t)=\sum_{i,j,k,l}H_{1,ijkl}\otimes E_{ij}\otimes E_{kl}$,
\begin{align}\begin{split}
    \label{Q1}
&Q_1(t)-\sum_{i,k}Q_{1,iikk}\otimes E_{ii}\otimes E_{kk}\\
=&\sum_{1\le i\ne j
\le n}\Big( \frac{e_{ij}\otimes E_{ij}\otimes 1}{u_j-u_i}+\frac{1}{t}\frac{e_{ij}\otimes 1\otimes E_{ij}}{u_j-u_i}+\frac{1}{t-1}\frac{-1\otimes E_{ij}\otimes E_{ji}}{u_j-u_i}\Big).\end{split}
\end{align}
Just as in the proof of Proposition \ref{formalsol}, the assumption $t\in D_u$ makes the right hand side of \eqref{Q1} well-defined. Since the diagonal part $\sum_{i,k}Q_{1,iikk}\otimes E_{ii}\otimes E_{kk}$ of $Q_1$ commutes with $u^{(2)}$, the identity \eqref{Q1} implies that
\[[Q_1(t), u^{(2)}]= \frac{T^{(2)}-\delta T^{(2)}}{t}+\frac{P-\delta P}{t-1}. \]
Therefore, the above identity and the asymptotics \eqref{simQ} of $Q_{d,d'}$ lead to
\[Q_{d,d'}\cdot h\Big(zu^{(2)}+\frac{\delta T^{(2)}}{t}+\frac{\delta P}{t-1} \Big)\cdot Q_{d,d'}^{-1}=h\Big(zu^{(2)}+\frac{T^{(2)}}{t}+\frac{P}{t-1}\Big)+O(z^{-1}),\]
as $z\rightarrow \infty$ within $\Sect_{d,d'}(t)$.
That is,
\begin{eqnarray*}
Y_{d,d'} \cdot C \cdot Y_{d,d'} ^{-1}=\frac{dY_{d,d'}}{dt}\cdot Y_{d,d'} ^{-1}-h\Big(zu^{(2)}+\frac{T^{(2)}}{t}+\frac{P}{t-1}\Big)=O(z^{-1}),\end{eqnarray*}
as $z\rightarrow\infty$ within $\Sect_{d,d'}(t)$.
The above identity can be rewritten as
\begin{align}\begin{split}
    &e^{hz(u^{(1)}+tu^{(2)})} z^{h{\delta T^{(1)}}}(tz)^{h{\delta T^{(2)}}}(tz-z)^{h{\delta P}}\cdot C \cdot e^{-hz(u^{(1)}-tu^{(2)})} z^{-h{\delta T^{(1)}}}(tz)^{-h{\delta T^{(2)}}}(tz-z)^{-h{\delta P}}\\ \label{sim}
    =& Q_{d,d'}(z,t)^{-1} \cdot O(z^{-1}) \cdot Q_{d,d'}(z,t).
    \end{split}
\end{align}
Since $Q_{d,d'}(z,t)\sim 1$, we get that the left hand side of the above identities is $O(z^{-1})$ as $z\rightarrow\infty$ within $\Sect_{d,d'}(t)$.
Let us write $C=Y_{d,d'} ^{-1}\cdot U_{d,d'}=\sum_{a,b,c,d}C_{i, abcd}\otimes E_{ab}\otimes E_{cd}$. Since the exponential term $e^{h{z}(u^{(1)}+tu^{(2)})}$ dominates, the asymptotics \eqref{sim} forces that for any $a,b,c,d$, 
\[e^{hz(u_a-u_b+tu_c-tu_d)}  C_{i, abcd}\otimes E_{ij}\otimes E_{kl}\rightarrow 0, \text{ as } z\rightarrow \infty \text{ within } \Sect_{d,d'}(t).\]
Since the sector $\Sect_i$ has an opening angle bigger than $\pi$, it implies that all the elements $C_{i,abcd}$ must vanish. In the meanwhile, for the diagonal elements, we have
\begin{equation}\label{diazero}
\Big(z^{h(e_{aa}+e_{cc})}\cdot  C_{i,  aacc}\cdot z^{-h(e_{aa}+e_{cc})}\Big)  \otimes E_{aa}\otimes E_{cc}=O(z^{-1}),
\end{equation}
as $z\rightarrow \infty$ within $\Sect_{d,d'}(t)$.
If $h$ is a purely imaginary number, the formula \eqref{diazero} implies $C_{i,  iikk}=0$. For a generic $h\in\mathbb{C}\setminus\mathbb{Q}$, $C_{i,  iikk}=0$ by means of analytic continuation. Since $C=Y_{d,d'}^{-1}\cdot U_{d,d'}$ is zero,
note the fact that $Y_{d,d'}$ is a fundamental solution implies $U_{d,d'}(z,t)=0$, and thus finishes the proof. 
\end{proof}

\subsubsection{Stokes factors}
For fixed $u$ and $h$, let us denote by $d(t), d'(t), d''(t)$ three adjacent anti-Stokes rays of the equation \eqref{gcKZ1} in anticlockwise order, at a point $t\in D_u$. Following \eqref{Ypm}, denote by $Y_{d,d'}$ and $Y_{d',d''}$ the solutions on the sectors $\Sect_{d,d'}(t)$ and $\Sect_{d',d''}(t)$ as in Definition \ref{tsector} respectively.

\begin{defi}
The Stokes factor $S_{d'}$ associated to the ray $d'(t)$ is defined by 
\begin{eqnarray}
  Y_{d',d''}(z,t)=  Y_{d,d'}(z,t)\cdot S_{d'} \hspace{3mm} \text{ in } \Sect_{d,d'}(t)\cap \Sect_{d',d''}(t).
\end{eqnarray}
\end{defi}

\begin{lem}\label{Sfactor}
If we write $S_{d'}=\sum_{1\le i,j,k,l\le n} c_{ijkl}\otimes E_{ij}\otimes E_{kl}$, then the element $c_{ijkl}\in{\rm End}(L(\lambda))$ must be zero unless $d'(t)$ coincides with the ray $-{\rm arg}\big(h (u_i-u_j+t(u_k-u_l))\big)$ with the same index $i,j,k,l$.
\end{lem}
\begin{proof} The proof of Lemma \ref{upper} can be applied. \end{proof} 

\begin{rmk}
In Section \ref{trivialS} we will prove that $S_{d'}$ must be the identity unless the anti-Stokes ray $d'(t)$ has the argument $-{\rm arg}\big(h(u_i-u_j)\big)$ for some $i\ne j$. That is only the fixed anti-Stokes rays of the form $-{\rm arg}\big(h(u_i-u_j)\big)$ (independent of $t$) can cause the jump of the solutions with prescribed asymptotics, while the anti-Stokes rays $-{\rm arg}\big(h (u_i-u_j+t(u_k-u_l))\big)$, varying with respect to $t$, do not cause any jump.
\end{rmk}
Note that as $t$ varies, the configuration of the anti-Stokes rays varies accordingly. In particular, as $t$ crosses from one side of some point to the other, some rays may collide and some new rays may emerge. In the following, we will determine all the Stokes factors of \eqref{gcKZ1} for $t\in\mathbb{R}$, and thus solve the connection problem between the solutions defined on different regions given in Theorem \ref{solution}.

\subsection{Preferred solutions of the system of equations \eqref{gcKZ1} and \eqref{gcKZ2} on certain domains}\label{secsolution}
For the region $R=(-\delta,0)$, $(0,\delta)$ or $(\frac{1}{\delta},\infty)$ with $\delta$ a sufficient small positive real number, we will introduce two adjacent anti-Stokes rays $d(t)$ and $d'(t)$ of the equation \eqref{gcKZ1} such that they continuously dependent on $t\in R$ and never collide. Then Theorem \ref{solution} enables us to take solution defined for $t\in\mathcal{R}$ and $z\in\Sect_{d,d'}(t)$. 

\begin{figure}[h]
\begin{center}
\includegraphics[scale=0.3]{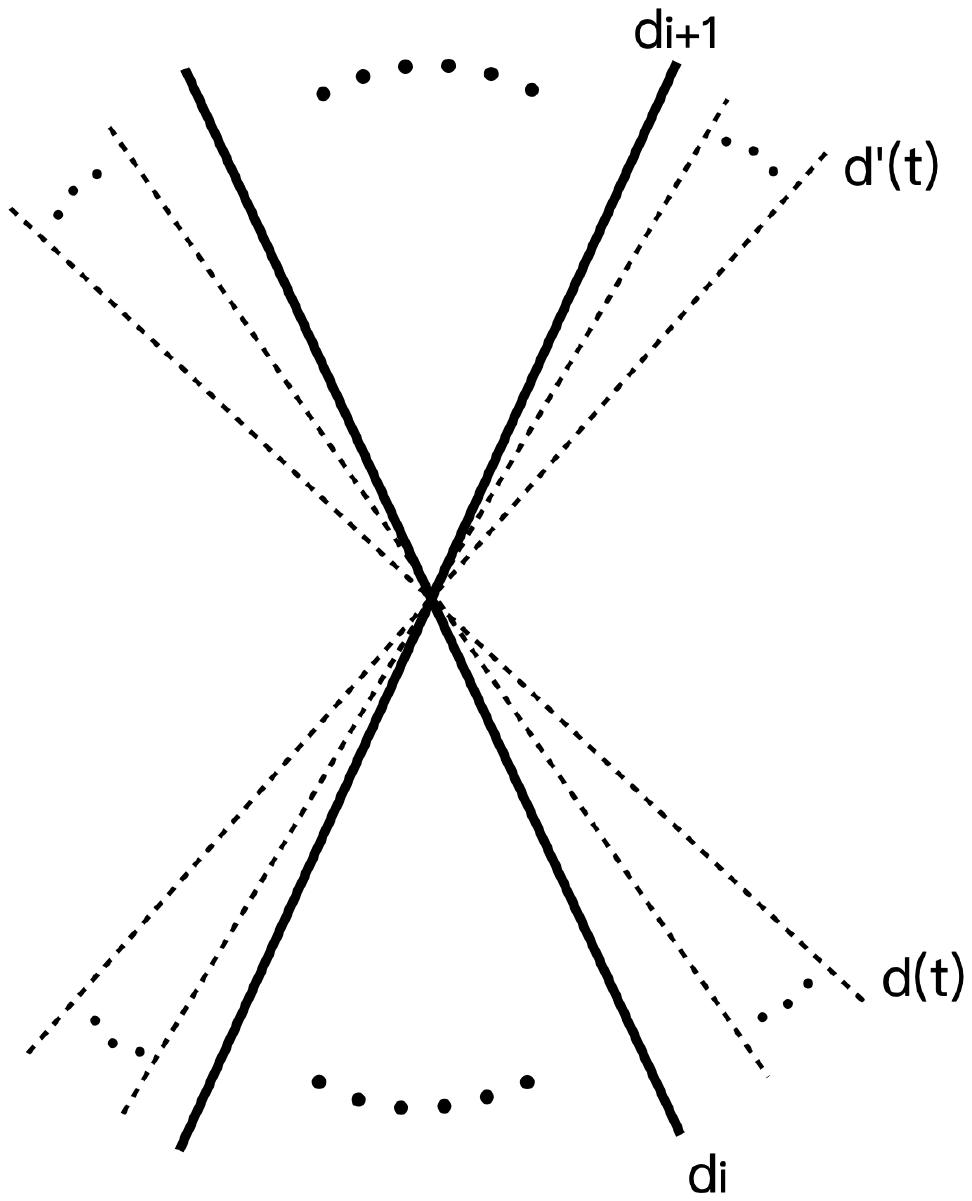}
\end{center}
\caption{The configuration of the anti-Stokes rays as $t$ close to $0$ (or $\infty$). The solid lines are the fixed rays (independent of $t$), and two adjacent such rays are labelled by $d_i$ and $d_{i+1}$. The dashed lines, the rays varying with respect to $t$, are close to the solid lines as $t$ near $0$ (or $\infty$). And $d(t), d'(t)$ represent the two adjacent rays approaching to $d_i$ and $d_{i+1}$ as $t\rightarrow 0$ (or $\infty$) respectively.}
\end{figure}
First recall from Definition \ref{tStokesray} that all the anti-Stokes rays of \eqref{gcKZ1} at $t$ have the arguments $-{\rm arg}\big(h(u_i-u_j)+ht(u_k-u_l)\big)$ for some index $i,j,k,l$. All these rays will approach to one of the fixed anti-Stokes rays $-{\rm arg}\big(h(u_i-u_j)\big)$ of \eqref{introqeq} as $t\rightarrow 0$ or $t\rightarrow \infty$ (without causing ambiguity we also denote a ray by its argument). 
Let us fix two (separate) adjacent anti-Stokes rays $d_i$ and $d_{i+1}$ of \eqref{introqeq}. Then associated to the two fixed rays $d_i$ and $d_{i+1}$ of \eqref{introqeq}, there exists a unique pair of (separate) adjacent anti-Stokes rays $d(t)$ and $d'(t)$ of \eqref{gcKZ1} defined for all $t\in (-\delta,0)\cup (0,\delta)\cup (1/\delta,+\infty)$ (where $\delta$ is a small enough real number), and such that $d(t)$ and $d'(t)$ approach the rays $d_i$ and $d_{i+1}$ of \eqref{introqeq} as $t\rightarrow 0-$, $t\rightarrow 0+$ and $t\rightarrow+\infty$ respectively. See Figure 2 for an illustration. We remark that the rays $d(t), d'(t)$ depend continuously on $t\in (-\delta,0)\cup (0,\delta)$, but may have a discontinuous jump as $t$ passes through $0$.

Then we can take a real number $\theta>0$ such that the sector 
\[\Sect_i^\theta:=(d_{i}+
\theta-\pi/2,d_{i+1}-\theta+\pi/2)\] has opening angle bigger than $\pi$, and is contained in $\Sect_{d,d'}(t)$ for all $t\in (-\delta,0)\cup (0,\delta)\cup (1/\delta, +\infty)$. By Theorem \ref{BorelY}, we get functions $Q_{d,d'}(z,t)$ defined on the three regions $\Sect_i^\theta \times (-\delta,0)$, $\Sect_i^\theta \times (0,\delta)$ and $\Sect_i^\theta \times (1/\delta,+\infty)$. According to Theorem \ref{solution}, we then get the corresponding solutions $Y_{d,d'}(z,t)$ of the system of equations \eqref{gcKZ1} and \eqref{gcKZ2} on the three different regions.

The analytic continuation of the solutions
\[Y_{d,d'}(z,t)= Q_{d,d'}(z;t)\cdot  e^{hz(u^{(1)}+tu^{(2)})}z^{h(\delta T^{(1)}}zt^{h\delta T^{(2)}}(zt-z)^{h\delta P)},\] 
from $\Sect_i^\theta \times (-\delta,0)$, $\Sect_i^\theta \times (0,\delta)$ and $\Sect_i^\theta \times (1/\delta,+\infty)$ to $\Sect_i^\theta \times (-\infty,0)$, $\Sect_i^\theta \times (0,1)$ and $\Sect_i^\theta \times (1,+\infty)$, with respect to the $t$ variable, will be denoted by
\begin{equation}\label{yik}
Y_{i1}(z,t), \hspace{2mm} Y_{i2}(z,t) \hspace{2mm} \text{and} \hspace{2mm} Y_{i3}(z,t)
\end{equation}
respectively.

We emphasize that for $1\le k\ne l\le 3$ the solution $Y_{ik}(z,t)$, continued to the domain of $Y_{il}(z,t)$ along a path on the complex $t$ plane, in general is not equal to $Y_{il}(z,t)$. See Lemma \ref{conn+2+1} for the connection formula between these solutions.

\subsection{Factorization properties of the solutions $Y_{ik}$ for $k=1,2,3$ at $t=0, 1, \infty$}\label{factorization}
Under the coordinate change $z_1=z$ and $z_2=zt$, the system of equations \eqref{gcKZ1} and \eqref{gcKZ2} becomes
\begin{eqnarray*}
&&\frac{1}{h}\frac{{\partial Y}}{\partial z_1}= \Big( u^{(1)}+\frac{T^{(1)}}{z_1}+
\frac{P}{z_1-z_2}\Big)\cdot Y, \\
&&\frac{1}{h}\frac{{\partial Y}}{\partial z_2}= \Big( u^{(2)}+\frac{T^{(2)}}{z_2}+
\frac{P}{z_2-z_1}\Big)\cdot Y.
\end{eqnarray*}
And under the coordinate change $\omega_1=z$ and $\omega_2=zt-z$, the system becomes
\begin{eqnarray*}
&&\frac{1}{h}\frac{{\partial Y}}{\partial \omega_1}= \Big( u^{(1)}+u^{(2)}+\frac{T^{(1)}}{\omega_1}+
\frac{T^{(2)}}{\omega_1+\omega_2}\Big)\cdot Y, \\
&&\frac{1}{h}\frac{{\partial Y}}{\partial \omega_2}= \Big( u^{(2)}+\frac{T^{(2)}}{\omega_1+\omega_2}+
\frac{P}{\omega_2}\Big)\cdot Y.
\end{eqnarray*}
Using the same method as in Section \ref{sec3first}, we can find canonical solutions of the system with prescribed asymptotics (in different asymptotic regions, see Theorem \ref{fact}) in terms of the new coordinates $(z_1, z_2)$ and $(\omega_1,\omega_2)$ respectively. 

In this subsection, we state the following theorem for the connection formula between these solutions. The proof of this theorem is long. Rather than to interrupt the present discussion we shall postpone that to the Appendix \ref{lastpiece}.

\begin{thm}\label{fact}
Let ${Y}_{ik}(z,t)$ for $k=1,2,3$ be the solutions given in the notation \eqref{yik}.

$(a).$ If $t\in (-\infty, 0)$, then 
\begin{eqnarray}\label{parta}
\widetilde{W}_{i+l}(z,zt)F_{i}^{(1)}(z)=Y_{i1}(z,t)=W_i(z,zt)F_{i+l}^{(2)}(zt)e^{\pi\mathi  h\delta P}.
\end{eqnarray}

$(b).$ If $t\in (0, 1)$, then
\begin{eqnarray}
W_i(z,zt)F_{i}^{(2)}(zt)e^{\pi\mathi  h\delta P}=Y_{i2}(z,t)=\overline{W}_i(z,zt-z)X^{(2)}_{l+i}(zt-z).
\end{eqnarray}

$(c).$ If $t\in (1,\infty)$, then
\begin{eqnarray}
\overline{W}_i(z,zt-z)X_{i}^{(2)}(zt-z)=Y_{i3}(z,t)=\widetilde{W}_{i}(z,zt)F_{i}^{(1)}(z).
\end{eqnarray}
\end{thm}
Here
\begin{itemize}
    \item $F_i^{(1)}(z)$ and $F_i^{(2)}(z)$ are the unique solutions, defined on $\Sect_i=S\Big(d_i-\frac{\pi}{2},d_{i+1}+\frac{\pi}{2}\Big)$, of the equations for ${\rm End}(L(\lambda))\otimes {\rm End}(\IC^n)\otimes {\rm End}(\IC^n)$ valued functions 
\begin{eqnarray}\label{app1}
\frac{1}{h}\frac{d F}{d z}= \Big(u^{(1)}+ \frac{T^{(1)}}{z}\Big)\cdot F, 
\end{eqnarray}
and 
\begin{eqnarray}\label{app2}
\frac{1}{h}\frac{d F}{d z}= \Big(u^{(2)}+ \frac{T^{(2)}}{z}\Big)\cdot F
\end{eqnarray}
respectively, with the prescribed asymptotics as $z\rightarrow \infty$ within  $\Sect_i$
\[F_i^{(1)}(z)\cdot e^{-h u^{(1)}z}z^{-h{\delta T^{(1)}}}\sim 1 \text{ and }F_i^{(2)}(z)\cdot e^{-h u^{(2)}z}z^{-h{\delta T^{(2)}}}\sim 1.\] 
By definition, $F_i^{(1)}$ and $F_i^{(2)}$ are respectively the two different extensions of the function $F_i(z)$ given in Theorem \ref{thmcanonicalsol} from ${\rm End}(L(\lambda))\otimes {\rm End}(\IC^n)$ to ${\rm End}(L(\lambda))\otimes {\rm End}(\IC^n)\otimes {\rm End}(\IC^n)$. Recall from Definition \ref{Stokesrays} that there are $2l$ anti-Stokes rays in total, and  $\Sect_i$ and $\Sect_{l+i}$ are two opposite Stokes supersectors of \eqref{introqeq}.

\item $X_i^{(2)}(z)=K^{(2)}_i(z) e^{hu^{(2)}z}z^{h\delta P}$ is the unique solution, defined on $\Sect_i$, of the equation for a ${\rm End}(L(\lambda))\otimes {\rm End}(\IC^n)\otimes {\rm End}(\IC^n)$ valued function
\begin{eqnarray}\label{Pequation}
\frac{1}{h}\frac{d X}{d z}= \Big(u^{(2)}+ \frac{P}{z}\Big)\cdot X,
\end{eqnarray}
with the holomorphic part $K^{(2)}_i(z)$ tending to $1$ as $z\rightarrow \infty$ within $\Sect_i$.

\item 
For any fixed $z_2$, $W_i(z_1,z_2)$ is the unique solution of the equation for an ${\rm End}(L(\lambda))\otimes {\rm End}(\IC^n)\otimes {\rm End}(\IC^n)$ valued function
\begin{eqnarray}\label{ratio1eq}
\frac{1}{h}\frac{d W}{d z_1}= \Big(u^{(1)}+\frac{T^{(1)}}{z_1}+\frac{P}{z_1- z_2}\Big)\cdot W,
\end{eqnarray}
defined for $z_1\in\Sect_i$ (the sector $S(d_i-{\pi}/{2},d_{i+1}+{\pi}/{2})$ of the $z_1$ plane) and $|z_1|>|z_2|$, and at the same time as $z_1\rightarrow \infty$ within $\Sect_i$
\begin{eqnarray}\label{simW1}
W_i(z_1,z_2)\cdot e^{-h u^{(1)}z_1}z_1^{-h{\delta T}}(z_1-z_2)^{-h\delta P}\sim 1. 
\end{eqnarray}
See Lemma \ref{esolW} for the existence of such a $W_i$. Similarly, for any fixed $z_1$, $\widetilde{W}_i(z_1,z_2)$ is the unique solution of \begin{eqnarray}\label{ratio2eq}
\frac{1}{h}\frac{d W}{d z_2}= \Big(u^{(2)}+\frac{T^{(2)}}{z_2}+\frac{P}{z_2- z_1}\Big)\cdot W,
\end{eqnarray}
defined for $z_2\in\Sect_i$ (the sector $S(d_i-{\pi}/{2},d_{i+1}+{\pi}/{2})$ of the $z_2$ plane) and $|z_2|>|z_1|$, and at the same time as $z_2\rightarrow \infty$ within $\Sect_i$
\begin{eqnarray}\label{simW2}
\widetilde{W}_i(z_1,z_2)\cdot e^{-h u^{(2)}z_2}z_2^{-h{\delta T}}(z_2-z_1)^{-h\delta P}\sim 1. 
\end{eqnarray}
\item For any fixed $\omega_2$, $\overline{W}_i(\omega_1,\omega_2)$ is the unique solution of the equation
\begin{eqnarray}\label{ratio3eq}
    \frac{1}{h}\frac{dW}{d\omega_1}= \Big(u^{(1)}+u^{(2)}+\frac{T^{(1)}}{\omega_1}+\frac{T^{(2)}}{\omega_1+\omega_2}\Big)\cdot W 
\end{eqnarray}
defined on $\Sect_i$ (the sector $S(d_i-{\pi}/{2},d_{i+1}+{\pi}/{2})$ of the $\omega_1$ plane), with the prescribed asymptotics as $\omega_1\rightarrow \infty$ within $\Sect_i$ 
\begin{eqnarray}\label{simW3}
\overline{W}_i(\omega_1,\omega_2)\cdot e^{-h(u^{(1)}+ u^{(2)})\omega_1}\omega_1^{-h{\delta T^{(1)}}}(\omega_1+\omega_2)^{-h\delta T^{(2)}}\sim 1.
\end{eqnarray}
\end{itemize}

\subsection{The connection formula on the complex $t$ plane}\label{connformula}
Section \ref{factorization} finds the connection matrices between the solutions defined in Theorem \ref{solution} for real $t\in \IR\setminus\{0, 1\}$. In this subsection, we find the connection matrices between these solutions on the complex $t$ plane. To match up with the definition of Stokes matrices (see the proofs of Lemma \ref{conn+3+2} to Lemma \ref{conn+1-3}), let us introduce the solutions modified by the corresponding formal half monodromy.

\begin{defi}\label{sixsols}
For any $k=1,2,3,$ we introduce the solutions $Y_{\pm k}$ on the regions $D_{\pm k}$
\begin{align*}
Y_{+3}(z,t)&=Y_{03}(z,t), & D_{+3}&=\{\ \ z\in\Sect_0,  \ \ \ 1<t\ \ \},\\
Y_{+2}(z,t)&= Y_{02}(z,t)\cdot e^{-\pi\mathi h\delta P},  & D_{+2}&=\{z\in\Sect_0, 0<t<1\},\\
Y_{+1}(z,t)&= Y_{01}(z,t)\cdot e^{-\pi\mathi h(\delta P+\delta T^{(2)})}, & D_{+1}&=\{\ \ z\in\Sect_0, \ \  \  t<0 \ \ \},\\
Y_{-1}(z,t)&= Y_{l1}(z,t)\cdot e^{-\pi\mathi h\delta T^{(1)}},  & D_{-1}&=\{z\in\Sect_l, \ \ \ \ \ \ t<0\},\\
Y_{-2}(z,t)&= Y_{l2}(z,t)\cdot e^{-\pi\mathi h(\delta T^{(1)}+\delta T^{(2)})}, & D_{-2}&=\{z\in\Sect_l, 0<t<1\},\\
Y_{-3}(z,t)&= Y_{l3}(z,t)\cdot e^{-\pi\mathi h(\delta P+\delta T^{(1)}+ \delta T^{(2)})}, & D_{-3}&=\{\ \ z\in\Sect_l, \ \ \ 1<t\ \ \}.
\end{align*} 
\end{defi}

From now on, let us switch to the complex $t$ plane. The solutions $Y_{\pm k}(z,t)$ are in general different, after being analytically continued to a common domain on the complex $t$ plane. In the following, we derive the connection formula between them.

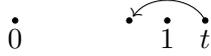
\begin{figure}[htb]
\[
  \begin{tikzpicture}
    [scale=0.5, baseline={([yshift=-.5ex]current bounding box.center)}]
\filldraw (-2,0) circle (2pt) node at (-2,-0.5) {$0$};
       
\filldraw (1,0) circle (2pt) node at (1,0) {};

\filldraw (2,0) circle (2pt) node at (2,-0.5) {$1$};

    \filldraw (3,0) circle (2pt) node at (3,-0.5) {$t$};
    
    \path[->, bend right=60] (3,0)  edge (1.1,0.1);
    \end{tikzpicture}
\]
   \caption{Transposition of $t$ near $1$ such that $t$ passes above $1$.}
\end{figure}

\begin{lem}\label{conn+3+2}
Under the assumption that $u\in\h_{\rm reg}$ and the sector $\Sect_0$ are chosen such that the Stokes matrices $S_{h+}(u)$ and $S_{h-}(u)$ are upper and lower triangular respectively, we have that
\begin{eqnarray}\label{connect2}
Y_{+3}(z,t)=Y_{+2}(z,t)\cdot R^{(12)},
\end{eqnarray}
after the analytic continuation of $Y_{+3}(z,t)$ to the defining domain $D_{+2}$ of $Y_{+2}$ along a simple positive path in the complex $t$ plane, as shown in Figure 1. Here $R$ is the standard R-matrix given in \eqref{sRmatrix} and recall that $R^{(12)}=1\otimes R\in {\rm End}(L(\lambda))\otimes {\rm End}(\IC^n)\otimes {\rm End}(\IC^n)$.
\end{lem}
\begin{rmk}
The extension of the Lemma \ref{conn+3+2} to the general case is simple, provided the permutation matrix is accounted for. 
\end{rmk}
\begin{proof} It follows from the first identity in part $(c)$ of Theorem \ref{fact} that \begin{eqnarray*}
Y_{+3}(z,t)=\overline{W}_0(z,zt-z)X_0^{(2)}(zt-z),
\end{eqnarray*} 
and from the second identity in part $(b)$ of Theorem \ref{fact} that
\begin{eqnarray*}
Y_{+2}(z,t)=\overline{W}_0(z,zt-z)X_l^{(2)}(zt-z)e^{-\pi\mathi\delta P}.
\end{eqnarray*}
Here the corresponding index $i$ of the anti-Stokes ray $d_i$ in Theorem \ref{fact} is just $i=0$. Note that $\overline{W}_0(z_1,z_2)$ does not contribute monodromy from the continuation on $z_2$ plane. Thus, after the analytic continuation of $Y_{+3}$, we have
\[Y_{+2}(z,t)^{-1}Y_{+3}(z,t)=e^{\pi\mathi\delta P}X_l^{(2)}(zt-z)^{-1}X_0^{(2)}(zt-z)=S_{h+}^{(2)}(P)^{-1}.\]
Here $S_{h+}^{(2)}(P)$ is simply the Stokes matrix of \eqref{Pequation} associated to $\Sect_0$, i.e., the ratio of $K^{(2)}_0$ and $K^{(2)}_l$ in the intersection $\Sect_0\cap\Sect_l$. 
The rest is to show that the inverse of the Stokes matrix $S_{h+}(P)^{-1}$ of the equation 
\begin{eqnarray}\label{Requation}
    \frac{dX}{dz}=h\Big(u+\frac{P}{z}\Big)X
\end{eqnarray}
for a ${\rm End}(L(\lambda))\otimes {\rm End}(\mathbb{C}^n)$ valued function, equals to $R$. The canonical solutions and the Stokes matrices of the equation \eqref{Pequation} are the extensions of the ones of \eqref{Requation} from ${\rm End}(L(\lambda))\otimes {\rm End}(\mathbb{C}^n)$ to ${\rm End}(L(\lambda))\otimes {\rm End}(\mathbb{C}^n)\otimes {\rm End}(\mathbb{C}^n)$.

Let us first consider the ${\frak {gl}}_2$ example. That is to consider the equation \eqref{Requation} with \[u=\scriptsize{\left(\begin{array}{cccc}
u_1&0&0&0\\
0&u_1&0&0\\
0&0&u_2&0\\
0&0&0&u_2
\end{array}\right)} \hspace{10mm} and \hspace{10mm} P=-\scriptsize{\left(\begin{array}{cccc}
1&0&0&0\\
0&0&1&0\\
0&1&0&0\\
0&0&0&1
\end{array}\right)}.\] It reduces to the computation of the $2\times 2$ system
\begin{eqnarray}\label{2b2}
    \frac{dX}{dz}={\footnotesize\left(
  \begin{array}{cc}
    hu_1 & 0  \\
    0 & hu_2
  \end{array}
\right)}X
+ \frac{1}{z}{\footnotesize\left(
  \begin{array}{cc}
    0 &-h \\
    -h & 0
  \end{array}
\right)}X.
\end{eqnarray}
Applying the formula of Stokes matrices of a general system of rank $2$ in \cite[Proposition 8]{BJL}, see also \cite{Xu} and Section \ref{2by2}, the Stokes matrices of \eqref{2b2} are expressed by the gamma functions of eigenvalues $\pm h$ of the residue matrix. Using the Euler's reflection formula $\Gamma (1-z)\Gamma (z)=\frac{\pi}{\sin(\pi z)}$ (for $z\notin\mathbb{Z}$) to rewrite the gamma function by hyperbolic function, we get the Stokes matrix $S_{h+}$ of \eqref{2b2}
\[S_{h+}^{-1}=\scriptsize{\left(\begin{array}{cc}
1&e^{\frac{h}{2}}-e^{-\frac{h}{2}}\\
0&1
\end{array}\right)}.\] Here our assumption ensures that $S_{h+}$ is upper triangular. Thus the Stokes matrix $S_{h+}(P)$ of the $4\times 4$ system is
\[S_{h+}(P)^{-1}=\scriptsize{\left(\begin{array}{cccc}
e^{\frac{h}{2}}&0&0&0\\
0&1&e^{\frac{h}{2}}-e^{-\frac{h}{2}}&0\\
0&0&1&0\\
0&0&0&e^{\frac{h}{2}}
\end{array}\right)}.\] 
Note that it doesn't depend on $u={\rm diag}(u_1,u_2)$, and coincides with the R-matrix $R$ defined in \eqref{sRmatrix} for ${\frak {gl}}_2$.

For general $n$, the equation \eqref{Requation} is a system of rank $n^2$. Let $\{v_i\}_{1\le i\le n}$ be the standard basis of $\mathbb{C}^n$, and let us take a basis $\{v_i\otimes v_j\}$ of $\mathbb{C}^n\otimes \mathbb{C}^n$ (with an order remained to be fixed). Note that $P$ is a permutation operator, i.e., $P(v_i\otimes v_j)=v_j\otimes v_i$. And $u$ is a diagonal operator, i.e., $u(v_i\otimes v_j)=u_j \cdot v_i\otimes v_j$ for all $i,j=1,...,n$. Then for fixed index $i\ne j$, the equation \eqref{Requation} can be restricted on the subspace spanned by the two vectors $v_i\otimes v_j$ and $v_j\otimes v_i$, reducing to a $2\times 2$ system of the form 
\begin{equation*}
    \frac{dX}{dz}={\footnotesize\left(
  \begin{array}{cc}
    hu_j & 0  \\
    0 & hu_i
  \end{array}
\right)}X
+ \frac{1}{z}{\footnotesize\left(
  \begin{array}{cc}
    0 &-h \\
    -h & 0
  \end{array}
\right)}X.
\end{equation*}
At the same time, the equation \eqref{Requation} can be restricted on the one dimensional space $\mathbb{C}v_i\otimes v_i$, reducing to a $1\times 1$ system of the form 
$\frac{df}{dz}=hu_if$ for a function $f(z)\in\mathbb{C}$. 

Therefore, the linear system \eqref{Pequation} can be decomposed to multiple $2\times 2$ and $1\times 1$ systems. It thus reduces to the computation of $\frak{gl}_2$ case. The lemma then follows by a direct computation, as long as the assumption on the orders of $u_1,...,u_n$ are accounted for. \end{proof}

\vspace{2mm}
Now let us take the analytic continuation of $Y_{+2}$ from $0<t<1$ to $t<0$ along a simple positive path in the complex $t$ plane, as shown in Figure 2.
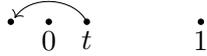
\begin{figure}[htb]
\[
  \begin{tikzpicture}
    [scale=0.5, baseline={([yshift=-.5ex]current bounding box.center)}]

    \filldraw (1,0) circle (2pt) node at (1,-0.5) {$ $};
    
\filldraw (2,0) circle (2pt) node at (2,-0.5) {$0$};

    \filldraw (3,0) circle (2pt) node at (3,-0.5) {$t$};

    \filldraw (6,0) circle (2pt) node at (6,-0.5) {$1$};
    \path[->, bend right=60] (3,0) edge (1.1,0.1);
    \end{tikzpicture}
\]
   \caption{Transposition of $t$ near $0$ such that $t$ passes above $0$.}
\end{figure}
\begin{lem}\label{conn+2+1}
For the connection formula from $Y_{+2}$ to $Y_{+1}$ along the above path, we have that,
\begin{eqnarray}\label{connect21}
Y_{+2}(z,t)=Y_{+1}(z,t)\cdot S^{(2)}_{h+}(u)^{-1}.
\end{eqnarray}
\end{lem}
\begin{proof} It follows from the first identity of part $(b)$ and the second identity of part $(a)$ of Theorem \ref{fact} that \begin{eqnarray*}
Y_{+2}(z,t)=W_0(z,zt)F_0^{(2)}(zt)e^{\pi\mathi  h\delta P}\cdot e^{-h\pi\mathi\delta P}=W_0(z,zt)F_0^{(2)}(zt),
\end{eqnarray*}
and 
\begin{align*}
Y_{+1}(z,t)=&W_0(z,zt)F_l^{(2)}(zt)e^{\pi\mathi  h\delta P}\cdot e^{-h\pi\mathi(\delta P+\delta T^{(2)})}\\
=&W_0(z,zt)F_0^{(2)}(zt)\cdot e^{-h\pi\mathi\delta T^{(2)}}.
\end{align*}
Note that $W_0(z_1,z_2)$ does not contribute monodromy from the continuation on $z_2$ plane. Thus, after the analytic continuation we have
\[Y_{+1}(z,t)^{-1}Y_{+2}(z,t)=e^{-h\pi\mathi\delta T^{(2)}}F_l^{(2)}(zt)^{-1}F_0^{(2)}(zt)=S_{h+}^{(2)}(u)^{-1}.\]
Here $S_{h+}^{(2)}(u)$ is simply the Stokes matrix of \eqref{app2} associated to $\Sect_0$, i.e., the ratio of $F^{(2)}_0$ and $F^{(2)}_l$ in the intersection $\Sect_0\cap\Sect_l$.  \end{proof}

\vspace{2mm}
For any $z\in \Sect_+$, let us take the analytic continuation of $Y_{+3}(z,t)$ from $t<0$ to $t>1$ along a simple positive path in the complex $t$ plane, as shown in Figure 3.
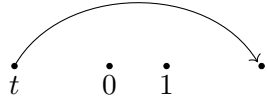
\begin{figure}[htb]
\[
  \begin{tikzpicture}
    [scale=0.5, baseline={([yshift=-.5ex]current bounding box.center)}]
    \filldraw (-2,0) circle (2pt) node at (-2,-0.5) {$t$};
    
\filldraw (0.5,0) circle (2pt) node at (0.5,-0.5) {$0$};
    
\filldraw (2,0) circle (2pt) node at (2,-0.5) {$1$};

\filldraw (4.5,0) circle (2pt) node at (4.5,0) { };

    \path[bend left=60,->] (-2,0) edge (4.4,0.1);
    \end{tikzpicture}
\]
   \caption{Transposition of $t$ near $\infty$ such that $t$ passes below $0$ and $1$.}
\end{figure}
\begin{lem}\label{conn+1-3}
After the analytic continuation of $Y_{+1}(z,t)$, we have
\begin{eqnarray}\label{connect3}
Y_{+1}(z,t)=Y_{-3}(z,t)\cdot S^{(1)}_{h+}(u)^{-1} \hspace{5mm} \text{in the domain}  \ \{ z\in\Sect_+\cap \Sect_-, \ 1<t\}.
\end{eqnarray}
\end{lem}
\begin{proof} It follows from the first identity of part $(a)$ and the second identity of part $(c)$ of Theorem \ref{fact} that
\begin{eqnarray*}
Y_{+1}(z,t)=\widetilde{W}_l(z,zt)F_0^{(1)}(z)e^{-h\pi\mathi(\delta P+\delta T^{(2)})},
\end{eqnarray*}
and
\begin{eqnarray*}
Y_{-3}(z,t)=\widetilde{W}_l(z,zt)F_l^{(1)}(z)e^{-h\pi\mathi(\delta P+\delta T^{(2)}+\delta T^{(1)})},
\end{eqnarray*}
respectively.
Note that $W_0(z_1,z_2)$ is defined for all $z_2\in\mathbb{C}$ and does not contribute monodromy from the continuation on $z_2$ plane. Thus, after the analytic continuation of $Y_{+1}$ we have
\begin{align*}
    Y_{-3}(z,t)^{-1}Y_{+1}(z,t)=&e^{h\pi\mathi(\delta P+\delta T^{(2)})} e^{\pi\mathi\delta T^{(1)}}F_l^{(1)}(z)^{-1}F_0^{(1)}(z)e^{-h\pi\mathi(\delta P+\delta T^{(2)})}\\=&S_{h+}^{(1)}(u)^{-1}.
\end{align*}
Here $S_{h+}^{(1)}(u)$ is simply the Stokes matrix of \eqref{app1} associated to $\Sect_0$, i.e., the (modified by $e^{\pi\mathi\delta T^{(1)}}$) ratio of $F^{(1)}_0$ and $F^{(1)}_l$ in the intersection $\Sect_0\cap\Sect_l$. And the second identity follows from the fact that $e^{h\pi\mathi(\delta P+\delta T^{(2)})}$ commutes with the coefficient matrix of the equation \eqref{app1} and the regularized initial value at $z=\infty$, and therefore commutes with the solutions $F^{(1)}_0$ and $F^{(1)}_l$ as well as the Stokes matrix $S_{h+}^{(1)}(u)$.  \end{proof}

\subsection{The $\rm RLL=LLR$ relation: a proof of Theorem \ref{thmRLL}}\label{RLLLLR}
As a consequence of Lemma \ref{conn+3+2}, \ref{conn+2+1} and \ref{conn+1-3},
\begin{pro}\label{end1}
After the analytic continuation of $Y_{+3}(z,t)$ from the domain $D_{+3}$ to $D_{-3}$, in the manner $D_{+3}\rightarrow D_{+2}\rightarrow D_{+1}\rightarrow D_{-3}$ specified in Lemma \ref{conn+3+2}-\ref{conn+1-3}, we have
\begin{eqnarray}\label{RLLid1}
Y_{+3}(z,t)=Y_{-3}(z,t)\cdot S^{(1)}_{h+}(u)^{-1}S^{(2)}_{h+}(u)^{-1} R^{(12)}, \ \ \text{  for  } (z,t)\in D_{-3}.
\end{eqnarray}
\end{pro}
In the meanwhile, we can also take the analytic continuation of $Y_{+3}(z,t)$ from the domain $D_{+3}$ to $D_{-3}$, in a manner $D_{+3}\rightarrow D_{-1}\rightarrow D_{-2}\rightarrow D_{-3}$ along the following three paths on the $t$ plane given in Figure 4-6 respectively

\begin{figure}[htb]
\[
  \begin{tikzpicture}
    [scale=0.5, baseline={([yshift=-.5ex]current bounding box.center)}]
    \filldraw (-2,0) circle (2pt) node at (-2,-0.5) {$0$};
    
\filldraw (-0.5,0) circle (2pt) node at (-0.5,-0.5) {$1$};
    
\filldraw (2,0) circle (2pt) node at (2,-0.5) {$t$};

\filldraw (-4.5,0) circle (2pt) node at (-4.5, 0) { };

    \path[->, bend left=60] (2,0) edge (-4.4,-0.1);
    \end{tikzpicture}
\]
   \caption{Transposition of $t$ connecting $D_{+3}$ and $D_{-1}$.}
\end{figure}
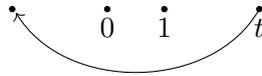

\begin{figure}[htb]
\[
  \begin{tikzpicture}
    [scale=0.5, baseline={([yshift=-.5ex]current bounding box.center)}]

    \filldraw (1,0) circle (2pt) node at (1,-0.5) {$t$};
    
\filldraw (2,0) circle (2pt) node at (2,-0.5) {$0$};

   \filldraw (3,0) circle (2pt) node at (3, 0) { };

    \filldraw (6,0) circle (2pt) node at (6,-0.5) {$1$};
    
    \path[->, bend right=90] (1,0) edge (3,-0.1);
    \end{tikzpicture}
\]
   \caption{Transposition of $t$ near $0$ connecting $D_{-1}$ and $D_{-2}$.}
\end{figure}
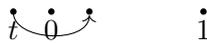
   
\begin{figure}[htb]
\[
  \begin{tikzpicture}
    [scale=0.5, baseline={([yshift=-.5ex]current bounding box.center)}]

    \filldraw (5,0) circle (2pt) node at (5,-0.5) {$t$};
    
\filldraw (2,0) circle (2pt) node at (2,-0.5) {$0$};

    \filldraw (7,0) circle (2pt) node at (7,0) { };

    \filldraw (6,0) circle (2pt) node at (6,-0.5) {$1$};
    
    \path[->, bend right=90] (5,0) edge (7,-0.1);
    \end{tikzpicture}
\]
   \caption{Transposition of $t$ near $1$ connecting $D_{-2}$ and $D_{-3}$.}
   \end{figure}
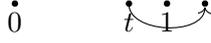
   
Then analog to Lemma \ref{conn+3+2}-\ref{conn+1-3}, we have
\begin{lem}
After the analytic continuation of $Y_{+3}(z,t)$ along the path given in Figure 4, we have
\begin{equation*}
Y_{+3}(z,t)=Y_{-1}(z,t)\cdot S^{(1)}_{h+}(u)^{-1} \ \text{in the domain}  \ \{ z\in\Sect_0\cap \Sect_l, \ t<0\}.
\end{equation*}
\end{lem}
\begin{proof} It follows from the second identity of part $(c)$ and the first identity of part $(a)$ of Theorem \ref{fact} that
\begin{eqnarray*}
Y_{+3}(z,t)=\widetilde{W}_0(z,zt)F_0^{(1)}(z),
\end{eqnarray*}
and
\begin{eqnarray*}
Y_{-1}(z,t)=\widetilde{W}_0(z,zt)F_l^{(1)}(z)\cdot e^{-\pi\mathi h\delta T^{(1)}},
\end{eqnarray*}
respectively. Here the corresponding index $i$ of the anti-Stokes ray $d_i$ in Theorem \ref{fact} is taken as $i=0$. The rest of the proof is similar to that of Lemma  \ref{conn+3+2}-\ref{conn+1-3}.  \end{proof}

\begin{lem}
After the analytic continuation of $Y_{-1}(z,t)$ along the path given in Figure 5, we have
\begin{eqnarray}
Y_{-1}(z,t)=Y_{-2}(z,t)\cdot S^{(2)}_{h+}(u)^{-1} \hspace{5mm} \text{in }  \ D_{-2}.
\end{eqnarray}
\end{lem}
\begin{proof} It follows from the identities in parts $(a)$ and $(b)$ of Theorem \ref{fact} that
\begin{eqnarray*}
Y_{-1}(z,t)=W_l(z,zt)F_0^{(2)}(zt),
\end{eqnarray*}
and
\begin{eqnarray*}
Y_{-2}(z,t)=W_l(z,zt)F_l^{(2)}(zt)\cdot e^{-\pi\mathi h\delta T^{(2)}},
\end{eqnarray*}
respectively. Here we remark that the corresponding index $i$ of the anti-Stokes ray $d_i$ in Theorem \ref{fact} is taken as $i=l$. \end{proof}

Similar, 
\begin{lem}
After the analytic continuation of $Y_{-2}(z,t)$ along the path given in Figure 6, we have
\begin{eqnarray}
Y_{-2}(z,t)=Y_{-3}(z,t)\cdot R^{(12)} \hspace{5mm} \text{in }  \ D_{-3}.
\end{eqnarray}
\end{lem}

\begin{pro}\label{end2}
After the analytic continuation $D_{+3}\rightarrow D_{-1}\rightarrow D_{-2}\rightarrow D_{-3}$, we have
\begin{equation}\label{RLLid2}
Y_{+3}(z,t)=Y_{-3}(z,t)\cdot R^{(12)}S^{(2)}_{h+}(u)^{-1}S^{(1)}_{h+}(u)^{-1}, \ \text{  for  } (z,t)\in D_{-3}.
\end{equation}
\end{pro}
Since the analytic continuation in the Proposition \ref{end1} and \ref{end2} are along two paths which are homotopic, the comparison of the identities \eqref{RLLid1} and \eqref{RLLid2} proves the identity \eqref{RLL1} in Theorem \ref{thmRLL}. 

The remaining part of Theorem \ref{thmRLL} can be proved in a similar way, i.e., by the comparison of the connection formula of $Y_{\pm k}$ along various homotopic paths, which we simply skip here.
%, we observe the following fact: if we take the analytic continuation of $Y_{\pm k}$ along a simple positive path, instead of the negative path, in the complex $t$ plane, then the upper triangular Stokes matrix $S_{h+}$ will be replaced by the lower one $S_{h-}$ in the connection formula. For example, let us take the analytic continuation of $Y_{+3}(z,t)$ where $t<0$ transports to $t>1$ along a simple positive path in the complex $t$ plane, as shown in Figure 7. 

\section{Representation of quantum groups arising from the Stokes matrices}\label{tworealization}
\subsection{The Faddeev-Reshetikhin-Takhtajan (FRT) realization of quantum groups}
Let us recall the FRT realization \cite{FRT} of the quantized universal enveloping algebra $U(R)$ by means of solutions of the Yang-Baxter equation. The algebra $U(R)$ is generated by elements $l_{ij}^{(+)}$, $l_{ji}^{(-)}$, $1\le i\le j\le n$: set $L_{\pm}=\sum_{i,j}l_{ij}^{(\pm)}\otimes E_{ij}\in U(R)\otimes {\rm End}(\mathbb{C}^n)$, with $l_{ij}^{(+)}=l_{ji}^{(+)}=0$ for $1\le j<i \le n$, then the defining relations are given in matrix form
\begin{eqnarray}\label{lRLL1}
&&R^{12}L^{(1)}_\pm L^{(2)}_\pm=L^{(2)}_\pm L^{(1)}_\pm R^{12},\\ \label{lRLL2}
&&R^{12}L^{(1)}_+ L^{(2)}_-=L^{(2)}_+ L^{(1)}_- R^{12},
\end{eqnarray}
and $l_{ii}^{(+)}l_{ii}^{(-)}=l_{ii}^{(-)}l_{ii}^{(+)}=1,$ for $i=1, ..., n.$

Following \cite{FRT}, the two different realization $U(R)$ and $U_q(\frak{gl}_n)$ (given in Theorem \ref{introthm1}) of quantum groups are isomorphic. In particular, following \cite[Theorem 2.1]{DF}, there is an explicit isomorphism $I: U_q({\frak {gl}}_n)\cong U(R)$, under which
\begin{equation}\label{twoiso}
l_{ii}^{(+)}=q^{-h_i}, \hspace{3mm} l_{i,i+1}^{(+)}=-(q-q^{-1})q^{-h_i}e_i, \hspace{3mm} l_{i+1,i}^{(-)}=(q-q^{-1})f_i q^{h_i}.
\end{equation}

\subsection{Proof of Theorem \ref{introthm1}}

Now we are ready to prove Theorem \ref{introthm1}, i.e., 
\begin{thm}\label{thm1}
For any fixed $h\notin\mathbb{Q}$ and $u\in\h_{\rm reg}$, the map (with $q=e^{\pi\mathi h}$)
\begin{equation}
\begin{split}
\mathcal{S}_q(u): U_q(\frak{gl}_n)&\rightarrow {\rm End}(L(\lambda))~;\\
e_i&\mapsto \frac{S_{h+}(u)_{i,i}^{-1}\cdot S_{h+}(u)_{i,i+1}}{q^{-1}-q}, \\
f_i&\mapsto \frac{S_{h-}(u)_{i+1,i}\cdot S_{h-}(u)_{i,i}^{-1}}{q-q^{-1}}, \\
q^{h_i}&\mapsto S_{h+}(u)_{i,i}
\end{split}
\end{equation}
defines a representation of the Drinfeld-Jimbo quantum group $U_q(\frak{gl}_n)$ on the vector space $L(\lambda)$. 
\end{thm}
\begin{proof}
Theorem \ref{thmRLL} states that the Stokes matrices $S_{h\pm}(u)$ satisfy the RLL relation \eqref{lRLL1}-\eqref{lRLL2}. Thus, as a consequence of the explicit isomorphism $I: U_q({\frak {gl}}_n)\cong U(R)$ given in \eqref{twoiso}, we get a representation of $U_q({\frak {gl}}_n)$ from the Stokes matrices. This concludes the proof.
\end{proof}

\section{The stokes phenomenon of the equation \texorpdfstring{\eqref{introqeq}}{d} in the resonant case}\label{seclast}
In the case of $h\in\mathbb{Q}$, the equation \eqref{introqeq} becomes resonant, and uniqueness of its formal solution is not valid. Accordingly, there exists a family of formal/holomorphic solutions, as well as Stokes matrices $S_{h\pm}(u;c)$, depending on a finite set $c$ of parameters. This section solves the recursive relation in the formal solution of equation \eqref{introqeq} for all $h\ne 0$, and discusses the semiclassical limit of the formal solution. It then introduces a pair of distinguished Stokes matrices $S_{h\pm}^o(u)$ among the family, and gives a description of all $S_{h\pm}(u;c)$ via $S_{h\pm}^o(u)$. In the end, it shows that $S_{h\pm}^o(u)$ gives rise to a representation of $U_q(\frak{gl}_n)$ at $q=e^{h/2}$ a root of unity, i.e., gives a proof of Theorem \ref{introthm1re}. 

\subsection{Solving the recursive relation in the formal solution of equation \eqref{introqeq} for all $h\ne 0$}
In the resonant case, i.e., $h\in \mathbb{Q}$, the proof of Proposition \ref{uniformal} fails. Thus, to study the existence of formal solutions, one needs a better understanding of the recursive relation \eqref{simHm}. For this purpose, let us solve explicitly the recursive relation \eqref{simHm} for all $h\ne 0$ (including $h\in \mathbb{Q})$.  

Let us define $T(\zeta):=\zeta {\rm Id}-hT$ for $\zeta$ an indeterminate. Given any elements $a,b$ and $c,d$ of $\{1,...,n\}$, let us introduce the $2\times 2$ quantum minor (formed from rows $a,b$ and columns $c,d$) of $T(\zeta)$,
\begin{equation}\label{qminor}
\Delta^{a,b}_{c,d}(T(\zeta)):=(\delta_{ac}\zeta -h e_{ac})(\delta_{bd}\zeta -h e_{bd})-(\delta_{bc}(\zeta+h) -h e_{bc})(\delta_{ad}(\zeta+h) -h e_{ad}), 
\end{equation}
which is a degree two polynomial in $\zeta$ with coefficients in $U(\frak{gl}_n)$. Here $\delta_{ab}=1$ for $a=b$; and $\delta_{ab}=0$ for $a\ne b$. Suppose $X\in U(\frak{gl}_n)$ is an element such that 
\begin{align*}
    \delta_{ac}X\cdot e_{bd}=e_{bd}\cdot \delta_{ac}X, \hspace{5mm} \delta_{bd}X\cdot e_{ac}=e_{ac}\cdot \delta_{bd}X,\\ \delta_{bc}X\cdot e_{ad}=e_{ad}\cdot \delta_{bc}X, \hspace{5mm} \delta_{ad}X\cdot e_{bc}=e_{bc}\cdot \delta_{ad}X,
\end{align*}
then $\Delta^{i,j}_{k,l}(T(X))$ is defined without ambiguity. For example, let us introduce the special quantum minors that are used in the following proposition: for any index $1\le k,j \le n$ and $k,j\ne l$, by the definition of the quantum minor in \eqref{qminor}, we have 
\[\Delta^{l,k}_{l,j}(T(\zeta))=\left\{\begin{array}{lr}  -(\zeta -h e_{ll})he_{kj}-h e_{kl}he_{lj},  & \text{ if } \ k\ne j \text{ and } k,j\ne l\\ 
    (\zeta+1 -h e_{ll})(\zeta+1-he_{kk})-h e_{kl}he_{lj}, & \text{ if } \  k=j\ne l.
             \end{array}
\right.\] 
Since for any integer $m$, $he_{ll}-m$ commutes with the elements $e_{ll}e_{kj}$ and $e_{kl}e_{lj}$, by replacing $\zeta$ by $he_{ll}-m$ in $ \Delta^{l,k}_{l,j}(T(\zeta))$, we get 
\begin{equation}\label{qkj}
    \Delta^{l,k}_{l,j}(T(he_{ll}-m))=\left\{\begin{array}{lr}  mhe_{kj}-h e_{kl}he_{lj},  & \text{ if } \ k\ne j \text{ and } k,j\ne l\\ 
    -m(he_{ll}-m-he_{kk})-h e_{kl}he_{lj}, & \text{ if } \  k=j\ne l.
             \end{array}
\right.
\end{equation}

\begin{pro}\label{Yangsol}\cite{TX}
For all $h\in \mathbb{C}\setminus\{0\}$, the ordinary differential equation \eqref{introqeq}
has a formal fundamental solution
\begin{equation*}
\widehat{F}(z)=\widehat{H}(z) e^{{huz}}z^{h\delta T}, \ \ \ {\it with} \ \widehat{H}=1+H_1z^{-1}+H_2z^{-1}+\cdot\cdot\cdot, \end{equation*}
where every coefficient $H_{m+1}=\sum_{i,j}(H_{m+1})_{ij}\otimes E_{ij}\in{\rm End}(L(\lambda))\otimes{\rm End}(\mathbb{C}^n)$ is recursively determined by
\begin{eqnarray}\label{Yangian1}
(H_{m+1})_{kl}=\left\{\begin{array}{lr}  \frac{1}{m}\frac{1}{hu_l-hu_k}\Big(\sum\limits_{j\ne l}{\Delta^{l,k}_{l,j}\big(T(he_{ll}-m)\big)}\cdot (H_{m})_{jl}\Big),  & \text{ if } \ k\ne l\\ -\frac{1}{m+1} \Big(\sum\limits_{j \neq k} he_{kj}\cdot  (H_{m+1})_{jk}\Big), & \text{ if } \  k=l.
             \end{array}
\right.
\end{eqnarray}
Furthermore, we have the commutative relations for all integers $m\ge 0$ and $k=1,...,n$,
\begin{equation}\label{commrel}
[e_{ll},H_{m,kl}]=\left\{\begin{array}{lr}  -H_{m,kl},  & \text{ if } \ l\ne k\\ 
    0, & \text{ if } \ l=k.
             \end{array}
\right.
\end{equation}
\end{pro}
\begin{proof}
Let us prove inductively that \eqref{Yangian1} provides a solution of the recursive relation \eqref{simHm} in Proposition \ref{uniformal} (thus a formal solution of the equation), and \eqref{commrel} are true. First for $m=0, 1$ we check that
\begin{equation}
    H_{0,kl}=\left\{\begin{array}{lr}  0,  & \text{ if } \ k\ne l\\ 
    1, & \text{ if } \  k=l
             \end{array}
\right. \text{   and   }  H_{1,kl}=\left\{\begin{array}{lr}  \frac{e_{kl}}{u_l-u_k},  & \text{ if } \ k\ne l\\ 
        \sum_{j=1, j\ne k}^n\frac{he_{kj}e_{jk}}{u_k-u_l}, & \text{ if } \  k=l
             \end{array}
\right. 
\end{equation} 
give a solution of the recursive relations \eqref{Hknel} and \eqref{Hkel} for $m=0$. Here note that $H_0$ is just the initial condition. In the meanwhile, one checks that they satisfy \eqref{Yangian1} and \eqref{commrel} for $m=0$ by hand. 

Now given a sequence of solutions $H_i$ for $i=0,1,...,m$ of the recursive relation \eqref{Hknel} and \eqref{Hkel} such that \eqref{Yangian1} and \eqref{commrel} hold, let us construct $H_{m+1}$ explicitly.

{\bf Case I:} $1\le k\ne l\le n$. Recall that in this case, $H_{m+1, kl}$ is uniquely determined by the recursive relation \eqref{Hknel}. Thus let us define
\begin{align}  
(u_l-u_k) H_{m+1, kl}:&=\frac{m}{h}H_{m,kl}+\sum_{j\ne l} e_{kj} H_{m, jl}+e_{kl}H_{m, ll} - H_{m, kl} e_{ll} \\
&=\frac{m}{h}H_{m,kl}+\sum_{j\ne l} e_{kj} H_{m, jl}+e_{kl} \cdot \Big(-\frac{h}{m}\sum\limits_{j \neq l} e_{lj}\cdot  H_{m,jl}\Big)-H_{m,kl} e_{ll} \\ \label{commHe}
&=\frac{m}{h}+\sum_{j\ne l} e_{kj} H_{m, jl}-\frac{h}{m}\sum\limits_{j \neq l} e_{kl}e_{lj}\cdot  H_{m,jl}-e_{ll} H_{m,kl}-H_{m,kl}\\
&= \frac{1}{mh}\Big(\sum\limits_{j\ne l}{\Delta^{l,k}_{l,j}\big(T(he_{ll}-m)\big)}\cdot (H_{m})_{jl}\Big)
\end{align}
Here in the second identity, by induction we replace $H_{m,ll}$ by the right hand side of \eqref{Yangian1} (with $k=l$); in the third identity, we use the induction assumption \eqref{commrel} to replace $H_{m,kl}e_{ll}$ by $e_{ll}H_{m,kl}+H_{m,kl}$; in the last identity, we just use the definition \eqref{qkj} of the quantum minors. 

On the other hand, taking the right hand side of the identity \eqref{commHe} as the definition of $H_{m+1,kl}$, then using the relations $[he_{ll}, e_{kj}]=0$ for $k,j\ne l$, $[e_{ll},e_{kl}e_{lj}]=0$, we see that 
\begin{equation}\label{llkl}
    [e_{ll},H_{m+1,kl}]=0, \ \text{ for } k\ne l.
\end{equation}

Here we stress that the extra term $-H_{m,kl}$, coming from the commutative relation \eqref{commrel}, just cancels with the (Capelli) shift, i.e., the number $1$ appearing in $\zeta+1$, in the definition of the quantum minor.

{\bf Case II:} $1\le k=l\le n$. Note that $H_{m+1, kl}$ for $k\ne l$ are known as in the Case I. Now let us construct $H_{m+1, kk}$. By \eqref{llkl} we have $[e_{kk}, e_{kj}\cdot  H_{m+1,jk}]=0$ for any $j\ne k$. Therefore, we get
\begin{equation}\label{kkjk}
\Big[e_{kk},-\frac{h}{m+1} \sum\limits_{j \neq k} e_{kj}\cdot  H_{m+1,jk}\Big]=0.
\end{equation}
Thus we see that the element
\[H_{m+1, kk}:=-\frac{h}{m+1} \Big(\sum\limits_{j \neq k} e_{kj}\cdot  H_{m+1,jk}\Big)\in {\rm End}(L(\lambda)) \]
solves the recursion relation \eqref{Hkel}, i.e.,
\[
0=\sum_{j\ne k} e_{kj} H_{m+1, jk}+\frac{m+1}{h} H_{m+1, kk}+[e_{kk}, H_{m+1, kk}].\]
And by \eqref{kkjk}, we have $[e_{kk}, H_{m+1, kk}]=0$. 

Combing the Case I and Case II, we have constructed a solution $H_{m+1}$ of the recursive relation \eqref{Hknel} and \eqref{Hkel} such that \eqref{Yangian1} and \eqref{commrel} hold for $m+1$. By induction, this concludes the proof. 
\end{proof}

In particular, this proposition states that the entries of the coefficients $H_m$ are polynomials of degree $m$ in $h$. Thus, the formal solution $\hat{F}(z)$ defined in Proposition \ref{uniformal} extends to $h\in\mathbb{Q}$. To distinguish between the resonant and nonresonant cases, we denote by 
\begin{eqnarray}
    \label{Fo}\widehat{F}^o(z)=\widehat{H}^o(z) e^{{huz}}z^{h\delta T}
\end{eqnarray}
the formal solution in Proposition \ref{Yangsol} as $h\in\mathbb{Q}$.
\begin{rmk}
The coefficient matrix of the recursive relation in \eqref{Yangian1} actually encodes the defining relation of $Y(\frak{gl}_{n-1})$. In \cite{TX}, we prove for any $u$ (not necessary regular), the recursive relation is equivalent to a difference equation
\[X(z+1)=L(z) X(z),\]
where $L(z)$ satisfies the defining relation of Yangian $Y({\rm gl}_m)$ for an integer $m$ determined by $u$. We also generalize the result to the analogous equation \eqref{introqeq} for classical Lie algebras and twisted Yangians. We refer the reader to \cite[Section 2.5]{Molev} and \cite[Section 3.3]{Molev} for a general theory on the quantum minors and Yangians.
\end{rmk}
\subsubsection*{Semiclassical limit of the formal power series}
Note that $hT=(he_{ij})$ is a collection of noncommutative variables $he_{ij}$ in $U(\frak{gl}_n)$, while $A=(a_{ij})$ is a collection of commutative variables $a_{ij}\in Sym(\frak{gl}_n)$. However, since $[he_{ij},he_{kl}]=h (\delta_{jk}he_{il}-\delta_{li} he_{kj})$, as $h\rightarrow 0$ any power series in $he_{ij}$ corresponds to a power series of the same form but in the commutative variables $a_{ij}$. That is under the isomorphism $U(h\frak{gl}_n)/hU(h\frak{gl}_n)\cong Sym(\frak{gl}_n)$, any word \[he_{i_1j_1}he_{i_2j_2}\cdots he_{i_kj_k}\in U(h\frak{gl}_n), \ \ i_1, j_1, ..., i_k,j_k\in \{1,...,n\} \] becomes the polynomial \[a_{i_1j_1}a_{i_2j_2}\cdots a_{i_kj_k}\in Sym(\frak{gl}_n),\ \ i_1, j_1, ..., i_k,j_k\in \{1,...,n\} .\]

Therefore, the $h\rightarrow 0$ limit of the formal solution of \eqref{introqeq} will become the formal solution of \eqref{introeq}. To be more precise, in \cite{TXu}, we prove
\begin{pro}\cite{TXu}\label{sclrmatrix}
For any $u\in \h_{\rm reg}$ and $A\in\frak{gl}_n$, the ordinary differential equation \eqref{introeq} has a formal fundamental solution
$\widehat{f}(z)=\widehat{h}(z) e^{{huz}}z^{h\delta A}$, with $\widehat{h}=1+h_1z^{-1}+h_2z^{-1}+\cdot\cdot\cdot$,
where the entries $h_{m+1,ij}\in \mathbb{C}$ of every coefficient matrix $h_{m+1}=(h_{m+1,ij})\in {\rm End}(\mathbb{C}^n)$ is recursively determined by
\begin{equation}\label{sclhm}
h_{m+1,ik}=\left\{\begin{array}{lr}  \frac{1}{m}\frac{1}{u_k-u_i}\Big(\sum\limits_{j\ne k}{\Delta^{k,i}_{k,j}\big((a_{kk}-m){\rm Id}-A\big)}\cdot h_{m,jk}\Big), & \text{ if } \ i\ne k  \\ -\frac{1}{m+1} \sum\limits_{j \neq k} a_{kj}\cdot  h_{m+1,jk}, & \text{ if } \ i=k,
             \end{array}
\right.
\end{equation}
where $\Delta^{k,i}_{k,j}\big((a_{kk}-m){\rm Id}-A\big)$ is the ordinary $2\times 2$ minor formed by the $k,i$ rows and $k,j$ columns of the $n\times n$ matrix $(a_{kk}-m){\rm Id}-A$. Here we denote by $(a_{kk}-m){\rm Id}$ the scalar matrix with all diagonal elements equal to $a_{kk}-m\in\mathbb{C}$.    
\end{pro}
Proposition \ref{sclrmatrix} can be proved in the same way as Proposition \ref{Yangsol}. Now by \eqref{Yangian1} the entry $H_{m,ik}$ of $H_{m}$ is an element of filtered degree $2m$ in $U(h\frak{gl}_n)$, while by \eqref{sclhm} the entry $h_{m,ik}$ of $h_m$ is a degree $2m$ polynomial in $Sym(\frak{gl}_n)$. Following Proposition \ref{sclrmatrix} and Proposition \ref{Yangsol}, the classical coefficient in \eqref{sclhm} has the same form of the quantum coefficient given in \eqref{Yangian1}, provided replacing the ordinary minors by quantum minors. It follows that under the isomorphism $U(h\frak{gl}_n)/hU(h\frak{gl}_n)\cong Sym(\frak{gl}_n)$, $H_{m,ik}$ becomes $h_{m,ik}$. That is the precise meaning of the fact that the $h\rightarrow 0$ limit of the formal solution of \eqref{introqeq} becomes the formal solution of \eqref{introeq}.

Taking the semiclassical limit as $h\rightarrow 0$ is also compatible with the Borel resummation, and thus the semiclassical limit of the Stokes matrices of \eqref{introqeq} becomes the Stokes matrices of \eqref{introeq}. See \cite{TLX} for a similar statement in a different setting.

\subsection{The family of formal fundamental solutions}\label{famforsol}
Recall the recursion relation \eqref{Hkel}, 
\begin{equation}\label{Hk}
0=\sum_{j\ne k} e_{kj} H_{m, jk}+\frac{m}{h} H_{m, kk}+[e_{kk}, H_{m, kk}].
\end{equation}
Knowing the off-diagonal elements $H_{m, jk}$ for $j\ne k$, the equation \eqref{Hk} is used to find $H_{m,kk}$. Assume $H_{m,kk}$ solves \eqref{Hk}. We see that if an element $U\in{\rm End}(L(\lambda))$ satisfies $mU-h[e_{kk}, U]=0$ (since $h\in\mathbb{Q}$ such an element may exists), then $H_{m,kk}+U$ also solves \eqref{Hk}. It reflects the non-uniqueness of the formal solution, and motivates the following parameterization of the family of formal solutions of \eqref{introqeq}. 

For any positive integer $m$, let $I_{m}\subset \{1,...,n\}$ denote the set of all $k$ such that $m{\rm Id}-h{\rm ad}_{e_{kk}}$ is not invertible on ${\rm End}(L(\lambda))$. Then let $U_{mk,s}\in{\rm End}(L(\lambda))$, for $s=1,2,....,d_{mk}$, be a set of linearly independent elements satisfying 
\begin{eqnarray}\label{mkk}
mU_{mk,s}-h[e_{kk}, U_{mk,s}]=0.
\end{eqnarray}
Let $I$ be set of the positive integers $m$ such that the above set $I_m$ is not empty. Associated to a family $c=\{c_{ik,s}\}$ of complex parameters for $m\in I$, $k\in I_m$ and $s=1,...,d_{mk}$, we introduce a Laurent polynomial valued in ${\rm End}(L(\lambda))\otimes {\rm End}(\mathbb{C}^n)$
\begin{equation}\label{polyU}
U(z;c)=1+\sum_{m\in I}\sum_{k\in I_{m}}\sum_{s=1}^{d_{mk}}c_{mk,s} U_{mk,s}\otimes E_{kk} z^{-m}.
\end{equation}

\begin{pro}\label{nuniformal}
For $h\in\mathbb{Q}$, the ordinary differential equation \eqref{introqeq} has a family of formal solutions taking the form \begin{eqnarray}\label{familyfsol}
\widehat{F}(z;c)=\widehat{H}^o(z)\cdot U(z;c) e^{{huz}}z^{h{\delta T}}. \end{eqnarray}
Here $\widehat{H}^o(z)$ is the formal series from \eqref{Fo}. Furthermore,
\begin{eqnarray}\label{Dc}
U(z;h, c) e^{{zu}}z^{h{\delta T}}=e^{{zu}}z^{h{\delta T}} D(c)
\end{eqnarray}
where $D(c)=\sum_k D_{kk}\otimes E_{kk}$ is the diagonal constant given by
\[ D(c)=1+\sum_{m\in I}\sum_{k\in I_{m}}\sum_{s=1}^{d_{mk}}c_{mk,s} U_{mk,s}\otimes E_{kk}\in{\rm End}(L(\lambda))\otimes {\rm End}(\mathbb{C}^n). \]
\end{pro}
\begin{proof} To verify the function in \eqref{familyfsol} is a solution, we only need to check
\begin{eqnarray}
\frac{dU(z;c)}{dz}+\Big[U(z;c), \frac{h\delta T}{z}\Big]=0.
\end{eqnarray}
Recall that $\delta T=\sum_k e_{kk}\otimes E_{kk}$, then the above identity follows from \eqref{mkk}. In the meanwhile, a direct computation using the identity \eqref{mkk} shows that \eqref{Dc} holds.\end{proof}

\subsection{A family of canonical solutions with a prescribed asymptotics and Stokes matrices}\label{famcanonicalF}
As in the case of $h\notin\mathbb{Q}$, the Borel resummation of $\hat{H}^o(z)$ can be used to construct the unique holomorphic function $H^o_i(z)$ with the asymptotics $\hat{H}^o(z)$ as $z\rightarrow\infty$  within each Stokes supersector $\Sect_i$.

\begin{thm}\label{famthmcanonicalsol}
For any choice of the parameter set $c$, the function $F_{i}(z;c)=H^o_i (z) U(z;c)\cdot  e^{huz}z^{h{\delta T}}$ is the unique (therefore called canonical) holomorphic solution of \eqref{introqeq} defined on $\Sect_i$, such that $F_{i}e^{-huz}\cdot z^{-h\delta T}\sim \hat{H}_i^o\cdot U(z;c)$ within $\Sect_i$.
\end{thm}

\begin{defi}\label{famdefiStokes}
Given $h\in\mathbb{Q}$, $u\in \h_{\rm reg}$ and any choice of the parameter set $c$, the {\it Stokes matrices} of the equation \eqref{introqeq} (with respect
to the chosen sector $\Sect_0$ and the branch of ${\rm log}(z)$) are the unique matrices such that:
\begin{itemize}
    \item  if $F_0$ is continued in a positive sense to $\Sect_l$ then $F_l(z;c)=F_0 (z;c)\cdot  S_{h+}(u;c)e^{\pi\mathi h\delta T}$, and 
    \item if $F_l$
is continued in a positive sense to $\Sect_0$ then $F_0(ze^{-2\pi \mathi };c)=F_l(z;c)\cdot  e^{-\pi\mathi h\delta T}  S_{h-}(u;c)^{-1}$.
\end{itemize}
\end{defi} 

\subsection{Comparison of Stokes matrices in the resonant and nonresonant cases}
On the one hand, for $h\notin \mathbb{Q}$, $u\in\h_{\rm reg}$ and a representation $L(\lambda)$ of $\frak{gl}_n$, there is only one pair of Stokes matrices $S_{h\pm}(u)$, based on the uniqueness of the formal solution of \eqref{introqeq}. Then following Proposition \ref{Yangsol}, for any fixed $u$ the unique formal solution (therefore its Borel resummation and the corresponding Stokes matrices) extends from $h\notin\mathbb{Q}$ to rational numbers. For $h\in\mathbb{Q}$ the resulting Stokes matrices from the natural extension are denoted by $S_{h\pm}^o(u)$.

On the other hand, for the same $h\in\mathbb{Q}$, $u$ and $L(\lambda)$, there is a family of Stokes matrices $S_{h\pm}(u;c)$ of \eqref{introqeq} parameterized by the set $c$, where $c$ depends on the chosen $h$ and $L(\lambda)$. If all the parameters $c_{mk,s}$ in \eqref{polyU} of the set $c$ are chosen to be zero, then the corresponding formal solution coincides with the one given in Proposition \ref{Yangsol}. Therefore, the corresponding Stokes matrices $S_{h\pm}(u;c=0)$ coincide with the natural extension $S_{h\pm}^o(u)$ from the resonant cases.

Following the definition of Stokes matrices and the identity
\[F_{i}(z;c)=H^o_i(z) U(z;c)\cdot  e^{huz}z^{h{\delta T}}=H^o_i (z)\cdot  e^{huz}z^{h{\delta T}}D(c)\]
on each $\Sect_i$, we have
\begin{cor}\label{reducec}
For $h\in\mathbb{Q}$, $u\in\h_{\rm reg}$ and any choice of parameter set $c$, the Stokes matrices 
\begin{eqnarray}
S_i(u;c)= D(c)\cdot S_i^o(u)\cdot D(c)^{-1},
\end{eqnarray}
where $D(c)$ is the diagonal element given in Proposition \ref{nuniformal}.
\end{cor}

Corollary \ref{reducec} gives a description of the family of Stokes matrices of \eqref{introqeq} in the case of $h$ a rational number, and reduces the family to the study of $S_{h\pm}^o(u)$.

\subsection{The $2\times 2$ example}
Let us first assume $h\notin \mathbb{Q}$ and consider the rank two case. In this case, up to some substitution, the parameters $u_1$ and $u_2$ can be changed to $0$ and $1$. That is
\begin{equation}\label{22}
\frac{dF}{dz}=h\left(\left(\begin{array}{cc}
    0 & 0  \\
    0 & 1
  \end{array}\right)
+{\left(
  \begin{array}{cc}
    e_{11} & e_{12}  \\
    e_{21} & e_{22}
  \end{array}
\right)} z^{-1} \right)\cdot F.
\quad
\end{equation} 
After a coordinate change, 

Given any finite dimensional representation $L(\lambda)$ of $\frak{gl}_n$ with the highest weight $\lambda=(\lambda^{(2)}_1, \lambda^{(2)}_2)$, where $\lambda^{(2)}_1-\lambda^{(2)}_2\in\mathbb{Z}_{\ge 0}$. Let $\{\xi_\Lambda\}$ be a Gelfand-Tsetlin basis of $L(\lambda)$, parameterized by the Gelfand-Tsetlin patterns $\Lambda$. Such a pattern is a collection of numbers $\{\lambda^{(1)}_1(\Lambda), \lambda^{(2)}_1(\Lambda), \lambda^{(2)}_2(\Lambda)\}$ with the fixed $\lambda^{(2)}_1(\Lambda)=\lambda^{(2)}_1$ and $\lambda^{(2)}_2(\Lambda)=\lambda^{(2)}_2$ satisfying the interlacing conditions
\begin{eqnarray}
\lambda_1^{(2)}(\Lambda)-\lambda^{(1)}_1(\Lambda)\in\mathbb{Z}_{\ge 0}, \hspace{5mm} \lambda^{(1)}_1(\Lambda)-\lambda^{(2)}_{2}(\Lambda)\in\mathbb{Z}_{\ge 0}.
\end{eqnarray}

Let us introduce the endomorphism on ${\rm End}(L(\lambda))$,
\begin{equation}\label{opl}
\ell^{(j)}_i(\xi_\Lambda)=\lambda^{(j)}_i(\Lambda)-i+1, \ \ \ \text{ for } 1\le i\le j\le 2.
\end{equation}
If ${\rm dim}(L(\lambda))=m$, then in terms of the Gelfand-Tsetlin basis, each $\ell^{(j)}_i$ is just a diagonal $m\times m$ matrix. Set $h$ a nonezero complex number. If (under the Gelfand-Tsetlin basis) the diagonal matrix \[1+h \ell^{(j_1)}_{i_1}-h \ell^{(j_2)}_{i_2}={\rm diag}(a_1,...,a_m) \ \text{ and } \ h\ell^{(j_3)}_{i_3}-h\ell^{(j_3)}_{i_3}={\rm diag}(b_1,...,b_m)\] for any index $1\le i_s\le j_s\le 2$ and $s=1,2,3, 4$, then we define the diagonal matrix (element in ${\rm End}(L(\lambda))$)
\[{}_{1}F_{1}\big(1+h\ell^{(j_1)}_{i_1}-h\ell^{(j_2)}_{i_2},hl^{(j_3)}_{i_3}-h\ell^{(j_3)}_{i_3};z\big):={\rm diag}\big({}_{1}F_{1}(a_1,b_1;z),...,{}_{1}F_{1}(a_m,b_m;z)\big).\]
Here for any parameters $\alpha\in \mathbb{C}$ and $\beta\in \mathbb{C}\backslash\{0,-1,-2,\cdots\}$,
the function 
\begin{equation}\label{kummer}
    {}_{1}F_{1}(\alpha,\beta;z)=\sum_{k=0}^{\infty}\frac{(\alpha)_k}{(\beta)_k}\frac{z^k}{k!}
\end{equation}
is the confluent hypergeometric function, where $(\alpha)_0=1$ and $(\alpha)_k=\alpha (\alpha+1) \cdots (\alpha+k-1)$, for $k\geq 1$. 

\begin{rmk}
In the $2\times 2$ case, the notations can be simplified. But we still use the convention from the Gelfand-Tsetlin basis, because the following proposition has a generalization to all $n$.
\end{rmk}
\begin{pro}\cite{Xu}\label{2by2}
The equation \eqref{22}
has a solution 
\begin{equation}\label{quantum diagFz}
    F(z)= \left(\begin{array}{cc}
   e_{12}\cdot Y_{11} & e_{12}\cdot Y_{12}  \\
    Y_{21} & Y_{22}
  \end{array}\right) \cdot z^{\frac{h}{2\pi \mathrm{i}}T_{d}},
\end{equation}
where $T_{d}=\operatorname{diag}(\ell^{(2)}_1+1,\ell^{(2)}_2+1)\in {\rm End}(L(\lambda))\otimes {\rm End}(\mathbb{C}^2)$, 
and the elements in ${\rm End}(L(\lambda))$ are
\begin{align*}
&Y(z)_{11}=\frac{1}{\ell^{(2)}_1-\ell^{(1)}_1} \cdot{ }_{1} F_{1}(\alpha_{11}; \beta_{11} ;hz
),\\
&Y(z)_{12}=\frac{1}{\ell^{(2)}_2-\ell^{(1)}_1} \cdot{ }_{1} F_{1}(\alpha_{12}; \beta_{12} ;hz
),\\
	&Y(z)_{21}= { }_{1} F_{1}(\alpha_{21}; \beta_{21};hz),
 \\
	&Y(z)_{22}= { }_{1} F_{1}(\alpha_{22}; \beta_{22};hz),
\end{align*}
with the diagonal actions in ${\rm End}(L(\lambda))$ given by
\begin{align*}
&\alpha_{11} =h( \ell_{1}^{(2)}-\ell_{1}^{(1)}), &\beta_{11}=1+h(\ell_{1}^{(2)}-\ell_{2}^{(2)}),\\
&\alpha_{12} =h( \ell_{2}^{(2)}-\ell_{1}^{(1)}), &\beta_{12}=1+h(\ell_{2}^{(2)}-\ell_{1}^{(2)}),\\
&\alpha_{21} =1+h( \ell_{1}^{(2)}-\ell_{1}^{(1)}),&\beta_{21}=1+h(\ell_{2}^{(1)}-\ell_{2}^{(2)}),\\
&\alpha_{22} =1+h( \ell_{2}^{(2)}-\ell_{1}^{(1)}),  &\beta_{22}=1+h(\ell_{2}^{(2)}-\ell_{1}^{(2)}).
\end{align*}
\end{pro}
Besides, let $F_0(z)$ and $F_1(z)$ be the canonical solutions in ${\rm Sect}_0$ and ${\rm Sect}_{1}$ given as in Section \ref{canonicalF} (in this case there are only two secotrs). These solutions are related by $2\times 2$ transition matrices $U_0$ and $U_{1}\in{\rm End}(L(\lambda))\otimes {\rm End}(\mathbb{C}^2)$ such that 
\begin{align}\label{U0}
&F(z) =F_0(z) \cdot U_{0},\quad z\in \operatorname{Sect}_0,\\ \label{U1}
&F(z)=F_{1}(z) \cdot U_{1},\quad z\in \operatorname{Sect}_{l},
\end{align}
By the definition, we know that the $2\times 2$ Stokes matrices of \eqref{22} is computed by
\begin{equation}\label{S2eq}
     S_{h+}(u) \cdot U_0=\mathrm{e}^{\frac{\delta T}{2}} U_{1}.
\end{equation} 
Using the following known asymptotic behavior of confluent hypergeometric function ${}_{1}F_{1}(\alpha,\beta;z)$ with complex parameters $\alpha, \beta$, see. e.g., \cite{OLBC}, one can compute explicitly the asymptotic expansion of $F(z)$ as $z\rightarrow \infty$,
\begin{align*}
    {}_{1}F_{1}(\alpha,\beta;z)\sim&\frac{\Gamma(\beta)}{\Gamma(\beta-\alpha)}\frac{1}{(z \mathrm{e}^{-\pi \mathi })^{{\alpha}}}+\frac{\Gamma(\beta)}{\Gamma(\alpha)}\operatorname{e}^zz^{(\alpha-\beta)}, \text{ as } z\rightarrow \infty \text{ within } S(0,\pi);\\
    {}_{1}F_{1}(\alpha,\beta;z)\sim&\frac{\Gamma(\beta)}{\Gamma(\beta-\alpha)}\frac{1}{(z \mathrm{e}^{\pi \mathi })^{{\alpha}}}+\frac{\Gamma(\beta)}{\Gamma(\alpha)}\operatorname{e}^zz^{(\alpha-\beta)}, \text{ as } z\rightarrow \infty \text{ within } S(-\pi,0).
\end{align*}
Comparing the asymptotic expansion of $F(z)$ with the prescribed asymptotics of $F_0(z)$ and $F_1(z)$ in \eqref{U0} and \eqref{U1}, we find 
  \begin{align*}
U_0&=\left(\begin{array}{cc}
   \frac{2\pi\mathi e_{12}\Gamma(\beta_{11})}{\big(\ell^{(2)}_1-\ell_1^{(1)}\big)\Gamma(\beta_{11}-\alpha_{11})\big(h\mathrm{e}^{-\frac{\pi \mathi }{2}}\big)^{\alpha_{11}}} & \frac{2\pi\mathi e_{12}\Gamma(\beta_{12})}{\big(\ell^{(2)}_1-\ell_1^{(1)}\big)\Gamma(\beta_{12}-\alpha_{12})\big(h\mathrm{e}^{-\frac{\pi \mathi }{2}}\big)^{\alpha_{12}}}  \\
   \frac{\Gamma(\beta_{21})}{\Gamma(\alpha_{21})\big(h\mathrm{e}^{\frac{\pi \mathi }{2}}\big)^{\alpha_{21}}} & \frac{\Gamma(\beta_{22})}{\Gamma(\alpha_{22})(h\mathrm{e}^{\frac{\pi \mathi }{2}})^{\alpha_{22}}}
  \end{array}\right)\\
U_1&=\left(\begin{array}{cc}
   \frac{2\pi\mathi e_{12}\Gamma(\beta_{11})}{\big(\ell^{(2)}_1-\ell_1^{(1)}\big)\Gamma(\beta_{11}-\alpha_{11})\big(h\mathrm{e}^{\frac{3\pi \mathi }{2}}\big)^{\alpha_{11}}} & \frac{2\pi\mathi e_{12}\Gamma(\beta_{12})}{\big(\ell^{(2)}_1-\ell_1^{(1)}\big)\Gamma(\beta_{12}-\alpha_{12})\big(h\mathrm{e}^{\frac{3\pi \mathi }{2}}\big)^{\alpha_{12}}}  \\
   \frac{\Gamma(\beta_{21})}{\Gamma(\alpha_{21})\big(h\mathrm{e}^{\frac{\pi \mathi }{2}}\big)^{\alpha_{21}}} & \frac{\Gamma(\beta_{22})}{\Gamma(\alpha_{22})\big(h\mathrm{e}^{\frac{\pi \mathi }{2}}\big)^{\alpha_{22}}}
  \end{array}\right).
\end{align*}
Here again we use the matrix function convention that if $\alpha_{ij}$ or $\beta_{ij}={\rm diag}(d_1,...,d_m)$ in ${\rm End}(L(\lambda))$ under the Gelfand-Tsetlin basis, then $\Gamma(\alpha_{ij}):={\rm diag}(\Gamma(d_1),...,\Gamma(d_m))$ and $e^{\alpha_{ij}}:={\rm diag}(e^{d_1},...,e^{d_m})$.

Using the identity $(1+e_{12})e_{11}=e_{11}e_{12}$ and \eqref{S2eq}, we get
\begin{equation}\label{Stokes22}
S_{h+}(u)=\left(\begin{array}{cc}
   e^{\frac{he_{11}}{2}} & \frac{- 2\pi \mathi  (-h)^{\frac{h(e_{11}-e_{22}-1)}{2\pi \mathi }}  \mathrm{e}^{\frac{he_{11}}{2}}}{\Gamma(1+h(\ell^{(1)}_1-\ell^{(2)}_1-1))\Gamma(1+h(\ell^{(1)}_1-\ell^{(2)}_2-1))}he_{12}.  \\
  0 & e^{\frac{he_{22}}{2}}
  \end{array}\right).
\end{equation}
In the above computation, we make the assumption that $h\notin \mathbb{Q}$. It is because if $h$ is a negative rational number, the eigenvalues of $\beta_{11}=1+h(\ell^{(2)}_1-\ell^{(2)}_2)\in {\rm End}(L(\lambda))$ can be a negative integer. Thus, there exists a nonzero vector $\xi$ in $L(\lambda)$ such that the action of ${ }_{1} F_{1}(\alpha_{12}, \beta_{12};hz)\in {\rm End}(L(\lambda))$ on $\xi$ is not well-defined (in the definition \eqref{kummer}, the complex parameter $\beta$ can not be a non-positive integer). That is the solution $F(z)$ given in Proposition \eqref{2by2} may not be well-defined.

However, after we have computed the Stokes matrices in case that $h\notin \mathbb{Q}$, we can treat the leftover cases by means of continuation of $h$ 
into $\mathbb{Q}$. Thus the formula \eqref{Stokes22} is valid for all $h$.
Using the formula, we get

\vspace{2mm}
{\bf Proof of Proposition \ref{prop2}.}
Assume that $h=-\frac{m}{p}$ is a negative rational number, if $\Lambda^o$ is the pattern determined by $ 
\lambda^{(1)}_1(\Lambda^o)-\lambda^{(2)}_{2}(\Lambda^o)=p-1$, then $e_{12}(\xi_{\Lambda^0})$, the action of the raising operator $e_{12}$ on the basis element $\xi_{\Lambda^0}$, is a scalar multiple of another basis element $\xi_{\Lambda'}$ with the pattern $\Lambda'$ satisfying $\lambda^{(1)}_1(\Lambda')-\lambda^{(2)}_{2}(\Lambda')=p$. Let us write $ e_{12}(\xi_{\Lambda^0})=c \xi_{\Lambda'}$. Thus, the action of the $12$-entry $S_{h+}(u)_{12}\in {\rm End}(L(\lambda))$ on $\xi_{\Lambda^0}$ is
\begin{align}
S_{h+}(u)_{12}(\xi_{\Lambda^0})&=  \frac{-  2\pi \mathi (-h)^{\frac{h}{2\pi \mathi }(e_{11}-e_{22}-1)}\mathrm{e}^{\frac{e_{11}}{2}}}{\Gamma(1+h(\ell^{(1)}_1-\ell^{(2)}_1-1))\Gamma(1+h(\ell^{(1)}_1-\ell^{(2)}_2-1))}he_{12}(\xi_{\Lambda^0})\\
&= \frac{-  2\pi \mathi (-h)^{\frac{h}{2\pi \mathi }(e_{11}-e_{22}-1)}}{\Gamma(1+h(\ell^{(1)}_1-\ell^{(2)}_1-1))\Gamma(1+h(\ell^{(1)}_1-\ell^{(2)}_2-1))}hc\xi_{\Lambda'}\\
&=0.
\end{align}
In the last identity, we use the definition \eqref{opl} of the diagonal elements $\ell^{(1)}_1$ and $\ell^{(2)}_2$ to get
\[\frac{1}{\Gamma(1+h(\ell^{(1)}_1-\ell^{(2)}_2-1))}(\xi_{\Lambda'})=\frac{1}{\Gamma(1+h(\lambda^{(1)}_1(\Lambda')-\lambda^{(2)}_2(\Lambda')))}\times \xi_{\Lambda'}=0.\]
Here we use the fact that the non-positive integer $1+h(\lambda^{(1)}_1(\Lambda')-\lambda^{(2)}_2(\Lambda'))=1-m$ is a zero of $1/\Gamma(z)$. \qed

\begin{rmk}
Since the Stokes matrices in $2\times 2$ case is given in a closed formula, we can see why Conjecture \ref{conj2} and Conjecture \ref{lastconj} hold, by studying the Stirling's formula of the gamma function as $h\rightarrow\infty$ and by studying the zeros of the gamma function respectively. 
One of the main results in \cite{Xu} is a construction of an $n\times n$ model of the Stokes matrices via the isomonodromy approach. The model makes it possible to study several basic analysis problems in the Stokes phenomenon. 
\end{rmk}

\subsection{A proof of Theorem \ref{introthm1re}} 
Following Theorem \ref{thmRLL}, for all $h\notin\mathbb{Q}$, the Stokes matrices $S_{h\pm}(u)$ satisfy the RLL relations. The Stokes matrices $S_{h\pm}^o(u)$, corresponding to a rational number $h_0\in\mathbb{Q}$, are the continuous extension of $S_{h\pm}(u)$ from $h\in\mathbb{C}\setminus\mathbb{Q}$ to $h_0$. Therefore, Theorem \ref{introthm1re} follows from the fact that the RLL relations are preserved under the continuous extension of $h$.

It can also be seen from the explicit expression in $n=2$ example, where clearly the case of $q$ is a root of unity also defines a representation.  

\section{Appendix: a proof of Theorem \ref{fact}}\label{lastpiece}
In this section, we only give a proof of $(b)$ in Theorem \ref{fact}. Parts $(a)$ and $(c)$ follow from a similar treatment. The proof is decomposed to five steps.

\subsection{Step 1: Existence of $W_i(z_1,z_2)$}\label{complete} 

\begin{lem}\label{esolW}
There exists a unique function $W_i(z_1,z_2)\in {\rm End}(L(\lambda))\otimes {\rm End}(\IC^n)\otimes {\rm End}(\IC^n)$, defined for $z_1\in\Sect_i$ (the sector $S(d_i-{\pi}/{2},d_{i+1}+{\pi}/{2})$ of the $z_1$ plane) and $|z_1|>|z_2|$, satisfying the equation \eqref{ratio1eq} and the prescribed asymptotics in \eqref{simW1}.
\end{lem}
\begin{proof} 
Let us follow the standard procedure: first to show that for any fixed $z_2\in \mathbb{C}\setminus\{0\}$, the equation \eqref{ratio1eq} has a unique formal fundamental solution taking the form \[
\widehat{W}(z_1;z_2)=\big(1+\sum_{m\ge 1}K_m(z_2)z_1^{-m})\cdot e^{ hz_1u^{(1)}}z_1^{h\delta T^{(1)}}(z_1-z_2)^{h\delta P},\] where each coefficient $K_m(z_2)\in {\rm End}(L(\lambda))\otimes {\rm End}(\IC^n)\otimes {\rm End}(\IC^n)$. 

Expanding $\frac{P}{z_1-z_2}=\frac{P}{z_1}(1+\sum_{k\ge 1}(\frac{z_2}{z_1})^k)$ and comparing the coefficients of $z_1^{-m-1}$ in the two sides of \eqref{ratio1eq}, we get the recursion relation
\begin{align}\label{formalK}\begin{split}
[K_{m+1}(z_2), u^{(1)}]=&\sum_{r=1}^m(P-\delta P) K_{m-r}(z_2)z_2^r-K_{m}(z_2)\cdot  \Big(\delta T^{(1)}+\delta P\Big)\\
&+\Big( m/h+T^{(1)}+P\Big)\cdot K_{m}(z_2).\end{split}\end{align}
Let us write $K_m=\sum_{i,j,k,l} K_{m, ijkl}\otimes E_{ij}\otimes E_{kl}$, where each $K_{m, ijkl}\in{\rm End}(L(\lambda))$. Then

$(1)$ For $i\ne j$ or $k\ne l$: 
\begin{align*}
(u_j-u_i) K_{m+1, ijkl}=&\sum_{j'} e_{ij'} K_{m, j'jkl}+\delta_{lj}K_{m, ijkl}-K_{m, kjil}-K_{m, ijkl}e_{jj}\\
&+\frac{m}{h}K_{m, ijkl}+\sum_{r=1}^m K_{m-r,kjil}z_2^r.
\end{align*}

$(2)$ For $i=j$ and $k=l$: 
\begin{align*}
0=&\sum_{i'\ne i} e_{ii'} K_{m, i'ikk}+\delta_{ik} K_{m, iiii}-K_{m, kiik}+[e_{ii}, K_{m, iikk}]\\&+\frac{m}{h}K_{m, iikk}+\sum_{r=1}^m (K_{m-r,kiik}-\delta_{ik} K_{m, iiii})z_2^r.
\end{align*}
Since $h\notin \mathbb{Q}$, a same argument as in the proof of Proposition \ref{uniformal}
shows the existence and uniqueness of $K_m$ for all $m\ge 1$. 

After knowing the formal solution, we note that the anti-Stokes rays of the equation \eqref{ratio1eq} are $-{\rm arg}(hu_i-hu_j)$ for all the possible $i\ne j$, which coincide with the rays of the equation \eqref{introqeq}. Again, we label the corresponding Stokes supersectors $\Sect_0,...,\Sect_{2l-1}$ in a same way.
Then, for any fixed $z_2$, the Borel resummation of $\widehat{K}=1+\sum_{m\ge 1} K_m z_1^{-m}$ defines holomorphic functions $K_i$ in each domain ${\rm Sect}_i\cap\{z_1\in \mathbb{C}~|~|z_1|>|z_2|\}$, which in turn determines the solutions of \eqref{ratio1eq} with the prescribed asymptotics \eqref{simW1}. \end{proof}

\begin{lem}\label{anothersol}
The product $W_i(z_1,z_2)F_i^{(2)}(z_2)$ satisfies the compatible system of equations
\begin{equation}\label{eqz1}
\frac{1}{h}\frac{\partial W_i(z_1,z_2)F_i^{(2)}(z_2)}{\partial z_1}= \Big(u^{(1)}+\frac{T^{(1)}}{z_1}+\frac{P}{z_1- z_2}\Big)\cdot W_i(z_1,z_2)F_i^{(2)}(z_2),
\end{equation}
\begin{equation}\label{eqz2}
 \frac{1}{h}   \frac{\partial W_i(z_1,z_2)F_i^{(2)}(z_2)}{\partial z_2}= \Big(u^{(2)}+\frac{T^{(2)}}{z_2}+\frac{P}{z_2- z_1}\Big)\cdot W_i(z_1,z_2)F_i^{(2)}(z_2)
\end{equation}
\end{lem}
\begin{proof} Since $F_i^{(2)}(z_2)$ is independent of $z_1$ The equation \eqref{eqz1} follows from the equation \eqref{ratio1eq} satisfied by $W_i$.

By the compatibility of the equations \eqref{eqz1} and \eqref{eqz2}, we have that the function
\begin{eqnarray*}
N_{d,d'}(z_1,z_2):= \frac{\partial (W_iF_i^{(2)})}{\partial z_2}-h\Big(u^{(2)}+\frac{T^{(2)}}{z_2}+\frac{P}{z_2-z_1}\Big)\cdot W_iF_i^{(2)}
\end{eqnarray*}
satisfies the equation \eqref{eqz1}. It implies that the ratio $C:=(W_iF_i^{(2)}(z_2))^{-1}\cdot N_{d,d'}$ is independent of $z_1$.
To show that $C$ is zero, let us rewrite 
\begin{eqnarray*}
W_i(z_1,z_2)\cdot F_i^{(2)}(z_2)&=&K_i(z_1,z_2)\cdot e^{ hz_1u^{(1)}}z_1^{h\delta T^{(1)}}(z_1-z_2)^{h\delta P} F_i^{(2)}(z_2)\\
&=&
K_i\cdot (1-z_2/z_1)^{h\delta P} F_i^{(2)}(z_2)e^{hu^{(1)}z_1} z_1^{h{\delta T^{(1)}}+h\delta P}.
\end{eqnarray*}
Here we use the fact that $e^{hu^{(1)}z_1}z_1^{h\delta T^{(1)}+h\delta P}$ commutes with the coefficient matrix of the equation \eqref{app2} and the initial value as $z_2\rightarrow\infty$, therefore commutes with the solution $F_i^{(2)}(z_2)$. 
Set $L_i=K_i\cdot (1-z_2/z_1)^{h\delta P}$. Differentiating the above identity with respect to $z_2$ gives
\begin{eqnarray*}
\frac{\partial W_iF_i^{(2)}}{\partial z_2}\cdot (W_iF_i^{(2)})^{-1}= \frac{\partial L_i}{\partial z_2} L_i ^{-1}+L_i\cdot\Big(u^{(2)}+\frac{T^{(2)}}{z_2}\Big)\cdot L_i^{-1}.
\end{eqnarray*}
It leads to 
\begin{align}\label{eqtW}\begin{split}
&\frac{1}{h}W_iF_i^{(2)}\cdot C\cdot (W_iF_i^{(2)})^{-1}\\
=& \frac{\partial L_i}{\partial z_2} L_i ^{-1}+L_i\cdot\Big(u^{(2)}+\frac{T^{(2)}}{z_2}\Big)\cdot L_i^{-1}- \left(u^{(2)}+\frac{T^{(2)}}{z_2}+\frac{P}{z_2-z_1}\right).
\end{split}
\end{align}
Recall that $K_i(z_1,z_2)\sim 1$ as $z_1\rightarrow\infty$ within $\Sect_i$. Then the right hand side of \eqref{eqtW} is $O(z_1^{-1})$ as $z_1\rightarrow\infty$ within $\Sect_i$.
Similar to the argument in Theorem \ref{solution}, one shows that $C=0$. This concludes the proof. \end{proof}

\subsection{Step 2: the $t\rightarrow 0$ asymptotic property of the function $Q_{d,d'}$ given in Theorem \ref{BorelY}}

Following Section \ref{secsolution}, let us take the functions $Q_{d,d'}(z,t)$ defined on the regions
$\Sect_i^\theta \times (-\delta,0)$ and $\Sect_i^\theta \times (0,\delta)$. Let us study the $t\rightarrow 0$ asymptotics of $Q_{d,d'}(z,t)$. First by the expression \eqref{offdia}, in the case of $i=j$ and $k\ne l$, the coefficient before $Q_{m+1,iikl}$ becomes $t(u_l-u_k)$, which will blow up as $t\rightarrow 0$. Thus, for fixed $z$ the function $Q_{d,d'}(z,t)$, obtained by the Borel-Laplace transform of $\hat{Q}$, may not have limit as $t\rightarrow 0$. However, in the following, let us show that $Q_{d,d'}(z,t)$ has limit as $t\rightarrow 0$ along the path $zt=z_2$ a constant. 

\begin{pro}\label{keypara}
$(a).$ Let $z_2\in \Sect_i^\theta $ be any none zero constant. Then the function $Q_{d,d'}$ satisfies
\begin{equation}
Q_{d,d'}(z,t)\sim f_{i+}(z_2), \ \text{as} \ t\rightarrow 0+ \text{ along the path $zt=z_2$ in $\Sect_i^\theta \times (0,\delta)$}
\end{equation}
where $f_{i+}(z_2)$ is a continuous function, and is asymptotic to $1$ as $z_2\rightarrow \infty$.

$(b).$ Let $z_2\in \overline{\Sect}_{i+l}$ be any none zero constant, where $\overline{\Sect}_{i+l}$ is the opposite sector of $\overline{\Sect}_{i}$. Then the function $Q_{d,d'}$ satisfies
\begin{equation}
Q_{d,d'}(z,t)\sim f_{i-}(z_2) \ \text{as} \ t\rightarrow 0- \text{ along the path $zt=z_2$ in $\Sect_i^\theta \times (-\delta,0)$}
\end{equation}
where $f_{i-}(z_2)$ is a continuous function, and is asymptotic to $1$ as $z_2\rightarrow \infty$.
\end{pro}
\begin{proof} The proof has three parts including the estimate of formal solutions, the Borel transform and the Laplace transform. It ends in the end of this subsection. Since the proofs of $(a)$ and $(b)$ are the same, here we only prove $(a)$. First note that the path $zt=z_2$ for $t\in (0,\delta)$ is in the domain $\Sect_i^\theta \times (0,\delta)$.

Let us introduce $Q'_m:=t^mQ_m$. Then multiplying by $t^{m}$ on the two sides of the identity \eqref{offdia} leads to 

\begin{align}\label{recurisont}\begin{split}
&\frac{u_j-u_i+tu_l-tu_k}{t} Q'_{m+1, ijkl}\\  
=&\sum_{j'} e_{ij'} Q'_{m, j'jkl}+\sum_{l'} e_{kl'} Q'_{m, ijl'l}+Q'_{m, kjil}\\
&-Q'_{m, ijkl}(\delta_{lj}+e_{jj})- Q'_{m, ijkl}e_{ll}+\frac{m}{h}Q'_{m, ijkl}.
\end{split}
\end{align}
Note that we can take a real number $r>0$ such that \begin{eqnarray}\label{estimatet}
r\le |(u_j-u_i+tu_l-tu_k)/t|, \ \ \text{for all } i,j,k,l \ \text{ and } \ 0<t< \delta.\end{eqnarray}
Then multiplying by $t^{m}$ on the two sides of the identity \eqref{offdia}, and using the above inequality, leads to
\begin{eqnarray*}
|Q'_{m+1, ijkl}|\le &\frac{1}{r}\Big|\sum_{j'} e_{ij'} Q'_{m, j'jkl}+\sum_{l'} e_{kl'} Q'_{m, ijl'l}+Q'_{m, kjil}
\\&-Q'_{m, ijkl}(\delta_{lj}+e_{jj})- Q'_{m, ijkl}e_{ll}+\frac{m}{h}Q'_{m, ijkl}\Big|.
\end{eqnarray*}
Again, from the inequality it is a standard fact that there exists constants $c'>0$ and $K>0$ independent of $t$ such that
\begin{eqnarray}\label{estimate}
    |Q'_m(t)|\le c' {K}^m m!, \ \ \text{ for all } t\in (0,\delta).
\end{eqnarray}
Let us think of $\hat{Q}=1+\sum_{m\ge 1}Q'_m(t)\cdot (zt)^{-m}$ as a power series in the new variable $z_2=zt$. In the following, let us first recall the Borel resummation $Q_{d,d'}(z_2;t)$ of $\hat{Q}$ as a holomorphic function of $z_2=zt$ on each sector $\overline{\Sect_i}$, and then prove that $Q_{d,d'}(z_2;t)\sim 1,$ as  $z_2\rightarrow \infty$ within $\Sect_i^\theta $ uniformly with respect to $t\in (0,\delta)$. 

\vspace{2mm}
{\bf Borel transform.}
Denote by $\mathcal{B}_{z_2}(\hat{Q})$ the formal Borel transform of the power series $\hat{Q}-1$ in the variable ${z_2}=zt$ (in the Borel plane with complex variable $\xi$), i.e., 
\[\mathcal{B}_{z_2}(\hat{Q})(\xi;t):=\sum_{m\ge 1}\frac{Q'_m(t)}{\Gamma(m)}\xi^{m-1}.\]

\begin{rmk}
Note that $\hat{Q}$ is a formal power series in the variable $z_2^{-1}$ (negative power), while its Borel transform is a power series in $\xi$ (positive power). Here, in order to be consistent with the convention in \cite{LR}, we implicitly make a coordinate change switching the zero and infinite point. Accordingly, the change of the argument should be accounted for. For example, it leads to the minus sign of the arguments $d_i$ and $d_{i+1}$ in terms of the variable $\xi$ in Lemma \ref{exgrow}.
\end{rmk}
Thus by \eqref{estimate} the Borel transform is convergent in a small neighborhood of $\xi=0$. Recall that $\theta>0$ is a chosen real number such that the sector $\Sect_i^\theta=S(d_{i}+
\theta-\pi/2,d_{i+1}-\theta+\pi/2)$ has opening angle bigger than $\pi$, and is contained in $\Sect_{d,d'}(t)$ for all $t\in (0,\delta)$. Then
\begin{lem}\label{exgrow}
For any $t\in (0,\delta)$, the Borel transform $\mathcal{B}_{z_2}(\hat{Q})(\xi;t)$ can be analytically continued to the sector ${S(-d_{i+1}+
\theta,-d_{i}-\theta)}$ of the $\xi$ plane, and there exist constants $\alpha, \beta>0$ (independent of $t$) such that
\begin{equation}\label{expontt}
    |\mathcal{B}_{z_2}(\hat{Q})(\xi;t)|\le \alpha e^{\beta|\xi|}, \ \text{for all} \ \xi\in   {S(-d_{i+1}+
\theta,-d_{i}-\theta)} \ \ \text{and} \ t\in (0,\delta).
\end{equation}
Furthermore, the limit $\mathcal{B}_{z_2}(\hat{Q})(\xi;0)$ of $\mathcal{B}_{z_2}(\hat{Q})(\xi;t)$ as $t\rightarrow 0+$ exists and is a holomorphic function on $ {S(-d_{i+1}+
\theta,-d_{i}-\theta)}$ satisfying the same inequality \eqref{expontt}.
\end{lem}
\begin{proof} Set $L=hT^{(1)}+hT^{(2)}+hP$, then the identity \eqref{Hhat} is formally equivalent to the integral equation
\begin{align}\begin{split}
&\mathcal{B}_{z_2}(\hat{Q})(\xi;t)\cdot (hu^{(1)}/t+hu^{(2)})-(hu^{(1)}/t+hu^{(2)}+\xi)\cdot \mathcal{B}_{z_2}(\hat{Q})(\xi;t)\\ \label{inteq} =&L-\delta L+\int_{x=0}^\xi\Big(L\cdot \mathcal{B}_{z_2}(\hat{Q})(x;t)-\mathcal{B}_{z_2}(\hat{Q})(x;t)\cdot \delta L\Big)dx.
\end{split}
\end{align}

To study the integral equation, we employ an iteration, by beginning with $\mathcal{B}_{z_2}^{(0)}(\hat{Q})(\xi;t)\equiv 0$, and plugging $\mathcal{B}_{z_2}^{(m)}(\hat{Q})(\xi;t)$ into the right hand side of \eqref{inteq} and determining $\mathcal{B}_{z_2}^{(m+1)}(\hat{Q})(\xi;t)$ from the left hand side. The sequence so obtained is holomorphic in $G$, the largest star-shaped region that does not contain any point $\xi$ for which $hu^{(1)}/t+hu^{(2)}$ and $hu^{(1)}/t+hu^{(2)}+\xi$ have an eigenvalue in common. In particular, the set $G$ contains the sector $S(-d'(t),-d(t))$ in the $\xi$ plane, and therefore contains the sector $ {S(-d_{i+1}+
\theta,-d_{i}-\theta)}$ for the sufficiently small $\theta>0$.

Now let $\rho$ be any ray in the sector ${S(-d_{i+1}+
\theta,-d_{i}-\theta)}$ of $\xi$ plane, and let $\varepsilon>0$ be small enough such that $S(\rho-\varepsilon,\rho+\varepsilon)$ is inside $ {S(-d_{i+1}+
\theta,-d_{i}-\theta)}$. Let us show that each $\mathcal{B}_{z_2}^{(m)}(\hat{Q})(\xi;t)$ has exponential growth
of order $1$ at $\infty$. For this purpose, let us set $a=|A|$. By comparing the coefficient of $\xi^{k-1}$ of the iterated formula of $\mathcal{B}_{z_2}^{(m)}(\hat{Q})(\xi;t)$ for all $m$, inductively we get the estimates of the form
\begin{eqnarray}\label{estBorel}
    |\mathcal{B}_{z_2}^{(m)}(\hat{Q})(\xi;t)|\le \sum_{k\ge 1}{b^{(m)}_k}|\xi|^{k-1}/\Gamma(k),
\end{eqnarray}
where $b^{(m)}_k$ are determined by the recursive relation
\begin{eqnarray}\label{trecu}
b^{(m+1)}_k=c\cdot (a+2ab^{(m)}_{k-1}),
\end{eqnarray}
where $c$ is a constant which arises when solving the left hand side of \eqref{inteq}. The constant $c$ can be chosen independent of $t\in (0,\delta)$ and the direction $\rho$, because for all $t$ the numbers $hu_i/t-hu_j/t+hu_k-hu_l$ keep a fixed positive distance from the sector $S(\rho-\varepsilon,\rho+\varepsilon)\subset  {S(-d_{i+1}+
\theta,-d_{i}-\theta)}$.

Now for every $k$, the
numbers $b^{(m+1)}_k$ monotonically increase with respect to $m$ and become constant when $m\ge k$ (due to the fact that the initials $b^{(0)}_k=0$ for all $k$). The limiting values $b_k$ satisfy the same recursion equation as \eqref{trecu}, which implies that $b_k$ cannot grow faster than $\beta^k$ for some constant $\beta>0$. Therefore, each $\mathcal{B}_{z_2}^{(m)}(\hat{Q})(\xi;t)$ is estimated by $\alpha e^{\beta|\xi|}$ with suitably $\alpha,\beta$ independent of $t\in (0,\delta)$. 

Next let us show the convergence of $\mathcal{B}_{z_2}^{(m)}(\hat{Q})(\xi;t)$ as $m\rightarrow\infty$. It can be done by deriving estimates for the differences
$\mathcal{B}_{z_2}^{(m)}(\hat{Q})(\xi;t)- \mathcal{B}_{z_2}^{(m-1)}(\hat{Q})(\xi;t)$ and turning the sequence into a telescoping sum: similar to the estimate \eqref{estBorel}, inductively we get, from the iterated formula \eqref{inteq} for $\mathcal{B}_{z_2}^{(m)}(\hat{Q})(\xi;t)- \mathcal{B}_{z_2}^{(m-1)}(\hat{Q})(\xi;t)$ and $\mathcal{B}_{z_2}^{(m-1)}(\hat{Q})(\xi;t)- \mathcal{B}_{z_2}^{(m-2)}(\hat{Q})(\xi;t)$, that 
\begin{eqnarray}\label{estBd}
    |\mathcal{B}_{z_2}^{(m)}(\hat{Q})(\xi;t)- \mathcal{B}_{z_2}^{(m-1)}(\hat{Q})(\xi;t)|\le \sum_{k\ge m}w^{k-1}|\xi|^{k-1}/\Gamma(k),
\end{eqnarray}
where $w>0$ is a sufficient large constant that can be chosen independent of $t$. 

The estimate \eqref{estBd} shows that $\mathcal{B}_{z_2}^{(m)}(\hat{Q})(\xi;t)$ converges
uniformly on every compact subset of the region $S(\rho-\varepsilon,\rho+\varepsilon)\subset  {S(-d_{i+1}+
\theta,-d_{i}-\theta)}$, and uniformly with respect to $t\in (0,\delta).$ Since the direction $\rho$ is arbitrary, it proves the first part. 

For the second part, note that the limit $\mathcal{B}_{z_2}^{(m)}(\hat{Q})(\xi;0)$ of each iterated term $\mathcal{B}_{z_2}^{(m)}(\hat{Q})(\xi;t)$ exists as $t\rightarrow 0+$, and is a holomorphic function on the sector $ {S(-d_{i+1}+
\theta,-d_{i}-\theta)}$. Furthermore, $\mathcal{B}_{z_2}^{(m)}(\hat{Q})(\xi;t)$ converges as $t\rightarrow 0+$
uniformly on every compact subset of the region $S(\rho-\varepsilon,\rho+\varepsilon)$. In the end, since the constants $b^{(m)}_k$ in \eqref{estBorel} and the constant $w$ in \eqref{estBd} are independent of $t\in (0,\delta)$, we get the convergence of $\mathcal{B}_{z_2}^{(m)}(\hat{Q})(\xi;0)$ as $m\rightarrow\infty$, as well as the estimate of the resulting function $\mathcal{B}_{z_2}(\hat{Q})(\xi;0)$ as $m\rightarrow\infty$. \end{proof}

\vspace{2mm}
{\bf Laplace transform.}
The Laplace transform of the function $\mathcal{B}_{z_2}(\hat{Q})(\xi;t)$ in the direction $\rho$ is a function (in the Laplace plane of the initial
variable ${z_2}$) defined by
\[Q_\rho({z_2};t)=1+\int_{\xi=0}^{+\infty\cdot e^{\mathi\rho}}e^{-{z_2}{\xi}} \mathcal{B}_{z_2}(\hat{Q})(\xi;t) d\xi.\]
Here the line of integration is the line with argument $\rho$.
The inequality \eqref{expontt} ensures that for any $-d_{i+1}+
\theta\le \rho\le -d_{i}-\theta$, the integrand is indeed defined on the integral path, and that for any fixed $\beta'>\beta$ the integral exists for all $z_2\in R(\rho,\beta')$ and $t\in (0,\delta)$, where the domain
\[R(\rho,\beta'):=\{{z_2}\in\mathbb{C}~|~{\rm Re}(z_2 e^{\mathi \rho})>\beta', |\rho+{\rm arg}({z_2})|<\pi/2\}.\]
By Lemma \ref{exgrow}, there exists a continuous function $f_{\rho+}(z_2)$ such that $Q_{\rho}(z,t)\sim f_{\rho+}(z_2)$ as $t\rightarrow 0+$. Next, let us show that $f_{\rho+}(z_2)$ approaches to $1$ as $z_2\rightarrow \infty$. To see this, we only need to show that 
\begin{lem}
There exist constants $C,D>0$ (independent of $t$) such that 
\begin{equation}\label{BLestineq}
|Q_\rho({z_2};t)-\sum_{m=0}^{N-1}Q'_m(t){z_2}^{-m}|\le CN^Ne^{-N}\frac{|{z_2}|^{-N}}{D^N}, 
\end{equation}
for all ${z_2}\in R(\rho,\beta'), \ t\in(0,\delta), \ \ N\in\mathbb{N}_+.$
\end{lem}
\begin{proof} For a fixed $t$, the proof of the inequality \eqref{BLestineq} is standard, see e.g., \cite[Theorem 5.3.9]{LR}. In the following, we will go through the proof given in \cite{LR}, and show that the involved constants $D$ and $L$ can be chosen independent of $t\in(0,\delta)$.

Now without lose of generality, let us assume that $\rho=0$. By the inequality \eqref{expontt}, there exist a constant $\varepsilon$ such that $\mathcal{B}_{z_2}(\hat{Q})(\xi;t)$ is holomorphic in the union of $\{\xi:|\xi|\le 1/K\}$ and the sector $\{-\varepsilon<{\rm arg}(\xi)<\varepsilon\}$.
Just as in \cite[Theorem 5.3.9]{LR}, let us take a point $b$ with argument $\pi/4$ and small enough norm $|b|<1/K$ such that the path, following a straight line from $0$ to $b$ and
continues along a horizontal line from $b$ to $+\infty$, lies in the domain $\{\xi:|\xi|\le 1/K\}\cup \{-\varepsilon<{\rm arg}(\xi)<\varepsilon\}$ of the $\xi$ plane. Since the path is homotopy to $[0,+\infty)$, 
by the Cauchy’s theorem, the Laplace integral $Q_{\rho=0}({z_2};t)$ decomposes to
$Q_{\rho=0}(z_2;t)=1+Q^b({z_2};t)+Y^b(z_2;u)$ along the path, where
\[Q^b(z_2;t)=\int_{\xi=0}^{b}e^{-z_2{\xi}}\mathcal{B}_{z_2}(\hat{Q})(\xi;t)d\xi, \hspace{5mm} Y^b({z_2};t)=\int_{\xi=b}^{+\infty}e^{-z_2{\xi}}\mathcal{B}_{z_2}(\hat{Q})(\xi;t)d\xi. \]
On the one hand, given $0<\gamma<\pi/2$, following \cite[Lemma 1.3.2]{LR}, we have for $-3\pi/4+\gamma<{\rm arg}(z_2) <\pi/4-\gamma$, $N\in\mathbb{N}_+$,
\begin{eqnarray*}
|Q^b(z_2;t)-\sum_{m=0}^{N-1}Q'_m(t)z_2^{-m}|\le C'N^Ne^{-N}\frac{|z_2|^{-N}}{{D'}^N}, 
\end{eqnarray*}
where for any given $t\in (0,\delta)$ the constants $C', D'$ are
\begin{eqnarray}
C':=\sum_{m\ge 1}\frac{|Q'_m(t)|}{\Gamma(m)}b^{m}, \hspace{5mm} D':=b\cdot {\rm sin}(\gamma).
\end{eqnarray}
On the other hand, let us take the constant $\beta'>\beta$, then following the proof of \cite[Theorem 5.3.9]{LR}, we have for $z\in \{-3\pi/4+\gamma<{\rm arg}(z_2) <\pi/4-\gamma\}\cap R(\rho=0,\beta')$, 
\begin{eqnarray}\label{exgrowest}
|Y^b(z_2;t)|\le pe^{-{c}{z_2}}, 
\end{eqnarray} with the constants given by
\begin{eqnarray}
p=\frac{\alpha e^{\beta|b|}}{\beta'-\beta}, \hspace{5mm} c=|b|{\rm cos}(\pi/2-\gamma).
\end{eqnarray}
The estimation \eqref{exgrowest} further implies, see e.g,. \cite[Proposition 1.2.17]{LR}, for $z_2\in \{-3\pi/4+\gamma<{\rm arg}(z_2) <\pi/4-\gamma\}\cap R(\rho=0,\beta')$
\begin{eqnarray*}
|Y^b(z_2;t)|\le C''N^Ne^{-N}\frac{|z_2|^{-N}}{{D''}^N}, \ \ N\in\mathbb{N}_+,
\end{eqnarray*}
with the constants $C''$ and $D''$ determined by $p$ and $c$. In conclusion, if we take $C = {\rm max}(C', C'')$ and $D = {\rm max}(D', D'')$, then $Q_{\rho=0}(z_2;t)$ satisfies the inequality \eqref{BLestineq} on the domain
\[\{-3\pi/4+\gamma<{\rm arg}(z_2) <\pi/4-\gamma\}\cap R(\theta=0,\beta').\]
According to \cite[Theorem 5.3.9]{LR}, an argument using the symmetry with respect to the real axis, i.e., by choosing $\bar{b}$ instead of $b$ and the corresponding path, shows that $Q_{\rho=0}(z_2;t)$ satisfies the inequality \eqref{BLestineq} on the symmetric domain 
\[\{-\pi/4+\gamma<{\rm arg}(z_2) <3\pi/4-\gamma\}\cap R(\theta=0,\beta')\]
with respect to the real axis. Since the union of the above two domains cover $R(\rho=0,\beta')$, we get $Q_{\rho=0}(z_2;t)$ satisfies the inequality \eqref{BLestineq} on $R(\rho=0,\beta')$. 

In the end, let us check the independence of constants $C = {\rm max}(C', C'')$ and $D = {\rm max}(D', D'')$ on $t\in (0,\delta)$. First, following the inequality \eqref{estimate}, we have $C'\le \frac{Kb}{1-Kb}$. Thus we can set $C'=\frac{Kb}{1-Kb}$, and then the constants $C,D$ are determined by $\varepsilon, K,\alpha,\beta,b$ and $\delta$. By \eqref{estimate} and \eqref{expontt}, as a sufficiently small $\varepsilon$ fixed, those constants can be chosen independent of $t\in (0,\delta)$. 
It therefore verifies \eqref{BLestineq}. \end{proof} 

\vspace{2mm}

Now the Borel-Laplace transform $Q_\rho$ and $Q_{\rho'}$ of $\hat{Q}(z_2;t)$, with respect to two directions satisfying $-d_{i+1}+
\theta\le \rho, \rho' \le -d_{i}-\theta$, coincide in the overlapping of their defining domains. See e.g., \cite[Proposition 5.3.7]{LR} or \cite[Section 6.2]{Balser}. (While if $\rho$ and $\rho'$ are not in a same sector bounded by two adjacent Stokes rays,  $Q_\rho$ is in general not equal to $Q_{\rho'}$ at the points where both functions are defined.)  Thus the functions $Q_\rho$ for all $-d_{i+1}+
\theta\le \rho \le -d_{i}-\theta$ in the $\xi$ plane glue together into a holomorphic function $Q_{d,d'}(z_2;t)$ defined on the domain $\Sect_i^\theta \times (0,\delta)$. By \eqref{estimate} and \eqref{expontt}, there exists a continuous function $f_{i+}(z_2)$ (the gluing of $f_{\rho+}(z_2)$ for all $\rho$ with $-d_{i+1}+
\theta\le \rho \le -d_{i}-\theta$) such that $Q_{d,d'}(z_2;t)\sim f_{i+}(z_2)$ as $t\rightarrow 0+$. Furthermore, by \eqref{BLestineq}, we have
\[Q_{d,d'}(z_2;t)\sim 1, \hspace{3mm} \text{as} \ z_2\rightarrow \infty \ \ \text{within} \ \ \Sect_i^\theta \]
uniformly for $t\in (0,\delta)$. It therefore proves the Proposition \ref{keypara}. 
\end{proof} 

\vspace{3mm}

\subsection{Step 3: Proof of the first identity in part (b) of Theorem \ref{fact}}

Let $Q_{d,d'}$ be the function as in Section \ref{secsolution}, and $Y_{d,d'}$ the corresponding solution given by \eqref{QtoY}. 

\begin{pro}\label{prob}
The solution satisfy
\begin{equation}\label{0del}
Y_{d,d'}(z,t)=W_i(z,zt)F_i^{(2)}(zt)e^{\pi\mathi  h\delta P}, \text{ for }  z\in \Sect_{i}^\theta  \text{ and } t\in (0, \delta),
\end{equation}
\begin{equation}
Y_{d,d'}(z,t)=W_i(z,zt)F_{l+i}^{(2)}(zt)e^{\pi\mathi  h\delta P}, \text{ for }  z\in\Sect_{i}^\theta  \text{ and } t\in(-\delta,0).
\end{equation}
\end{pro}
\begin{proof} Let us only prove \eqref{0del}. Note that the system of equations \eqref{gcKZ1} and \eqref{gcKZ2} becomes the system of \eqref{eqz1} and \eqref{eqz2} in terms of the coordinates transform $z_1=z$ and $z_2=zt$. Following Lemma \ref{anothersol}, $Y_{d,d'}(z,t)$ and $W_i(z,zt)F_i^{(2)}(zt)e^{\pi\mathi  h\delta P}$ satisfy the same system. Thus, we have $Y_{d,d'}=W_iF_i^{(2)}e^{\pi\mathi  h\delta P}C$ for a constant $C$ (independent of $z_1$ and $z_2$).

Recall $F^{(2)}_i(z_2)=H^{(2)}_i(z_2) \cdot e^{hu^{(2)}z_2}z_2^{h\delta T^{(2)}}$ and $W_i=K_i\cdot e^{ hz_1u^{(1)}}z_1^{h\delta T^{(1)}}(z_1-z_2)^{h\delta P}$. To proceed, let us write
\begin{eqnarray*}
&&Y_{d,d'}(z_1,z_2)\cdot e^{-\pi\mathi  h\delta P}F_{i}^{(2)}(z_2)^{-1}W_i(z_1,z_2)^{-1}\\
%&=&Q_{d,d'}(z_1,z_2)\cdot (z_2/z_1-{1})^{h\delta P} e^{-\pi\mathi  h\delta P}e^{hu^{(1)}z_1}z_1^{h\delta T^{(1)}+h\delta P} H_i^{(2)}(z_2)^{-1}e^{ -hz_1u^{(1)}}z_1^{-h\delta T^{(1)}}(z_1-z_2)^{-h\delta P} K_i^{-1}\\
&=&Q_{d,d'}(z_1,z_2)(1-z_2/z_1)^{h\delta P} H_i^{(2)}(z_2)^{-1}\cdot (1-z_2/z_1)^{-h\delta P} K_i^{-1}.
\end{eqnarray*}
Here in the identity we use the fact that $e^{hu^{(1)}z_1}z_1^{h\delta T^{(1)}+h\delta P}$ commutes with the coefficient matrix of the equation \eqref{app2} and the initial value as $z_2\rightarrow\infty$, therefore commutes with the function $H_i^{(2)}(z_2)$.

Now following the prescribed asymptotics in the definition of canonical solutions, for sufficiently large $z_2$ the function $H_i^{(2)}(z_2)$ is sufficiently close to $1$. Furthermore, for the sufficiently large but fixed $z_2$, both $Q_{d,d'}$ (following from Proposition \ref{keypara}), $K_i(z_1,z_2)$ (following from Lemma \ref{esolW}) and $(1-{z_2}/{z_1})^{h\delta P}$ are asymptotic to $1$ as $z_1\rightarrow \infty$ within $\Sect_{i}^\theta \subset \Sect_i$. It gives $W_iF_{i}^{(2)}e^{\pi\mathi  h\delta P}Ce^{-\pi\mathi  h\delta P}(W_iF_{i}^{(2)})^{-1}\rightarrow 1$ for $z_2\rightarrow\infty$ and $z_2/z_1=t\rightarrow 0$ in proper domains. Since the exponential terms dominate and the opening of $\Sect_{i}^\theta $ is bigger than $\pi$, the constant $C$ must be $1$. \end{proof}

\vspace{2mm}

Note that the solution $W_i(z,zt)F_{i}^{(2)}(zt)e^{\pi\mathi  h\delta P}$ is associated to the two adjacent rays $d_i$ and $d_{i+1}$ of \eqref{introqeq}, while the solution $Y_{d,d'}(z,t)$ is associated to the two rays $d(t)$ and $d'(t)$ of \eqref{gcKZ1}. However, there are many other adjacent pair of rays of \eqref{gcKZ1} that are bounded by $d_i$ and $d_{i+1}$, $d(t)$ and $d'(t)$ is only one of them. Actually, we can generalize Proposition \ref{prob} to
\begin{pro}\label{extend0} 
Assume that $l(t),l'(t)$ are any pair of adjacent rays of \eqref{gcKZ1} that are bounded by $d_i$ and $d_{i+1}$ for all $t\in (-\delta,0)\cup (0,\delta)$, and $Y_{l,l'}(z,t)$ the corresponding solutions on the domain $t\in (-\delta,0)$, $z\in \Sect_{l,l'}(t)$ and $t\in (0,\delta)$, $z\in \Sect_{l,l'}(t)$ as in Theorem \ref{solution}. Then there exists a sufficiently small positive real number $1/\beta'>0$ such that for $z\in \Sect_{l,l'}(t)$ and $t\in (0, \delta)\cap (0,1/\beta')$,
\begin{equation}\label{l'conn}
Y_{l,l'}(z,t)=W_i(z,zt)F_{i}^{(2)}(zt)\cdot e^{\pi\mathi  h\delta P},\end{equation} 
and for $z\in \Sect_{l,l'}(t)$ and $t\in (-\delta,0)\cap(-1/\beta',0)$,
\begin{equation}
    Y_{l,l'}(z,t)=W_i(z,zt)F_{l+i}^{(2)}(zt)\cdot e^{\pi\mathi  h\delta P}.
\end{equation}    
\end{pro}
\begin{proof} The proof of this proposition include two lemmas, and it ends in the end of this subsection.
As a consequence of Proposition \ref{prob} (mainly based on Proposition \ref{keypara}), the function $W_i(z_1,z_2)$ given in Lemma \ref{esolW} has the following asymptotics
\begin{eqnarray}\label{Wasym}
    W_i(z_1,z_2)\cdot e^{- hz_1u^{(1)}}z_1^{-h\delta T^{(1)}}(z_1-z_2)^{-h\delta P}\sim 1
\end{eqnarray}
as $z_1\rightarrow \infty$ within ${\Sect_i}^\theta$ along the path $z_2/z_1$ being a constant inside $(-\delta,0)\cup (0,\delta)$. Actually we can prove the asymptotics \eqref{Wasym} holds in a larger region. To see this, let us check the Borel-Laplace transform of $\hat{W}$ in details. Following Lemma \ref{esolW}, the unique formal fundamental solution of the equation \eqref{ratio1eq} takes the form \[
\widehat{W}(z_1;z_2)=\big(1+\sum_{m\ge 1}K_m(z_2)z_1^{-m})\cdot e^{ hz_1u^{(1)}}z_1^{h\delta T^{(1)}}(z_1-z_2)^{h\delta P},\] 
where $K_m$ is recursively determined by \eqref{formalK}. Let us introduce $K'_m:=z_2^{-m}K_m$. Then multiplying by $z_2^{-m}$ on the two sides of the identity \eqref{formalK} leads to 
\begin{align*}
z_2\cdot [K'_{m+1}(z_2), u^{(1)}]=&\Big( m/h+T^{(1)}+P\Big)\cdot K'_{m}(z_2)-K'_{m}(z_2)\cdot  \Big(\delta T^{(1)}+\delta P\Big)\\
&+\sum_{r=1}^m(P-\delta P) K'_{m-r}(z_2).\end{align*}
Let us think of $\hat{K}=1+\sum_{m\ge 1}K'_m(z_2)\cdot (z_2/z_1)^{m}$ as a power series in the new variable $t=z_2/z_1$. From the above recursive relation, we see that the norm of its coefficient $K'_{m+1}$ is proportional to $1/|z_2|$. 

Denote by $\mathcal{B}_{t}(\hat{K})$ the formal Borel transform of the power series $\hat{K}-1$ in the variable ${t}=z_2/z_1$ (in the Borel plane with complex variable $\xi$), i.e., 
\[\mathcal{B}_{t}(\hat{K})(\xi;z_2):=\sum_{m\ge 1}\frac{K'_m(z_2)}{\Gamma(m)}\xi^{m-1}.\]
Set $V=hT^{(1)}+hP$, then the identity \eqref{formalK} is formally equivalent to the integral equation
\begin{align*}
&\mathcal{B}_{t}(\hat{K})(\xi;z_2)\cdot hz_2u^{(1)}-(hz_2u^{(1)}+\xi)\cdot \mathcal{B}_{t}(\hat{K})(\xi;z_2)=V-\delta V\\&+\int_{x=0}^\xi\Big(V\cdot \mathcal{B}_{t}(\hat{K})(\xi;z_2)-\mathcal{B}_{t}(\hat{K})(\xi;z_2) \delta V+\mathcal{B}_{t}(\hat{K})(\xi;z_2)\cdot \sum_{k\ge 1}\frac{(\xi-x)^k}{k!} \Big)dx.
\end{align*}

Here we remark that $\hat{K}$ is a formal power series in the variable $z_1^{-1}$ (negative power), while its Borel transform is a power series in $t$ (positive power). It explains the minus sign of the arguments $d_i$ and $d_{i+1}$ in the following lemma.
\begin{lem}
For any $z_2=|z_2|e^{\mathi \delta}$ and any small positive real number $\varepsilon$, the Borel transform $\mathcal{B}_{t}(\hat{K})(\xi;z_2)$ can be analytically continued to the sector \[{S(-d_{i+1}+\delta+\varepsilon,-d_{i}+\delta-\varepsilon)}\] of the $\xi$ plane, and if $z_2$ keeps a fixed positive distance from $0$ then there exist constants $\alpha, \beta>0$ (independent of $|z_2|$) such that
\begin{equation}\label{estWexp}
    |\mathcal{B}_{t}(\hat{K})(\xi;z_2)|\le \frac{1}{|z_2|}\alpha e^{\beta|\xi|}, \ \text{ for all } \ \xi\in   S(-d_{i+1}+\delta+\varepsilon,-d_{i}+\delta-\varepsilon).
\end{equation}
\end{lem}
\begin{proof}
It follows from the similar argument as in Lemma \ref{exgrow}: one applies the same estimate of the integral equation for the product function $z_2\cdot \mathcal{B}_{t}(\hat{K})(\xi;z_2)$ to find the constants $\alpha$ and $\beta$. Here we just stress that the constant $\alpha$ in \eqref{estWexp} relies on the distance from the numbers $hu_iz_2-hu_jz_2$ and the sector $S(-d_{i+1}+\delta+\varepsilon,-d_{i}+\delta-\varepsilon)$, which is independent of the norm $|z_2|$ (as long as $|z_2|$ is bigger than a fixed positive real number as in our assumption) but depends on the number $\varepsilon$ and the argument $\delta$ of $z_2$. \end{proof}

Now the Laplace transform of the function $\mathcal{B}_{t}(\hat{K})(\xi;z_2)$ along the direction $\rho$, with $-d_{i+1}+\delta+\varepsilon<\rho<-d_{i}+\delta-\varepsilon$, is a function (in the Laplace plane of the initial
variable $t=z_2/z_1$) defined by
\[K_\rho(t;z_2)=1+\int_{\xi=0}^{+\infty\cdot e^{\mathi\rho}}e^{-{\xi}{t}} \mathcal{B}_{t}(\hat{K})(\xi;z_2) d\xi.\]
The inequality \eqref{estWexp} ensures that for any fixed real number $\beta'>\beta>0$, the integral exists for all $t\in R(\rho,\beta')$, where the domain
\[R(\rho,\beta'):=\{t\in\mathbb{C}~|~{\rm Re}(t e^{\mathi \rho})>\beta', |\rho+{\rm arg}(t)|<\pi/2\}.\]

Similar to the standard estimate \eqref{BLestineq}, we can use the estimate \eqref{estWexp} to determine the asymptotics of $K(t;z_2)$ as $t\in R(\rho,\beta')$ and the constants involved. In particular, we get (the proof of the following lemma is skipped)
\begin{lem}
There exist constants $C,D>0$ (independent of $|z_2|$) such that for all $t\in R(\rho,\beta')$, $N\in\mathbb{N}_+$
\begin{eqnarray}\label{extK}
|K_\rho(t;z_2)-\sum_{m=0}^{N-1}K'_m(z_2){t}^{-m}|\le \frac{1}{|z_2|}CN^Ne^{-N}\frac{|t|^{-N}}{D^N}.
\end{eqnarray}
\end{lem}
All the functions $K_\rho(t;z_2)$ coincide with each other in their common domain as $-d_{i+1}+\delta+\varepsilon<\rho<-d_{i}+\delta-\varepsilon$, and thus glue to a function $K_i(t;z_2)$ which satisfies the inequality \eqref{extK} as \[t\in S(-d_{i+1}+\delta+\varepsilon-\pi/2,-d_{i}+\delta-\varepsilon+\pi/2) \ \text{ and } \ |t|<\varepsilon_1,\]
with $\varepsilon_1$ a fixed sufficiently small number. Therefore, when the argument $\delta$ of $z_2$ satisfies 
\begin{eqnarray}\label{deltavar}
    d_i+\varepsilon-\pi/2<\delta<d_{i+1}-\varepsilon+\pi/2,
\end{eqnarray} 
the above sector in the $t$ plane includes the positive real axis.

Since the function $W_i(z_1,z_2)$ given in Lemma \ref{esolW} is obtained by the Borel-Laplace transform, thus by definition 
\[W_i(z_1,z_2)=K_i({z_2}/{z_1};z_2)e^{ hz_1u^{(1)}}z_1^{h\delta T^{(1)}}(z_1-z_2)^{h\delta P}.\]
By \eqref{extK} and \eqref{deltavar}, for any fixed $t=z_2/z_1\in (0,1/\beta')$ the function \[K({z_2}/{z_1};z_2)\rightarrow 1, \ \text{ as } \ z_2\rightarrow \infty \ \text{  within } S(d_i+\varepsilon-\pi/2,d_{i+1}-\varepsilon+\pi/2).\] 
Note that $\varepsilon$ can be any small positive real number. In particular, for any fixed $t\in(0,\delta)\cap (0,1/\beta')$ the number $\varepsilon$ can be small enough such that $S(d_i+\varepsilon-\pi/2,d_{i+1}-\varepsilon+\pi/2)$ contains the sector $\Sect_{l,l'}(t)$. Using the same argument as in Proposition \ref{prob}, it implies that for any fixed $t\in(0,\delta)\cap (0,1/\beta')$, the two solutions $Y_{l,l'}(z,t)$ and $W_i(z,zt)F_{i}^{(2)}(zt)\cdot e^{\pi\mathi  h\delta P}$ have the same prescribed asymptotics as $z\rightarrow\infty$ within $\Sect_{l,l'}(t)\subset \Sect_i$. Then by the uniqueness statement of the solution in sectorial region $\Sect_{l,l'}(t)$ with opening bigger than $\pi$, the identity \eqref{l'conn} holds. \end{proof}

\subsection{Step 4: Triviality of the Stokes factors and extension with respective to the $t$ variable}\label{trivialS}
Let us denote by $d(t), d'(t), d''(t)$ three separate adjacent anti-Stokes rays of the equation \eqref{gcKZ1} (in anticlockwise order) for $t\in (0,\delta)$ that are bounded between $d_i$ and $d_{i+1}$. Then Proposition \ref{extend0} implies that for any $t\in(0,\delta)$
\begin{eqnarray}\label{extension}
Y_{d,d'}=Y_{d',d''} \hspace{3mm} \text{on} \ \ \Sect_{d,d'}\cap \Sect_{d',d''}.\end{eqnarray}
Recall that the Stokes factor $S_{d'}$ of \eqref{gcKZ1} is defined by $Y_{d',d''}=Y_{d,d'}\cdot S_{d'}$. As a consequence, we get
\begin{cor}\label{pStrivial}
If $t\in (0,\delta)$, then $S_{d'}$ must be ${\rm Id}$ unless the anti-Stokes ray $d'(t)$ has the argument $-{\rm arg}\big(h(u_i-u_j)\big)$ for some $i\ne j$.
\end{cor}

The proof of Proposition \ref{prob} implies the triviality of some Stokes factors of the equation \eqref{gcKZ1}. That is
\begin{pro}\label{Strivial}
If $t\in (0,1)$, then $S_{d'}$ must be ${\rm Id}$ unless the anti-Stokes ray $d'(t)$ equals to $-{\rm arg}\big(h(u_i-u_j)\big)$ for some $i\ne j$.
\end{pro}

\begin{proof} 
For any fixed $t$, we call an anti-Stokes ray $d(t)$ effective if its corresponding Stokes matrix $S_{d}$ is not an identity matrix. We call an anti-Stokes ray $d(t)$ ineffective if its Stokes matrix $S_{d'}$ is an identity matrix. In the following, let us prove that effective and ineffective rays will still be effective and ineffective respectively, as $t$ moves from a point near $0$ to $1$.

Let us move $t$ from the neighbourhood of $0$ to $1$. Let us assume that as $t$ crosses from one side of a point $t_0\in (0,1)$ to the other, only two rays $ d_-(t), d'_-(t)$ collide and two new rays $d_+(t), d'_+(t)$ emerge. 
That is if we denote by $d_l(t), d_-(t), d'_-(t),d_r(t)$ four separate adjacent anti-Stokes rays of \eqref{gcKZ1} in anticlock order for $t\in (t_-,t_0)$, and $d_l(t), d_+(t), d'_+(t),d_r(t)$ four separate adjacent rays for $t\in (t_0,t_+)$, then $d_l(t_0)\ne d_\pm(t_0)=d'_\pm(t_0)\ne d_r(t_0) $ at $t=t_0$. We can think of $d_-(t)$ becomes $d'_+(t)$ after $t$ passes through $t=t_0$. In the following, we show that if the Stokes factor $S_{d_-}$ associated to $d_-(t)$ for $t_-<t<t_0$ is identity, then so is the Stokes factor $S_{d'_+}$ associated to $d'_+(t)$ for $t_0<t<t_+$.

Let us recall that the Stokes factors $S_{d_\pm}$, $S_{d'_\pm}$, independent of $t$ and $z$, are determined by 
\begin{align*}
  Y_{d_-,d'_-}=&  Y_{d_l,d_-}\cdot S_{d_-} \text{ in } \Sect_{d_l,d_-}(t)\cap \Sect_{d_-,d'_-}(t), \\
  Y_{d'_-,d_r}=&  Y_{d_-,d'_-}\cdot S_{d'_-} \text{ in } \Sect_{d_-,d'_-}(t)\cap \Sect_{d',d_r}(t),
\end{align*}
and 
\begin{align*}
 Y_{d_+,d'_+}=&  Y_{d_l,d_+}\cdot S_{d_+} \text{ in } \Sect_{d_l,d_+}(t)\cap \Sect_{d_+,d'_+}(t), \\ Y_{d'_+,d_r}=&  Y_{d_+,d'_+}\cdot S_{d'_+} \text{ in } \Sect_{d_+,d'_+}(t)\cap \Sect_{d',d_r}(t).
\end{align*}
Thus we have $Y_{d'_-,d_r}=Y_{d_l,d_-}\cdot S_{d_-}S_{d'_-}$ and $Y_{d'_+,d_r}=Y_{d_l,d_+}\cdot S_{d_+}S_{d'_+}$. Since $d_l(t),d_-(t)$ (resp. $d_l(t),d_+(t)$) are adjacent and do not collide for $t\in (t_-,t_0]$ (resp. $t\in [t_0,t_+)$), by Theorem \ref{solution} we must have $Y_{d_l,d_-}=Y_{d_l,d_+}$ after the continuation of $t$ from $(t_-,t_0)$ to $(t_0,t_+)$. Similarly, we have $Y_{d'_-,d_r}=Y_{d'_+,d_r}$ after the continuation of $t$. In the end, we get \[S_{d_-}S_{d'_-}=S_{d_+}S_{d'_+}.\]

Using the uniqueness of the factorization of a matrix into an ordered product of two matrices satisfying the triangular property (see Lemma \ref{Sfactor}), we conclude that if $S_{d_-}$ is identity matrix, so is $S_{d'_+}$. We remark that in the above we assume that there are only two rays $d_-(t), d'_-(t)$ colliding at $t_0$. The extension to the case of collision of many anti-Stokes rays at a point $t_0$ is direct. 

In a summary, as $t$ varies from $0$ to $1$, the effective anti-Stokes rays $-{\rm arg}\big(h(u_i-u_j)\big)$ are preserved. In the meanwhile, let us assume some rays are ineffective anti-Stokes rays as $t$ near $0$. As $t$ moves toward $1$ along the real axis, these rays will move on the complex $z$-plane. They can cross each other and cross the ($t$-independent) effective rays as $t$ crosses through some particular point. And after crossing such a particular point, the ineffective rays are still ineffective. It therefore finishes the proof.
\end{proof} 

\vspace{2mm}
Proposition \ref{extend0}, or the proof of Proposition \ref{Strivial}, leads to a direct corollary.
\begin{cor}\label{Stokesconti}
The solution $Y_{i2}(z,t)$ of \eqref{gcKZ1} and \eqref{gcKZ2} defined on ${\Sect_{i}\times (0,1)}$, as the analytic continuation of the solution $Y_{d,d'}(z,t)$ given in \eqref{yik}, has the following asymptotics for all $t\in (0,1)$
\[Y_{i2}(z,t)\cdot e^{-hz(u^{(1)}+tu^{(2)})} z^{-h{\delta T^{(1)}}}(tz)^{-h{\delta T^{(2)}}}(tz-z)^{-h{\delta P}}\rightarrow 1,\]
as $z\rightarrow \infty \text{ within } \Sect_i$. 
\end{cor}

\subsection{Step 5: Proof of the second identity in part $(b)$ of Theorem \ref{fact}}

Let us consider the function $Y_{i2}(z,t)$ on $(z,t)\in\Sect_i\times (0,1)$. The second identity in part $(b)$ of Theorem \ref{fact} is equivalent to
\begin{pro}\label{prob2}
For fixed $zt-z\in \Sect_{l+i}$, the function $\overline{W}_i(z,zt-z):=Y_{i2}(z,t)X_{l+i}^{(2)}(zt-z)^{-1}$ is the unique solution of \eqref{ratio3eq} with the prescribed asymptotics in \eqref{simW3}.
\end{pro}
\begin{proof} 
Let us take the coordinates transform $\omega_1=z$ and $\omega_2=zt-z$. The anti-Stokes rays of \eqref{ratio3eq} take the form of $-{\rm arg}\big(h(u_i-u_j+u_k-u_l)\big)$. Let us take any Stokes supersector $\Sect_i'$ of \eqref{ratio3eq} which is contained inside $\Sect_{i}$ (the sector $S(d_i-{\pi}/{2},d_{i+1}+{\pi}/{2})$ of the $\omega_1$ plane), and denote by $\overline{W'}_i(\omega_1=z,\omega_2=zt-z)$ the corresponding solution of \eqref{ratio3eq} with prescribed asymptotics \eqref{simW3} on $\Sect_i'$. Then based on the equation satisfied by $\overline{W'}_i$ and $X_{l+i}^{(2)}$, we check that the function
\[\overline{W'}_i(z,zt-z) X_{l+i}^{(2)}(zt-z),\]
defined for $z\in \Sect_i'\subset \Sect_i$ and $t\in (0,1)$, satisfies the system of equations \eqref{gcKZ1} and \eqref{gcKZ2}. Similar to the proof of Proposition \ref{extend0}, a detailed analysis of the Borel-Laplace transform shows that there exists a sufficiently small real number $\varepsilon>0$ such that 
\[ \overline{W'}_i(z,zt-z)\cdot e^{-hz(u^{(1)}+u^{(2)})}z^{-h\delta T^{(1)}-h\delta T^{(2)}}\sim 1\] as $z\rightarrow\infty$ within $\Sect'_i$ for any fixed $t\in (1-\varepsilon,1)$.

Then Corollary \ref{Stokesconti} and the same argument as in Proposition \ref{prob} lead to the consequence that the solutions $Y_{i2}(z,t)$ and $\overline{W'}_i(z,zt-z) X_{l+i}^{(2)}(zt-z)$ have the same asymptotics, and actually
\[Y_{i2}(z,t)=\overline{W'}_i(z,zt-z) X_{l+i}^{(2)}(zt-z) \ \ \text{ for } (z,t)\in \Sect'_i\times (1-\varepsilon,1). \]
Since $Y_{i2}(z,t)$ is defined for $z\in \Sect_i$ and $t\in (0,1)$. It particularly implies that for any two Stokes supersector $\Sect_i'$ and $\Sect''_i$ of \eqref{ratio3eq} that are contained inside the same $\Sect_{i}$, the corresponding solutions $\overline{W'}_i$ and $\overline{W''}_i$ of \eqref{ratio3eq}, defined on $\Sect_i'$ and $\Sect''_i$, coincide in their common domain $\Sect_i'\cap\Sect''_i$, and therefore glue together to a unique function, denoted by $\overline{W}_i$, with the prescribed asymptotics in \eqref{simW3}. (It shows that the Stokes factors associated to the anti-Stokes rays $-{\rm arg}\big(h(u_i-u_j+u_k-u_l)\big)$ of  \eqref{ratio3eq}, that do not coincide with $d_i$ for $i=0,1,...,2l-1$, are the identity.) This concludes the proof.
\end{proof}

\begin{rmk}
From the above proof, we see that there exist sufficiently small positive real numbers $\varepsilon_1,\varepsilon_2,\varepsilon_3>0$ such that \begin{itemize}
    \item the $z\rightarrow\infty$ asymptotics of $W_i(z,zt)$ as $t$ fixed is valid for $t\in (-\varepsilon_1,0)\cup (0,\varepsilon_1)$, as in Proposition \ref{extend0};
    \item the $z\rightarrow \infty$ asymptotics of $\overline{W}_i(z,zt-z)$ is valid for fixed $t\in (1-\varepsilon_2,1)\cup (1,1+\varepsilon_2)$, as in Proposition \ref{prob2};
    \item the $z\rightarrow$ asymptotics of $\widetilde{W}_i(z,zt)$ as $t$ fixed is valid for $t\in (-\infty,-1/\varepsilon_3)\cup (1/\varepsilon_3,+\infty)$.
\end{itemize}
Since the interval $(-\infty,-1/\varepsilon_3)$ may not overlap with $(-\varepsilon_1,0)$, we can not directly state that the left hand side $\widetilde{W}_{i}(z,zt)F_{i+l}^{(1)}(z)$ of the identity \eqref{parta}, in part $(a)$ of Theorem \ref{fact}, is equal to the right hand side $W_i(z,zt)F_{i+l}^{(2)}(zt)e^{\pi\mathi  h\delta P}$ (after the necessary analytic continuation) by comparing their asymptotics. So it is necessary to introduce the solution $Y_{i1}(z,t)$ connecting them, similarly, the solutions $Y_{i2}(z,t)$ and $Y_{i3}(z,t)$.
\end{rmk}

\subsection*{Acknowledgements}
\noindent We would like to thank Anton Alekseev for useful comments. The author is supported by the National Key Research and Development Program of China (No. 2021YFA1002000) and by the National Natural Science Foundation of China (No. 12171006).

\Addresses
\end{document}